\documentclass[11pt]{article} %%[R1]

\usepackage[utf8]{inputenc}

\usepackage{float}
\usepackage{graphicx}
\usepackage[ruled,vlined]{algorithm2e}
\usepackage{algorithmic}
\usepackage{epsfig}
\usepackage{epstopdf}
\usepackage{amssymb,amsmath}
\usepackage{amsmath}
\usepackage{amsthm}
\usepackage{amsfonts}
\usepackage{psfrag}
\usepackage{rotating}
\usepackage{latexsym}
\usepackage{dsfont}
\usepackage{stmaryrd}
\usepackage{amsthm}
\usepackage{amssymb}
\usepackage{amsmath}
\usepackage{mathtools}
\usepackage{mathrsfs}   % Weitere Mathematik-Symbole, z.B. Laplace-L
\usepackage{siunitx}
\usepackage{pgfpages}
\usepackage[Euler]{upgreek} %Geradestehende griechische Buchstaben
\usepackage{isomath}
\usepackage{physics}
\usepackage{pifont} %For ding-Symbols
\usepackage{bbm}
\usepackage{empheq}
\usepackage[thicklines]{cancel}
\usepackage{subcaption}
\usepackage[authoryear,semicolon]{natbib}
\usepackage[force,almostfull]{textcomp}
\usepackage{fancyvrb}
\usepackage{textcomp}
\usepackage{url}
\usepackage{ifdraft}
\usepackage{url}
\usepackage{multirow}
\usepackage{rotating}
\usepackage{xspace}
\usepackage{array}
\usepackage{multido}
\usepackage{enumerate}
\usepackage[utf8]{inputenc}
\newlength\figureheight 
\newlength\figurewidth 
\usepackage{graphics}
\usepackage{pgfplots}
\pgfplotsset{compat=newest}
\pgfplotsset{plot coordinates/math parser=false}
\allowdisplaybreaks

\usepackage{scalefnt}

\newtheoremstyle{specialcasestyle}{1mm}{1mm}{\upshape}{}{\bfseries\upshape}{.}{0mm}{}
\theoremstyle{specialcasestyle}

\newtheorem{rem}{Remark}

\newtheorem{constraint}{Constraints}
\newtheorem{prblm}{Problem}

\theoremstyle{remark}          % Avoid cursive text in example environment

\newcommand*{\E}[1]{\ensuremath{\mathbb{E}\left[#1\right]}}
\newcommand*{\prob}[1]{\ensuremath{\mathbb{P}\left[#1\right]}}
\newcommand*{\tol}{\ensuremath{\mathrm{TOL}}}

\newcommand{\argmin}{\operatorname*{argmin}}
\newcommand{\argmax}{\operatorname*{argmax}}

\begin{document} %% [R2]

\title{Optimal power procurement for green cellular wireless networks under uncertainty and chance constraints}
\author{Nadhir Ben Rached \thanks{Department of Statistics, School of Mathematics, University of Leeds ({\tt N.BenRached@leeds.ac.uk}).}, Shyam Mohan Subbiah Pillai \thanks{Corresponding author; Chair of Mathematics for Uncertainty Quantification, Department of Mathematics, RWTH Aachen University, Aachen, Germany({\tt subbiah@uq.rwth-aachen.de}).} \\
and Ra\'ul Tempone \thanks{Computer, Electrical and Mathematical Sciences \& Engineering Division (CEMSE), King Abdullah University of Science and Technology (KAUST), Thuwal, Saudi Arabia ({\tt raul.tempone@kaust.edu.sa}). Alexander von Humboldt Professor in Mathematics for Uncertainty Quantification, RWTH Aachen University, Aachen, Germany ({\tt tempone@uq.rwth-aachen.de}).}.
        \thanks{This work was supported by the KAUST Office of Sponsored Research (OSR) under Award No. URF/1/2584-01-01 and the Alexander von Humboldt Foundation. This work was also partially performed as part of the Helmholtz School for Data Science in Life, Earth and Energy (HDS-LEE) and received funding from the Helmholtz Association of German Research Centres. For the purpose of open access, the author has applied a Creative Commons Attribution (CC BY) licence to any Author Accepted Manuscript version arising from this submission.}
}
\date{}
\maketitle
\thispagestyle{empty}
\begin{abstract}
Given the increasing global emphasis on sustainable energy usage and the rising energy demands of cellular wireless networks, this work seeks an optimal short-term, continuous-time power procurement schedule to minimize operating expenditure and the carbon footprint of cellular wireless networks equipped with energy storage capacity, and hybrid energy systems comprising uncertain renewable energy sources. Despite the stochastic nature of wireless fading channels, the network operator must ensure a certain quality-of-service (QoS) constraint with high probability. This probabilistic constraint prevents using the dynamic programming principle to solve the stochastic optimal control problem. This work introduces a novel time-continuous Lagrangian relaxation approach tailored for real-time, near-optimal energy procurement in cellular networks, overcoming tractability problems associated with the probabilistic QoS constraint. The numerical solution procedure includes an efficient upwind finite-difference solver for the Hamilton--Jacobi--Bellman equation corresponding to the relaxed problem, and an effective combination of the limited memory bundle method (LMBM) for handling nonsmooth optimization and the stochastic subgradient method (SSM) to navigate the stochasticity of the dual problem. Numerical results, based on the German power system and daily cellular traffic data, demonstrate the computational efficiency of the proposed numerical approach, providing a near-optimal policy in a practical timeframe.
\end{abstract}
\textbf{Keywords:} Stochastic optimal control; chance constraints; Lagrangian relaxation; dynamic programming; wireless networks. \\
%\textbf{2010 Mathematics Subject Classification} 60H35. 65C30. 65C05. 65C35.

\section{Introduction}
\label{sec:intro}

Since 2021, worldwide mobile broadband traffic has increased by an average of 19.6\% annually, reaching 1 zettabyte (ZB) in 2023 \citep{UTI2024}. 4G network coverage increased from 41\% of the world population in 2015 to 92\% of the world population in 2024 \citep{UTI2024}. The estimated transmission network electricity consumption in 2022 was 240 to 340~TWh, an increase of 64\% from 2015 \citep{IEA2023}. Next-generation transmission technology (e.g., 5G) requires even more power to ensure sufficient mobile network coverage \citep{strielkowski20215g}. A consequence of the high energy demand is increased greenhouse gas emissions. Data transmission networks accounted for around 330~Mt $\text{CO}_2$-equivalent in 2020, around 0.9\% of all energy-related emissions \citep{IEA2023}. About 83\% of the energy consumed by cellular networks came from fossil fuels in 2021, whereas only 9\% came from renewables \citep{GSMA2023}. This problem motivates the coupling of cellular wireless networks with renewable energy to reduce the costs and carbon footprint of the telecommunications industry. Hybrid energy systems provide a consistent power supply to wireless networks by combining renewable energy sources (e.g., solar and wind) \citep{han2014powering}. However, modeling these systems and deriving an optimal power procurement policy can be challenging, due to the stochasticity of solar/wind power and wireless channels \citep{electronics12092014}. 

Previous studies have considered the optimal power procurement problem in a discrete-time setting for isolated energy systems under an uncertain energy demand and incoming renewable energy~\citep{Kloppel:2013aa,Wu:2014aa,Huang:2018aa}. The corresponding optimization problem for wireless networks has additional difficulties due to the uncertainty in wireless channels and probabilistic constraints due to the network quality of service (QoS) \citep{Azgin:2003aa,farooq2016_stochastic,Classen:2015aa,Challita:2016aa}. The solution method for this problem in the current literature can be generally categorized into robust optimization and stochastic optimization. The robust optimization approach constructs a bounded uncertainty set encompassing all realizations of the stochastic variables and optimizes the worst-case objective function value \citep{Ding:2018aa,Chang:2021aa,Zhou:2023aa,Du:2023aa}. Robust optimization, while ensuring feasibility under the worst-case scenario, often leads to overly conservative policies that could underutilize renewable energy resources. The stochastic optimization approach assumes a known probability distribution for the random variables and optimizes the expected value of the objective function over the probability distribution of the uncertain parameters \citep{An:2017aa,Bhattacharya:2018aa}. 

Chance-constrained programming is a subclass of stochastic optimization methods, first introduced in~\citep{Charnes:1958aa,Charnes:1959aa}. This approach  allows constraints to be violated with a certain prescribed probability or risk level. The standard approach to solving discrete chance-constrained problems is to replace the chance constraint with a conservative but tractable approximation using the Bernstein inequality~\citep{Nemirovski:2006aa,Ma:2013aa}. This approach is well-suited for wireless networks, since the associated QoS constraint is commonly defined probabilistically, allowing the operator to balance between conservatism and optimality. However, the computational complexity of the method in the discrete-time setting increases exponentially for large datasets and complex problems. For example, the authors in~\citep{Liu:2013aa,benrached17_energy,Liu:2022aa} solved the resulting discrete optimization problem using heuristic techniques, yielding suboptimal solutions.

This work proposes a novel mathematical modeling and numerical framework for the optimal management of hybrid energy systems to power cellular wireless networks under uncertainty and chance constraints in continuous-time for a short-term planning horizon, where the system comprises locally installed renewable energy sources, fossil-fuel power stations, and a storage capacity modeled by a single battery. We employ stochastic differential equations (SDEs) to model the incoming instantaneous renewable power and wireless fading channels. Data-driven SDEs have been employed to model instantaneous wind~\citep{Moller:2016aa,Elkantassi:2017aa,Iversen:2017aa,caballero21_derivative} and solar power~\citep{Iversen:2014aa,Iversen:2017ab,Badosa:2018aa}. These SDE models follow a mean-reversion to deterministic wind speed/solar irradiance forecasts. The parameters of the SDE characterizing its uncertainty are calibrated from historic discrepancies (data) between the actual and forecasted quantities. Moreover, we employ SDEs to model stationary wireless fading channels. The authors in~\citep{Feng:2007aa,Charalambous:2008aa,Mossberg:2009aa,Amar:2024aa} have constructed and configured SDEs to reproduce characteristics of various fading channel models, primarily concerning their stationary distribution. 

This work formulates and solves a time-continuous stochastic optimal control problem to derive a near-optimal policy minimizing the  expected operating expenditure and carbon footprint of cellular base stations. We optimize for both cost minimization and environmental impact, introducing a trade-off parameter $w \in [0,1]$ to balance these objectives. The optimal control is subject to demand constraints, a probabilistic and time-pointwise QoS constraint, and the dynamics and capacity limits of energy systems and wireless channels. The probabilistic constraint prevents the direct application of the dynamic programming principle because the constraint controls the joint distribution of the random variables. The optimal policy must be a function of the current state and time for the standard dynamic programming lemma to hold, not a function of its distribution~\citep{Bertsekas:2012aa}. To overcome this problem, we pose the corresponding dual problem by introducing a time-continuous Lagrangian relaxation to penalize violations of the probabilistic constraint~\citep{Lemarechal:2001aa}. This dual formulation yields a standard time-continuous stochastic optimization problem solved iteratively by optimizing the dual function. Each iteration involves solving the Hamilton--Jacobi--Bellman (HJB) partial differential equation (PDE)~\citep{Kushner:1990aa,Pham:2009aa} using an upwind finite-difference scheme~\citep{Sun:2015aa} to compute the dual objective value and subgradient. With access to only noisy subgradients, we solve the convex, nonsmooth dual problem using a novel combination of the limited memory bundle method (LMBM)~\citep{Lemarechal:1977aa,Makela:2002aa,Karmitsa:2012aa} and the stochastic subgradient method~\citep{Boyd:2003aa,Boyd:2008aa}. Moreover, we design a refinement strategy for the Lagrangian multipliers to control the pointwise violation in the chance constraint. The numerical experiments and results based on the German power system and daily cellular traffic profiles validate the efficiency of the proposed models and approach.

The contributions can be summarized as follows:

\begin{enumerate}
	\item We propose a novel time-continuous optimization framework for optimal power procurement for green wireless cellular networks, subject to SDE dynamics and chance constraints. Compared to the discrete-time formulation in previous studies~\citep{Liu:2013aa,benrached17_energy,Liu:2022aa}, the proposed approach decouples the model development from the numerical approximation, enhancing the model fidelity (see Remark~\ref{rem:disc_vs_cont}). This formulation also yields a continuous control curve over time, allowing its application for any time discretization scheme and eliminating the need for ad hoc interpolations~\citep{Ben-Hammouda:2024aa}.
	\item We calibrate the data-driven SDE model developed for instantaneous wind power in~\citep{caballero21_derivative} using German wind power data from the year 2023. The calibrated SDE is a driving dynamic for the stochastic optimal control problem.
	\item We derive an SDE to model the instantaneous Nakagami wireless fading channel~\citep{nakagami19603} with a shifted-gamma invariant distribution, using the Pearson class of diffusions~\citep{forman2008_pearson}. This SDE is a driving dynamic for the stochastic optimal control problem.
	\item We apply Lagrangian relaxation to the probabilistic QoS constraint, transforming the problem into a standard time-continuous stochastic optimization problem. Studies have explored numerical methods for time-continuous stochastic optimization with final-time chance constraints using Lagrangian relaxation~\citep{Andrieu:2011aa,pfeiffer:pastel-00881119,Pfeiffer:2021aa} or reformulation as a stochastic target problem~\citep{Bouchard:2010aa,Bouchard:2009aa}. However, the proposed approach is novel in addressing a chance constraint that must be satisfied at every time point. Moreover, we implement this within the context of cellular wireless networks.
	\item We develop an iterative algorithm to optimize the dual function within a finite-dimensional function class numerically. Each iteration involves solving the HJB PDE to compute the dual function value and its noisy subgradient. The proposed approach extends the work in~\citep{Ben-Hammouda:2024aa} on the time-continuous deterministic optimization of coupled hydrothermal power systems to the stochastic setting.
	\item With access only to noisy realizations of the subgradient, we combine the LMBM method~\citep{Haarala:2004aa} with the stochastic subgradient method~\citep{Boyd:2008aa} to solve the dual optimization problem. We iteratively refine the Lagrange multiplier to ensure compliance with the chance constraint at every time point.
\end{enumerate}

The outline of the paper is as follows. Section~\ref{sec:model} introduces the complete system model, detailing the costs, power generation, wireless network, and state dynamics. Next, Section~\ref{sec:soc_formulation} formulates the primal stochastic optimal control problem and the associated dual problem by Lagrangian relaxation of the chance constraint. Then, Section~\ref{sec:numerical_methods} introduces the numerical approach for solving the dual problem and presents all components of the final algorithm in detail. This section discusses the numerical solution of the HJB PDE, subgradient estimation, dual-problem optimization, and numerical errors arising from the proposed approach. Finally, Section~\ref{sec:results} demonstrates the numerical results of applying the above approach when solving the stochastic optimization problem for a cellular base station based on the German power system and daily cellular traffic profiles.

\section{System Model}
\label{sec:model}

\begin{figure}[h!]
	\centering
	\includegraphics[width=0.8\textwidth]{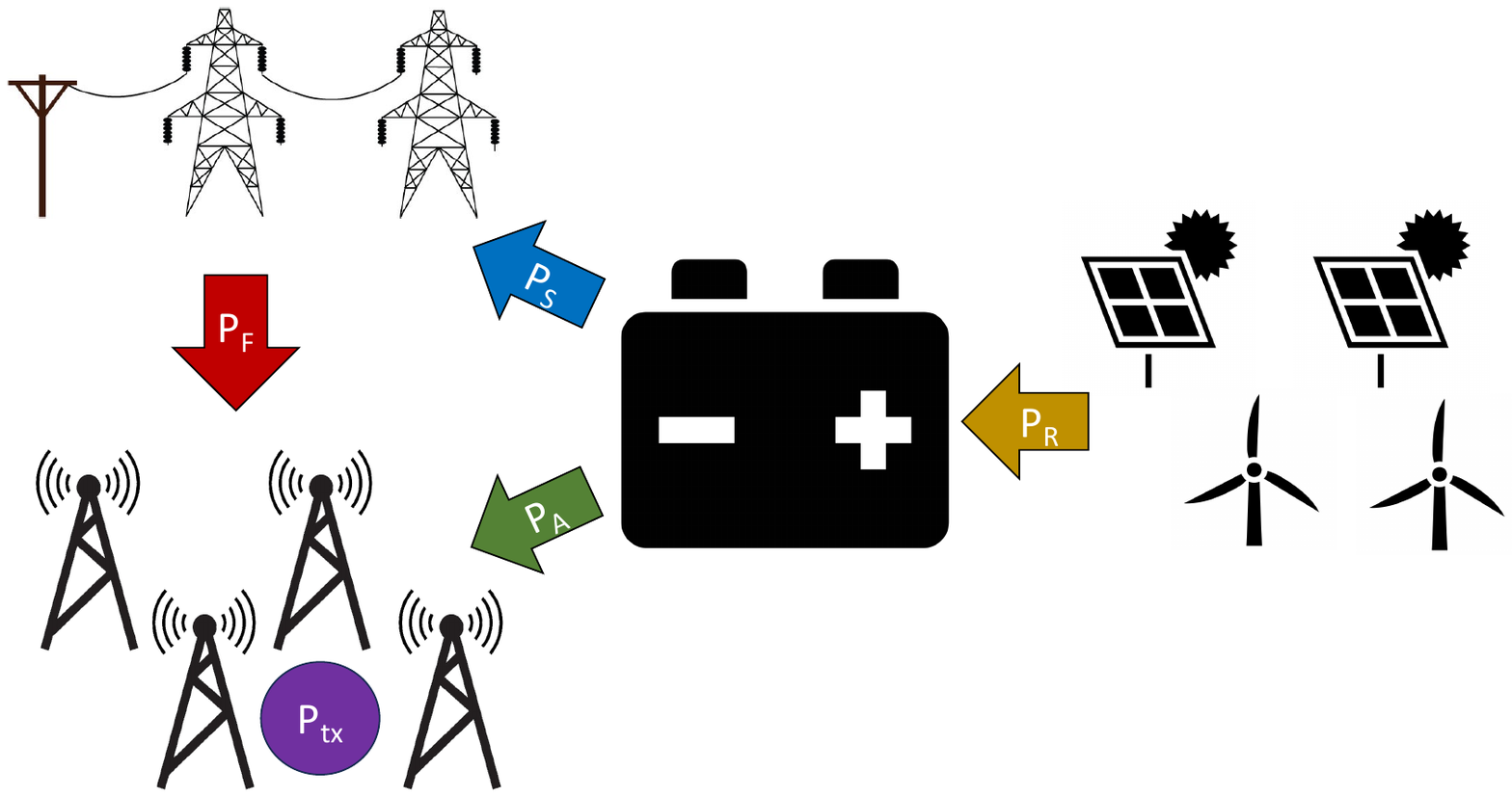}
	\caption{Schematic illustration of the power flow in a base station (Section~\ref{sec:base_station_model}) in a cellular wireless network.}
	\label{fig:base_station}
\end{figure}

This work considers a cellular network of multiple noninteracting base stations, each powered by renewable energy sources (solar panels or wind turbines) and battery storage. Figure 1 in~\citep{benrached17_energy} depicts the schematic network. Figure~\ref{fig:base_station} illustrates the electrical power flow in each base station, each with a battery that stores the power generated from renewables ($P_R$). The power transmitted by the base station to serve cellular users consists of power procured from the battery ($P_A$) and bought from the grid ($P_F$). Any extra energy stored in the battery could be sold back to the grid  ($P_S$) for revenue. We aim to optimize the energy procurement at each base station for a one-day operation cycle.

\subsection{Base Station Model}
\label{sec:base_station_model}

The instantaneous power procured by the mobile operator from the traditional grid to power the base station is denoted by $\{P_F(t) \in \mathbb{R}^+: t \in [0,T]\}$, where $0$ and $T$ denote the time at the beginning and end of the day, respectively. The corresponding unit price is denoted by $\{K_b(t) \in \mathbb{R}^+: t \in [0,T]\}$. This price can vary during the day depending on strategies followed by the stakeholders in the energy market. Moreover, the base station has its own internal source of power (i.e. renewable power generator) whose generated power is denoted by $\{P_R(t) \in \mathbb{R}^+: t \in [0,T]\}$. This power is assumed to be free of charge. The internal power source is assumed to have a maximum power output capacity denoted by $\bar{P}_R \in \mathbb{R}^+$. The base station has a battery, allowing the mobile operator to store the incoming renewable power. As the base station is interconnected with the traditional grid, the mobile operator chooses to use the stored renewable power to run the base station or sell the extra energy back to the grid for revenue. The instantaneous power drawn from the battery to run the base station is denoted by $\{P_A(t) \in \mathbb{R}^+: t \in [0,T]\}$. The instantaneous power from the battery sold back to the grid is denoted by $\{P_S(t) \in \mathbb{R}^+: t \in [0,T]\}$, and the corresponding unit price is denoted by $\{K_s(t) \in \mathbb{R}^+: t \in [0,T]\}$. We reasonably assume $K_s(t) \leq K_b(t)$ for all $t \in [0,T]$. The instantaneous power transmitted by the base station to serve its users is denoted by $\{P_\mathrm{tx}^\mathrm{tot}(t)\in \mathbb{R}^+: t \in [0,T]\}$. We assume the total transmitted power is equally divided among each average mobile user connected to the network. The instantaneous power transmitted by the base station per average user is denoted by $\{P_\mathrm{tx}(t)\in \mathbb{R}^+: t \in [0,T]\}$ and $P_\mathrm{tx}^\mathrm{tot}(t) = N_u(t) P_\mathrm{tx}(t)$ for all $t \in [0,T]$, where $\{N_u(t) \in \mathbb{R}^+: t \in [0,T]\}$ is the number of connected mobile users at time $t$. The instantaneous power balance equation for the base station is

\begin{equation}
	\label{eqn:power_balance_constraint}
	P_A(t) + P_F(t) = C_\text{scal} N_u(t) P_\mathrm{tx}(t) + C_\text{offset}, \quad 0 \leq t \leq T ,
\end{equation} 

where $C_\text{scal} \in \mathbb{R}^+$ represents a factor that scales with the transmitted power due to the amplifier and feeder losses in the network, and $C_\text{offset} \in \mathbb{R}^+$ models the offset power required to operate the base station irrespective of the transmitted power~\citep{benrached17_energy}. We assume that the base station can transmit a maximum power of $\bar{P}_\mathrm{tx} \in \mathbb{R}^+$.

\begin{equation}
	P_\mathrm{tx}^\mathrm{tot}(t) \leq \bar{P}_\mathrm{tx}, \quad 0 \leq t \leq T \cdot
\end{equation}

The term $N_u(t)$ in~\eqref{eqn:power_balance_constraint} is determined by the choice of the daily mobile user traffic model, heavily influencing the power demand in the base station. One such model is introduced in Section~\ref{sec:model_description}.

\subsection{Cellular Network Model}

According to the power law, the power transmitted by the base station decays with increased distance from the base station. Moreover, the transmitted power is uncertain due to atmospheric conditions, natural and human-made environmental obstacles, and interference from other signals~\citep{farooq2016_stochastic}. Hence, the instantaneous power received by a user $u$ at position $\mathbf{x}_u(t) = (x_u(t),y_u(t)) \in \mathbb{R}^2$ from a base station at $\mathbf{x}_\text{BS} = (x_\text{BS},y_\text{BS}) \in \mathbb{R}^2$ is given as $P_\mathrm{rx}^u(t) = P_\mathrm{tx}(t) \xi(t) \kappa \norm{\mathbf{x}_u(t)-\mathbf{x}_\text{BS}}^{-\eta}$, where $\kappa, \eta \in \mathbb{R}^+$ denotes the path loss constant and exponent, respectively. Moreover, $\norm{\cdot}$ denotes the Euclidean norm. $\xi: [0,T] \cross \Omega \rightarrow [\underline{\xi},\infty)$ represents a stochastic process, that is greater than $\underline{\xi} \in \mathbb{R}^+$ almost surely (a.s.), modeling the wireless fading channel. 

This work assumes the channel follows Nakagami fading~\citep{nakagami19603}. The Nakagami fading model is a good data-fit model for cellular systems in urban and suburban areas~\citep{nakagami_datafit}. In this model, the instantaneous power is gamma distributed, a  Pearson type III distribution~\citep{Pearson:1997aa}. Hence, Pearson diffusion~\citep{forman2008_pearson} is employed to model the wireless fading channel. Thus, this work constructs an SDE for the process $\xi$ \eqref{eqn:gamma_fading} with a mean-reverting linear drift and a squared diffusion, which is a second-order polynomial of the state, and whose invariant  distribution is the shifted-gamma distribution~\citep{forman2008_pearson} with shape parameter $\mu>1$ and scale parameter $\theta \in \mathbb{R}^+$. We consider the probability space $\{\Omega,\mathcal{F},\{\mathcal{F}_t\}_{t \geq 0}, P\}$, where $\mathcal{F}_t$ is the filtration of the standard Wiener process $W:[0,T] \cross \Omega \rightarrow \mathbb{R}$:

\begin{empheq}[left=\empheqlbrace, right = \cdot]{equation}
\label{eqn:gamma_fading}
	\begin{alignedat}{2}
		\dd \xi(t) &= -\theta \left( \xi(t) - \underline{\xi} - \mu \right) \dd t + \sqrt{2 \theta (\xi(t)-\underline{\xi})} \dd W(t), \quad t>0 \\
		\xi(0) &\sim \mu^\xi_0,
	\end{alignedat}
\end{empheq}

where $\mu^\xi_0$ is the invariant distribution of the process $\xi$.

\begin{rem}{\textit{ (Computational efficiency of Pearson diffusions). }}
	Pearson diffusion naturally models the continuous-time wireless fading channel, where process $\xi$ lies in $\left[ \underline{\xi},\infty \right)$ for all $t \in [0,T]$. This model also has the advantage of analytical tractability, providing closed-form solutions for the invariant distribution of the process $\xi$ and offering better insight into its behavior. This approach makes the model choice computationally more efficient than other Markov-based models (e.g. discrete-state Markov chains), whose complexity increases quadratically as the number of possible states increases.
\end{rem}

In this work, the network performance is characterized by the offered QoS, determined by the signal-to-noise ratio (SNR) threshold $\text{SNR}_\mathrm{th}$. The SNR, in decibels (dB), at user $u$ at $\mathbf{x}_u(t)$ served by the base station can be expressed as follows:

\begin{equation}
	\label{eqn:snr_definition}
		\text{SNR}^u(t) = 10 \log_{10} \left( \frac{P_\mathrm{rx}^u(t)}{\sigma_0} \right) = 10 \log_{10} \left( \frac{P_\mathrm{tx}(t) \xi(t) \kappa \norm{\mathbf{x}_u(t)-\mathbf{x}_\text{BS}}^{-\eta}}{\sigma_0} \right) \cdot
\end{equation}

The ambient noise level is denoted by $\sigma_0 \in \mathbb{R}^+$. The QoS of the base station is determined by the proportion of users whose SNR is above $\text{SNR}_\mathrm{th}$. A user $u$ is defined to be in outage at time $t$ if $\text{SNR}^u(t) < \text{SNR}_\mathrm{th}$. The proportion of users in outage is denoted by $\{\phi_\text{out}(t) \in [0,1]: t \in [0,T]\}$ and can be expressed as follows:

\begin{equation}
	\label{eqn:outage_proportion_defn}
	\phi_\text{out}(t) = \frac{1}{N_u(t)} \sum_{u=1}^{N_u(t)} \mathbbm{1}_{\{\text{SNR}^u(t) < \text{SNR}_\mathrm{th}\}},
\end{equation}

where $\phi_\text{out}(t)$ can be interpreted as an empirical average over a user distribution. The following approximation assumes that each user is identical and that sufficient number of mobile users exist at all times:

\begin{align}
	\phi_\text{out}(t) &= \frac{1}{N_u(t)} \sum_{u=1}^{N_u(t)} \mathbbm{1}_{\{\text{SNR}^u(t) < \text{SNR}_\mathrm{th}\}} \nonumber \\
	&= \frac{1}{N_u(t)}\sum_{u=1}^{N_u(t)} \mathbbm{1}_{\left\{10 \log_{10} \left( \frac{P_{\mathrm{tx}}(t) \xi(t) \kappa \norm{\mathbf{x}_u(t)-\mathbf{x}_\text{BS}}^{-\eta}}{\sigma_0} \right) < \mathrm{SNR}_{\mathrm{th}}\right\}} \nonumber \\
	\label{eqn:mean_field_approx}
	&\approx \int_{\mathbb{R}^2} \mathbbm{1}_{\left\{ \norm{z - \mathbf{x}_{\mathrm{BS}}} > \left( \frac{P_{\mathrm{tx}}(t) \xi(t) \kappa}{\sigma_0 10^{\frac{\mathrm{SNR}_{\mathrm{th}}}{10}}} \right)^{\frac{1}{\eta}} \right\}} \rho_z(t) \dd z,
\end{align}

where $\rho_z(t)$ denotes the distribution of mobile users at time $t$. In \eqref{eqn:mean_field_approx}, $\phi_\text{out}(t)$ is the complementary cumulative distribution function (ccdf) of the user distribution. This quantity can be computed (even analytically) for many well-known distributions. 

The objective of the network operator is to ensure $\phi_\text{out}(t) < \phi_\mathrm{th}$ for all $t \in [0,T]$, where $0 < \phi_\mathrm{th} \ll 1$ represents a threshold ratio. Due to the uncertainty associated with channel fading, achieving the objective in an a.s. sense is infeasible. Instead, we introduce the following chance constraint for the network QoS.

\begin{equation}
	\label{eqn:qos_constraint}
	\prob{\phi_\text{out}(t) \geq \phi_\mathrm{th}} \leq \epsilon, \quad 0 \leq t \leq T,
\end{equation}

where $1-\epsilon$, for $0 < \epsilon \ll 1$, is the confidence level of satisfying the constraint in \eqref{eqn:qos_constraint}. The network operator also generates revenue from servicing the users throughout the day: 

\begin{equation}
	\label{eqn:revenue_network}
	\mathcal{R}_1 = \E{\int_0^T K_\text{net}(t) N_u(t) (1- \phi_\text{out}(t)) \dd t},
\end{equation}

where $\{K_\text{net}(t) \in \mathbb{R}^+: t \in [0,T]\}$ denotes the unit price paid by users to connect to the network, and $N_u(t) (1- \phi_\text{out}(t))$ represents the number of users that are connected to the network (not in outage).

\subsection{Renewable Power Model}

This work also incorporates the uncertainty in renewable power generation. We apply a data-driven parametric SDE, whose solution defines a stochastic process that models the error in a given forecast for renewable power~\citep{caballero21_derivative}. Let the process $R:[0,T] \cross \Omega \rightarrow [0,1]$ denote the normalized generated renewable power.

\begin{empheq}[left=\empheqlbrace, right = \cdot]{equation}
	\label{eqn:re_dynamics}
	\begin{alignedat}{2}
		\dd R (t) &= \left( \dot{p}(t) - \theta(t) (R(t) - p(t)) \right) \dd t + \sqrt{2 \alpha \theta_0 R(t) (1 - R(t))} \dd W(t), \quad t>0  \\ 
		R(0) &\sim \mu^R_0,
	\end{alignedat}
\end{empheq}

where $\{p(t) \in [0,1]:t \in [0,T]\}$ is a deterministic forecast for the normalized power provided by an official source and $\{\dot{p}(t):t \in [0,T]\}$ is its time derivative. In addition, $\{\theta(t): t \in [0,T]\}$ denotes the mean-reversion time-varying parameter ensuring that the process $R$ is the unique strong solution to~\eqref{eqn:re_dynamics} for all $t \in [0,T]$ with range $[0,1]$, a.s. Moreover, $\alpha,\theta_0 \in \mathbb{R}^+$ are parameters inferred from discrepancies between historical real-world  observations and their forecasts. In addition, $W:[0,T] \cross \Omega \rightarrow \mathbb{R}$ is a standard Wiener process. The drift function in \eqref{eqn:re_dynamics} is designed such that the process reverts to its mean $p(t)$, with a time-varying speed $\theta(t)$  proportional to the deviation of the process $R(t)$ from its mean, and it tracks the time derivative $\dot{p}(t)$. The diffusion function in \eqref{eqn:re_dynamics} is designed so that the process $R(t)$ avoids exiting from the range $[0,1]$. The net renewable power generated is then given by $P_R(t) = \bar{P}_R R(t)$ for all $t \in [0,T]$.

\subsection{Battery Model}

This work uses a simple model for the battery based on a circuit model of the state of the charge of a lithium-ion battery~\citep{brucker2021_grey}. We assume that it has no energy loss when used and that no associated operational costs or aging effects exist. Process $A:[0,T] \rightarrow [0,1]$ denotes the normalized charge in the battery:

\begin{empheq}[left=\empheqlbrace, right = \cdot]{equation}
	\label{eqn:battery_dynamics}
	\begin{alignedat}{2}
		\dd A(t) &= \frac{\left( P_R(t) - P_S(t) - P_A(t) \right)}{\bar{A}} \dd t \\
		A(0) &= A_0 \\
		-\underline{P}_A(A(t)) &\leq P_A + P_S - P_R \leq \bar{P}_A(A(t)),
	\end{alignedat}
\end{empheq}

where $\bar{A}$ represents the maximum battery charge capacity, $\underline{P}_A$ denotes the maximum power the battery can absorb, and $\bar{P}_A$ indicates the maximum power that the battery can supply. These quantities are characteristic of the battery. Generally, the battery can supply more power than it can absorb; hence, $\underline{P}_A < \bar{P}_A$ for all $t \in [0,T]$. In addition, $A_0 \in [0,1]$ represents the normalized charge in the battery at the start of the day. The requirement that $0 \leq A(t) \leq 1$ for all $t \in [0,T]$ is enforced via constraints on charging and discharging at the extremes:

\begin{equation}
\label{eqn:characteristic}
	\begin{cases}
		P_R(t) - P_S(t) - P_A(t) &\geq 0, \quad \text{if } A(t) = 0 \\
		P_R(t) - P_S(t) - P_A(t) &\leq 0, \quad \text{if } A(t) = 1  \cdot
	\end{cases}
\end{equation}

The network operator generates revenue by selling some energy stored in the battery throughout the day:

\begin{equation}
	\label{eqn:revenue_power}
	\mathcal{R}_2 = \E{\int_0^T K_s(t) P_S(t) \dd t} \cdot
\end{equation}

\begin{rem}{\textit{ (Battery model improvements)}.}
	 Although the battery model~\eqref{eqn:battery_dynamics} is straightforward and idealistic, realistic battery efficiency losses in the model can be incorporated by adding simple energy loss terms in the objective function of the optimal control problem. For example, this adjustment could include Ohmic efficiency losses~\citep{Tamilselvi:2021aa} as a polynomial function of the battery discharge ($P_A + P_S - P_R$), and state-of-charge losses~\citep{Edge:2021aa} as an exponential function of the normalized battery charge $A$, which is  outside the current scope of the work.
\end{rem}

\subsection{Grid Power Model}

The network operator incurs the following cost of buying power from the traditional grid throughout the day:

\begin{equation}
	\label{eqn:cost_power}
	\mathcal{C}_1 = \E{\int_0^T K_b(t) P_F(t) \dd t} \cdot
\end{equation}

The power from the traditional grid is primarily generated by fossil fuel stations, and has a highly negative environmental effect. This work does not impose an upper limit on how much power could be bought from the grid. To restrict the green cellular network from buying substantial amounts of grid power, we also impose an additional "environmental" cost of buying fossil fuel-based power~\citep{farooq2015_optimized}.

\begin{equation}
	\label{eqn:cost_environment}
	\mathcal{C}_2 = \E{\int_0^T \left( C_1 P_F(t) + C_2 P_F(t)^2 \right) \dd t},
\end{equation}

where $C_1 \in \mathbb{R}^+$ and $C_2 \in \mathbb{R}^+$ denote the pollutant emission coefficients of grid power.

\subsection{Running Horizon Framework}
\label{sec:running_horizon}

Although this work aims to formulate an optimal power procurement problem for a one-day operation cycle, such a short-term approach could lead to spurious results at the end of the day. For example, the resulting optimal policy might have operators selling all stored energy at the end of the day to generate revenue. This scenario is unrealistic because the operator would prefer to keep some charge in the battery for the next day. This work devises a running horizon framework to avoid producing such impractical optimal policies and still work in a short-term finite horizon. The framework consists of two facets. First, we add fictitious revenue from storing battery charge at the end of the day. More specifically, we introduce a third revenue source, as follows:

\begin{equation}
\label{eqn:fictitious_revenue}
	\mathcal{R}_3 = \E{P_K \bar{A} A(T)},
\end{equation}

where $P_K \in \mathbb{R}^+$ represents the fictitious cost per unit of stored battery charge. For example, publicly available deterministic forecasts of normalized renewable power and day-ahead energy price forecasts can be applied to determine a suitable value of $P_K$. 

Next, we formulate the optimal control problem for a two-day operation cycle (instead of one) but apply the optimal power procurement policy only for one day. This approach requires two-day-ahead forecasts of renewable power, energy prices, and cellular user demand. The optimal power procurement policy is updated daily, as and when the operators receive the two-day-ahead forecasts of these quantities. Regarding the problem formulation, the only change is that the time horizon is denoted by $[0,T]$, where $0$ denotes beginning of the day, and $T$ denotes the end of the second day (instead of at the end of the current day). Henceforth, this notation is applied while formulating the optimization problem. 

\subsection{Model Summary}
\label{sec:model_summary}

The state and control vectors are denoted by $\mathbf{X}(t)$ and $\boldsymbol\phi(t)$, respectively:

\begin{align}
	\mathbf{X}(t) &:= \left[ A(t),R(t),\xi(t) \right], \\
	\boldsymbol\phi(t) &:= \left[ P_A(t), P_F(t), P_\mathrm{tx}(t), P_S(t) \right] \cdot
\end{align}

The state vector $\mathbf{X}(t)$ is defined on a probability space $(\Omega,\mathcal{F},P)$ and each of its components satisfy the controlled dynamics in \eqref{eqn:gamma_fading}, \eqref{eqn:re_dynamics}, and \eqref{eqn:battery_dynamics}, respectively in $[0,T]$. The initial conditions are also random. The control $\boldsymbol\phi(t) := \boldsymbol\phi(t,\mathbf{X}(t))$ is an $\mathcal{F}_t$-adapted Markovian control that is measurable, where $\{\mathcal{F}_t,t \geq 0\}$ represents the filtration on the probability space. 

\constraint{ \textit{ (A.s. instantaneous control constraints)}} 
\label{constraint:as}
The controls must satisfy the following constraints for all $t \in [0,T]$:

\begin{align}
	\label{eqn:constraint1}
	0 &\leq P_A(t), \\
	\label{eqn:constraint2}
	0 &\leq P_F(t), \\
	\label{eqn:constraint3}
	0 &\leq P_\mathrm{tx}(t) \leq \frac{\bar{P}_\mathrm{tx}}{N_u(t)}, \\
	\label{eqn:constraint4}
	0 &\leq P_S(t) \cdot
\end{align}

The controls must also follow the following state-dependent constraint:

\begin{equation}
	\label{eqn:constraint5}
	-\underline{P}_A(A(t)) \leq P_A(t) + P_S(t) - \bar{P}_R R(t) \leq \bar{P}_A(A(t)) \cdot
\end{equation}

Furthermore, for a given deterministic forecast of the number of cellular users connected to the network, $N_u(t)$, the controls are constrained to satisfy the following:

\begin{equation}
\label{eqn:constraint6}
	P_A(t) + P_F(t) = C_\text{scal} N_u(t) P_\mathrm{tx}(t) + C_\text{offset}, \quad 0 \leq t \leq T \cdot
\end{equation}

\constraint{ \textit{ (Probabilistic instantaneous control constraints)}}
\label{constraint:prob}
The controls must satisfy the following probabilistic constraint for all $t \in [0,T]$:

\begin{equation}
\label{eqn:constraint7}
	\prob{\phi_\text{out}(t) \geq \phi_\mathrm{th}} \leq \epsilon \cdot
\end{equation}

The shorthand notation $\phi_\text{out}(t) := \phi_\text{out}(t,\mathbf{X}(t),\boldsymbol\phi(t))$ is used here. The constraints in  \eqref{eqn:constraint1} to \eqref{eqn:constraint7} define the set of admissible Markov controls $\mathcal{A}(t,\mathbf{X}(t))$ at each time $t \in [0,T]$. Then, we define the set of admissible policies as follows:

\begin{equation}
\label{eqn:policy_defn}
	\bar{\mathcal{A}} = \left\{ \boldsymbol\phi : \boldsymbol\phi \text{ is } \mathcal{F}_t \text{-adapted}, \quad \boldsymbol\phi(t,\omega) \in \mathcal{A}(t,\mathbf{X}(t,\omega)), \quad \forall \omega \in \Omega, \quad \forall t \in [0,T] \right\} \cdot
\end{equation} 

\section{Stochastic Optimal Control Formulation}
\label{sec:soc_formulation}

\prblm{ \textit{ (Primal problem with the probabilistic constraint)}.}
\label{prob:primal}
Given the initial data $\mathbf{X}(0) = [A_0, R_0 \sim \mu_0^R, \xi_0 \sim \mu_0^\xi]$, deterministic forecasts for daily cellular user traffic profile $N_u(t)$, energy spot prices $K_b(t)$ and $K_s(t)$, renewable energy forecast $p(t)$, and mobile network usage price $K_\text{net}(t)$ for $0 \leq t \leq T$, we solve the following:

\begin{equation}
\label{eqn:primal_problem}
	\boldsymbol\phi^* = \argmin_{\{\boldsymbol\phi \in \bar{\mathcal{A}}\}} \mathcal{U}(\boldsymbol\phi) = w \left( \mathcal{C}_1 - \mathcal{R}_1 - \mathcal{R}_2 - \mathcal{R}_3 \right) + (1-w) \mathcal{C}_2,
\end{equation}

where the minimization is done over all controls satisfying the constraints in \ref{constraint:as} and \ref{constraint:prob} with dynamics  \eqref{eqn:gamma_fading}, \eqref{eqn:re_dynamics}, and \eqref{eqn:battery_dynamics}. We define $\mathcal{U}(\boldsymbol\phi)$ as the net utility of the base station. In addition, $w \in [0,1]$ is a Pareto parameter to weigh the two objectives: financial cost $(\mathcal{C}_1 - \mathcal{R}_1 - \mathcal{R}_2)$ and environmental cost ($\mathcal{C}_2$), where $w$ can be chosen by the base station operators, depending on how environmentally conscious they are. The dependence of the costs on the controls $\boldsymbol\phi$  in~\eqref{eqn:primal_problem} is suppressed for ease of notation.

\begin{rem}{\textit{ (Discrete vs. continuous-time formulation). }}
\label{rem:disc_vs_cont}
	The optimal solution to a discrete-time formulation leads to a discrete set of control decisions at each time step, dependent on the selected discretization. Applying the dynamic programming principle to a discrete-time formulation leads to a stochastic optimization problem at each time step. Finer grids substantially increase the computational work, making the problem computationally infeasible without guaranteeing that the optimal policies would converge. In a continuous-time formulation, the limit of the dynamic programming principle exists as the time step becomes infinitesimally small, leading to the derivation of the associated HJB PDE. Numerical discretization methods are applied just once to solve this PDE and obtain a solution at all time points. This approach facilitates discretization error control, ensuring convergence to the true optimal solution, and enables the flexible application of adaptive and higher-order numerical approximation schemes.  
\end{rem}

The probabilistic constraint in~\ref{constraint:prob} does not allow classical dynamic programming to solve  Problem~\ref{prob:primal} because the constraint controls the joint distribution of the Markovian dynamics, rather than each realization. The standard dynamic programming lemma~\citep{Bertsekas:2012aa} cannot be applied if the Markovian control function depended additionally on the state's distribution. We overcome this problem using Langrangian relaxation of the probabilistic constraint for all $t \in [0,T]$, yielding a continuous-time Lagrangian relaxation. For this approach, we define the time-continuous deterministic Lagrange multiplier function $\lambda: [0,T] \rightarrow \mathbb{R}^+$, associated with the relaxed constraint in~\ref{constraint:prob}. We let $\Lambda$ denote the admissible space of Lagrange multiplier functions. The new set of admissible controls is only defined by the constraints in~\eqref{eqn:constraint1} to~\eqref{eqn:constraint6} and is denoted by $\mathcal{A}^\mathrm{rel}(t,\mathbf{X}(t))$ at each time $t \in [0,T]$. Then, we define the new set of admissible policies as follows:

\begin{equation}
\label{eqn:policy_defn_v2}
	\bar{\mathcal{A}}^\mathrm{rel} = \left\{ \boldsymbol\phi : \boldsymbol\phi \text{ is } \mathcal{F}_t \text{-adapted}, \quad \boldsymbol\phi(t,\omega) \in \mathcal{A}^\mathrm{rel}(t,\mathbf{X}(t,\omega)), \quad \forall \omega \in \Omega, \quad \forall t \in [0,T] \right\} \cdot
\end{equation}

We define the Lagrangian using terminology from the classical Lagrangian relaxation technique for constrained optimization problems~\citep{boyd2004convex} as follows:

\begin{align}
	\label{eqn:lagrangian}
	\mathcal{L}(\boldsymbol\phi,\lambda) &= \mathcal{U}(\boldsymbol\phi) + \int_0^T \lambda(t) \left( \prob{\phi_\text{out}(t,\mathbf{X}(t),\boldsymbol\phi(t)) \geq \phi_\mathrm{th}} - \epsilon \right) \dd t \\
	&= \mathcal{U}(\boldsymbol\phi) + \int_0^T \E{\lambda(t) \left( \mathbbm{1}_{\{\phi_\text{out}(t,\mathbf{X}(t),\boldsymbol\phi(t)) \geq \phi_\mathrm{th}\}} - \epsilon \right)} \dd t \cdot \nonumber
\end{align}

\prblm{ \textit{ (Relaxed problem)}.}
\label{prob:relaxed}
Given the initial data $\mathbf{X}(0) = [A_0, R_0 \sim \mu_0^R, \xi_0 \sim \mu_0^\xi]$, Lagrange multiplier function $\lambda$,  deterministic forecasts for daily cellular user traffic profile $N_u(t)$, energy spot prices $K_b(t)$ and $K_s(t)$, renewable energy forecast $p(t)$, and mobile network usage price $K_\text{net}(t)$ for $0 \leq t \leq T$, we solve the following:

\begin{equation}
\label{eqn:relaxed_problem}
	\boldsymbol\phi^*(\lambda) = \argmin_{\{\boldsymbol\phi \in \bar{\mathcal{A}}^\mathrm{rel}\}} \mathcal{L}(\boldsymbol\phi,\lambda),
\end{equation}

where the minimization is done over all controls satisfying the constraint in~\ref{constraint:as} with dynamics \eqref{eqn:gamma_fading}, \eqref{eqn:re_dynamics}, and \eqref{eqn:battery_dynamics}.

With the relaxation of the probabilistic constraint in~\ref{constraint:prob}, the dynamic programming principle can be applied to solve Problem~\ref{prob:relaxed}. The dual problem associated with the primal Problem~\ref{prob:primal} is formulated below.

\prblm{ \textit{ (Dual problem)}.}
\label{prob:dual}
We determine

\begin{equation}
\label{eqn:dual_problem}
	\lambda^* = \argmax_{\lambda \in \Lambda} \Theta(\lambda),
\end{equation}

where the dual function $\Theta(\lambda)$ is given by

\begin{equation}
	\label{eqn:dual_function}
	\Theta(\lambda) = \mathcal{L}(\boldsymbol\phi^*(\lambda),\lambda),
\end{equation}

and $\boldsymbol\phi^*(\lambda)$ solves Problem~\ref{prob:relaxed} for a given $\lambda$.

The admissible set $\Lambda$ of Lagrange multiplier functions depends on the problem structure. For the considered optimization problem, $\lambda$ denotes a positive real-valued function of time $t$ due to the deterministic nature of the relaxed constraint~\eqref{eqn:constraint7}.

\subsection{HJB Equation Related to Problem~\ref{prob:relaxed}}
\label{sec:hjb}

As a first step to solving the dual problem (Problem~\ref{prob:dual}), Problem~\ref{prob:relaxed} is solved to obtain $\Theta(\lambda)$ for a given $\lambda$. Using terminology from standard continuous-time stochastic optimal control~\citep{Kushner:1990aa,Bertsekas:2012aa}, we first define the cost-to-go function $u:[0,T] \cross \Gamma \rightarrow \mathbb{R}$ associated with Problem~\ref{prob:relaxed} as follows:

\begin{align}
	\label{eqn:value_function}
	u(t,\mathbf{x}) &= \min_{\boldsymbol\phi \in \bar{\mathcal{A}}^\mathrm{rel}} \mathcal{U}_{t,\mathbf{x}}(\boldsymbol\phi) \\
	&= \min_{\boldsymbol\phi \in \bar{\mathcal{A}}^\mathrm{rel}} \mathbb{E} \Bigg[ \Bigg. \int_t^T \Bigg( \Bigg. w \Bigg( K_b(s) P_F(s) - K_s(s) P_S(s) - K_\text{net}(s) N_u(s) (1 - \phi_\text{out}(s,\xi(s),P_\mathrm{tx}(s))) \Bigg) \nonumber \\
	&+ (1-w) \left( C_1 P_F(s) + C_2 P_F(s)^2 \right) + \lambda(s) \left( \mathbbm{1}_{\{\phi_\text{out}(s,\xi(s),P_\mathrm{tx}(s)) \geq \phi_{\mathrm{th}}\}} - \epsilon \right) \Bigg. \Bigg) \dd s \nonumber \\
	&\qquad - P_K \bar{A} A(T) \Bigg| \mathbf{X}(t) = \mathbf{x} \Bigg. \Bigg] \cdot \nonumber
\end{align}

We define the components of $\mathbf{x} = [a,r,\chi]$, where $a$, $r$, and $\chi$ denote the variables corresponding to components of the state vector $A(t)$, $R(t)$, and $\xi(t)$ respectively. In addition, $\Gamma$ denotes the domain of $\mathbf{x}$, in this case $[0,1] \cross [0,1] \cross [\underline{\xi},\infty)$, because $A(t)$ and $R(t)$ are $[0,1]$-valued random variables and $\xi(t)$ is a $[\underline{\xi},\infty)$-valued random variable for all $t \in [0,T]$ (see Equations \eqref{eqn:gamma_fading}, \eqref{eqn:re_dynamics} and~\eqref{eqn:battery_dynamics}). Using \eqref{eqn:value_function}, the dual function is defined as $\Theta(\lambda) = u(0,\mathbf{X}(0))$, where $u$, for given $\lambda$, solves the following second-order nonlinear HJB final value PDE.

\begin{empheq}[left=\empheqlbrace, right = \cdot]{equation}
	\label{eqn:hjb_pde}
	\begin{alignedat}{2}
		&\frac{\partial u}{\partial t} + \mathcal{H}(t,\mathbf{x},\nabla u, \nabla^2 u; \lambda) = 0 \\
		&u(T,\mathbf{x}) = -P_K \bar{A} a, \quad \forall \mathbf{x} \in \Gamma,
	\end{alignedat}
\end{empheq}

where $\nabla \cdot$ denotes the gradient vector, and $\nabla^2 \cdot$ represents the Hessian matrix of a scalar-valued function. The Hamiltonian $\mathcal{H}$ associated with Problem~\ref{prob:relaxed} in \eqref{eqn:hjb_pde} is defined as follows:

\begin{align}
	\mathcal{H}(t,\mathbf{x},\nabla u, \nabla^2 u; \lambda) &= \min_{\boldsymbol\phi \in \mathcal{A}^\mathrm{rel}(t,\mathbf{x})} \Bigg. \Bigg[ \frac{\left( \bar{P}_R r - P_S - P_A \right)}{\bar{A}} \frac{\partial u}{\partial a} + \left( \dot{p}(t) - \theta(t) (r - p(t)) \right) \frac{\partial u}{\partial r} \nonumber \\
	& + \alpha \theta_0 r (1 - r) \frac{\partial^2 u}{\partial r^2} - \theta \left( \chi - \underline{\xi} - \mu \right) \frac{\partial u}{\partial \chi} + \theta (\chi - \underline{\xi}) \frac{\partial^2 u}{\partial \chi^2} + w \Bigg( \Bigg. K_b(t) P_F  \nonumber \\
	&- K_s(t) P_S - K_\text{net}(t) N_u(t) (1 - \phi_\text{out}(t,\chi,P_\mathrm{tx})) \Bigg. \Bigg) + (1-w) \left( C_1 P_F + C_2 P_F^2 \right) \nonumber \\
	\label{eqn:hamiltonian}
	& + \lambda(t) \left( \mathbbm{1}_{\{\phi_\text{out}(t,\chi,P_\mathrm{tx}) > \phi_{\mathrm{th}}\}} - \epsilon \right) \Bigg. \Bigg] \cdot
\end{align}

We obtain an approximation $\bar{\Theta}(\lambda)$ of the dual function $\Theta(\lambda)$ by numerically approximating the solution to \eqref{eqn:hjb_pde}. Section~\ref{sec:numerical_methods} discusses this approach in detail.

\subsection{Finite-dimensional Approximation of Problem~\ref{prob:dual}}
\label{sec:finite_dim_dual}

Problem~\ref{prob:dual} is an infinite-dimensional optimization problem to satisfy the constraint in~\ref{constraint:prob} pointwise in time. We discretize $\lambda(t)$ into a piecewise constant function to approximate the infinite-dimensional dual problem, allowing for practical numerical optimization. This approach implies satisfying the constraint for a finite number of subintervals in $[0,T]$. That is, the integral of the pointwise violation of the constraint in~\ref{constraint:prob} over each subinterval must be $0$. We consider the time discretization $0 = t_0 < t_1 < \ldots < t_{\ell-1} < t_\ell = T$ of the time domain $[0,T]$. The Lagrange multiplier is approximated as follows:

\begin{equation}
	\label{eqn:finite_dim_multiplier}
	\lambda(t) \approx \lambda^\ell(t) = \sum_{i=1}^{\ell} \Upsilon^\ell_i \mathbbm{1}_{\{t \in [t_{i-1},t_{i}]\}}, \quad \forall t \in [0,T],
\end{equation}

where $\{\Upsilon_i^\ell\}_{i=1}^{\ell}$ denotes components of the vector of Lagrange multipliers $\Upsilon^\ell$. For the given time discretization, we define 

\begin{equation}
	\label{eqn:finite_dim_lambda_space}
	\bar{\Lambda}^\ell = \{ \lambda^\ell: \text{with } \lambda^\ell(t) \text{ defined in \eqref{eqn:finite_dim_multiplier} } \text{ and } \Upsilon^\ell_i \in \mathbb{R}^+ \forall i \in [1,\ldots,\ell]  \} \cdot
\end{equation}
With this approximation, the finite-dimensional approximation of Problem~\ref{prob:dual} is formulated.

\prblm{ \textit{ (Finite-dimensional dual problem)}.}
\label{prob:finite_dim_dual}
For a given $\ell$, we determine the following:

\begin{equation}
\label{eqn:finite_dual_problem}
	\bar{\lambda}^\ell = \argmax_{\lambda^\ell \in \bar{\Lambda}^\ell} \Theta(\lambda^\ell) \cdot
\end{equation}

The following Lagrangian is written for a given $\lambda^\ell \in \bar{\Lambda}^\ell$:

\begin{equation}
	\label{eqn:finite_dim_lagrangian}
	\mathcal{L}(\boldsymbol\phi,\lambda^\ell) = \mathcal{U}(\boldsymbol\phi) + \sum_{i=1}^{\ell} \Upsilon^\ell_i \left( \mathfrak{D} \Theta (\lambda^\ell) \right)_i 
\end{equation}

with 

\begin{equation}
	\label{eqn:subgradient}
	\left( \mathfrak{D} \Theta (\lambda^\ell) \right)_i =  \E{\int_{t_{i-1}}^{t_{i}} \left( \mathbbm{1}_{\{\phi_\text{out}(t,\mathbf{X}(t),\boldsymbol\phi(t)) \geq \phi_\mathrm{th}\}} - \epsilon \right) \dd t}, \quad \forall i \in [1,\ldots,\ell] ,
\end{equation}

where $\left( \mathfrak{D} \Theta (\lambda^\ell) \right)_i$ denotes the $i^\mathrm{th}$ component of an $\ell$-dimensional subgradient vector $\mathfrak{D} \Theta$ of the dual function $\Theta$ evaluated at $\lambda^\ell \in \bar{\Lambda}^\ell$. To solve  Problem~\ref{prob:finite_dim_dual}, the optimal $\ell$-dimensional vector $\Upsilon^\ell \in \left( \mathbb{R}^+ \right)^\ell$ that maximizes the dual function $\Theta$ is selected. Hence, Problem~\ref{prob:finite_dim_dual} is a convex, $\ell$-dimensional optimization problem but is not generally smooth. Therefore, we apply subgradient methods to solve Problem~\ref{prob:finite_dim_dual}. Moreover, the choice of $\ell$ and the corresponding time discretization grid is set to control the pointwise violation in the constraint in~\ref{constraint:prob}. Section~\ref{sec:numerical_methods} details this approach.

\section{Numerical Approach}
\label{sec:numerical_methods}

\subsection{Numerically Solving the HJB PDE}
\label{sec:hjb_numerics}

We employ the explicit Euler upwind finite-difference scheme to numerically approximate the solution to the HJB equation in~\eqref{eqn:hjb_pde}. This method provides a convergent numerical scheme under certain conditions that can be proven for the considered problem~\citep{Sun:2015aa}. 

We consider the following time and space discretization of $(t,\mathbf{x}) \in [0,T] \cross [0,1] \cross [0,1] \cross [\underline{\xi},\infty)$, where $\mathbf{x} = (a,r,\chi)$, and the following finite-difference grid $\tau = \tau_t \times \tau_a \times \tau_r \times \tau_\chi$:

\begin{align}
\label{eqn:fd_discretization}
	&\tau_t : 0=t_0 < t_1 < \ldots < t_{N_t-1} < t_{N_t} = T, \quad \Delta t = \frac{T}{N_t} \\
	&\tau_a : 0=a_0 < a_1 < \ldots < a_{N_1-1} < a_{N_1} = 1, \quad \Delta a = \frac{1}{N_1} \nonumber \\
	&\tau_r : 0=r_0 < r_1 < \ldots < r_{N_2-1} < r_{N_2} = 1, \quad \Delta r = \frac{1}{N_2} \nonumber \\
	&\tau_\chi : 0 = \chi_0 < \chi_1 < \ldots < \chi_{N_3-1} < \chi_{N_3} = 1, \quad \Delta \chi = \frac{1}{N_3} \cdot \nonumber
\end{align}

For the dimension $\chi$, we approximate the domain $[\underline{\xi},\infty)$ by $[\underline{\xi},\bar{\chi}]$ for an appropriately chosen $\bar{\chi}$, which is scaled to $[0,1]$. We employ uniform grids in all dimensions; that is, 

\begin{align*}
	t_n &= n \Delta t, \quad n=0,\ldots,N_t, \\
	a_i &= i \Delta a, \quad i=0,\ldots,N_1, \\
	r_j &=j \Delta r, \quad j=0,\ldots,N_2, \\
	\chi_k &= k \Delta \chi, \quad k=0,\ldots,N_3 \cdot \\
\end{align*}

We combine the scaled drift and diffusion terms of dynamics~\eqref{eqn:gamma_fading}, \eqref{eqn:re_dynamics}, and \eqref{eqn:battery_dynamics} into vector-valued functions $F$ and $G$ respectively.

\begin{equation*}
	F(t,\mathbf{x},\boldsymbol\phi) = \begin{bmatrix}
		\frac{\bar{P}_R r - P_S - P_A}{\bar{A}} \\
		\dot{p}(t) - \theta(t) \left( r - p(t) \right) \\
		-\theta \left( \chi - \frac{\mu}{\bar{\chi} - \underline{\xi}} \right)
		\end{bmatrix}, \quad
	G(t,\mathbf{x},\boldsymbol\phi) = \begin{bmatrix}
		0 \\
		\alpha \theta_0 r (1-r) \\
		\frac{\theta \chi}{\bar{\chi} - \underline{\xi}}
		\end{bmatrix} \cdot
\end{equation*}

For convenience, we define the following function representing the running cost in~\eqref{eqn:hamiltonian}. 

\begin{align*}
	H(t,\mathbf{x},\boldsymbol\phi) &= w \left( K_b(t) P_F - K_s(t) P_S - K_\text{net}(t) N_u(t) (1 - \phi_\text{out}(t,\chi,P_\mathrm{tx})) \right) \\
	&+ (1-w) \left( C_1 P_F + C_2 P_F^2 \right) + \lambda(t) \left( \mathbbm{1}_{\{\phi_\text{out}(t,\chi,P_\mathrm{tx}) \geq \phi_{\mathrm{th}}\}} - \epsilon \right) \cdot
\end{align*}

For a function $f$, we define $f^+(t,\mathbf{x},\boldsymbol\phi) := \max \left( f(t,\mathbf{x},\boldsymbol\phi),0 \right)$ and $f^-(t,\mathbf{x},\boldsymbol\phi) := \max \allowbreak \left( -f(t,\mathbf{x} ,\boldsymbol\phi), 0 \right)$. For a function $f$, we define $f_{(i,j,k)}^n := f(t_n,[a_i,r_j,\chi_k],\boldsymbol\phi(t_n,[a_i,r_j,\chi_k]))$. Let $\bar{u}^n_{(i,j,k)}$ denote the numerical approximation of the cost-to-go function $u(t,\mathbf{x})$ at the point $(t_n,a_i,r_j,\chi_k)$ on the grid $\tau$. Then, the explicit update rule for the upwind finite-difference scheme is given as follows:

\begin{align}
	\label{eqn:linear_update}
	\bar{u}_{(i,j,k)}^{n-1} &= \bar{u}_{(i,j,k)}^n \Bigg. \Bigg( 1 - 2 (G_1)_{(i,j,k)}^n \frac{\Delta t}{\Delta r^2} - 2 (G_2)_{(i,j,k)}^n \frac{\Delta t}{\Delta \chi^2} \nonumber \\
	&\qquad - \frac{\Delta t}{\Delta a} \abs{(F_1)_{(i,j,k)}^n} - \frac{\Delta t}{\Delta r} \abs{(F_2)_{(i,j,k)}^n} - \frac{\Delta t}{\Delta \chi} \abs{(F_3)_{(i,j,k)}^n} \Bigg. \Bigg) \\
	&+ \bar{u}_{(i+1,j,k)}^n \frac{\Delta t}{\Delta a} (F_1^+)_{(i,j,k)}^n + u_{(i-1,j,k)}^n \frac{\Delta t}{\Delta a} (F_1^-)_{(i,j,k)}^n \nonumber \\
	&+ \bar{u}_{(i,j+1,k)}^n \left((G_2)_{(i,j,k)}^n \frac{\Delta t}{\Delta r^2} + (F_2^+)_{(i,j,k)}^n \frac{\Delta t}{\Delta r} \right) \nonumber \\
	&+ \bar{u}_{(i,j-1,k)}^n \left( (G_2)_{(i,j,k)}^n \frac{\Delta t}{\Delta r^2} + (F_2^-)_{(i,j,k)}^n \frac{\Delta t}{\Delta r} \right) \nonumber \\
	&+ \bar{u}_{(i,j,k+1)}^n \left( (G_3)_{(i,j,k)}^n \frac{\Delta t}{\Delta \chi^2} + (F_3^+)_{(i,j,k)}^n \frac{\Delta t}{\Delta \chi} \right) \nonumber \\
	&+ \bar{u}_{(i,j,k+1)}^n \left( (G_3)_{(i,j,k)}^n \frac{\Delta t}{\Delta \chi^2} + (F_3^-)_{(i,j,k)}^n \frac{\Delta t}{\Delta \chi} \right) + \Delta t H_{(i,j,k)}^n \cdot \nonumber \\
	\bar{u}_{(i,j,k)}^{N_t} &= -P_K \bar{A} a_i, \quad \forall i,j,k \cdot
\end{align}

In \eqref{eqn:linear_update}, we apply the notation $\{F_i\}_{i=1}^3$ and $\{G_i\}_{i=1}^3$ to denote the $i^\mathrm{th}$ component of the vector-valued functions $F$ and $G$, respectively. Optimal controls at grid point $(t_n,a_i,r_j,\chi_k)$, denoted by $(\boldsymbol\phi^*)_{(i,j,k)}^n$, are computed using the known derivatives (up to the second order) of the numerical cost-to-go-function as follows:

\begin{align}
	\label{eqn:numerical_controls}
	(\boldsymbol\phi^*)_{(i,j,k)}^n &= \argmin_{\boldsymbol\phi \in \mathcal{A}^\mathrm{rel}(t_n,[a_i,r_j,\chi_k])} \Bigg. \Bigg[ F^+_1(t_n,[a_i,r_j,\chi_k],\boldsymbol\phi) (\Delta^+_a \bar{u})_{(i,j,k)}^n + F_1^-(t_n,[a_i,r_j,\chi_k],\boldsymbol\phi) (\Delta^-_a \bar{u})_{(i,j,k)}^n \\
	&\quad + H(t_n,[a_i,r_j,\chi_k],\boldsymbol\phi) \Bigg. \Bigg] , \nonumber
\end{align}

where $(\Delta^+_a \bar{u})_{(i,j,k)}^n$ and $(\Delta^-_a \bar{u})_{(i,j,k)}^n$ represent the upwind derivatives given by 

\begin{equation}
\label{eqn:numerical_derivative}
	(\Delta^+_a \bar{u})_{(i,j,k)}^n = \frac{\bar{u}_{(i+1,j,k)}^n - \bar{u}_{(i,j,k)}^n}{\Delta a}, \quad (\Delta^+_a \bar{u})_{(i,j,k)}^n = \frac{\bar{u}_{(i,j,k)}^n - \bar{u}_{(i-1,j,k)}^n}{\Delta a} \cdot
\end{equation}

In general, boundary conditions must be imposed on the HJB PDE~\eqref{eqn:hjb_pde} to obtain a unique numerical solution. However, boundary conditions at the boundaries $a=0,1$, $r=0,1$, and $\chi=0$ are unnecessary in this setting because the drift terms $\{F_i\}_{i=1}^3$ naturally have a sign at the boundaries ensuring the use of an interior point to approximate first-order derivatives. Similarly, the diffusion terms $\{G_i\}_{i=1}^3$ are naturally $0$ at these boundaries, ensuring no exterior points are necessary to approximate second-order derivatives. At the boundary $\chi=\bar{\chi}$, we impose a nonreflective boundary condition to indicate that the actual domain is not bounded at $\chi=\bar{\chi}$. We impose the Courant--Friedrichs--Lewy (CFL) condition at each grid point in $\tau$ to ensure the stability of the numerical scheme

\begin{equation}
	\label{eqn:cfl}
	\Delta t \left( \frac{\abs{(F_1)_{(i,j,k)}^n}}{\Delta a} + \frac{\abs{(F_2)_{(i,j,k)}^n}}{\Delta r} + \frac{\abs{(F_3)_{(i,j,k)}^n}}{\Delta \chi} + \frac{2 (G_2)_{(i,j,k)}^n}{\Delta r^2} + \frac{2 (G_3)_{(i,j,k)}^n}{\Delta \chi^2}  \right) \leq 1 \cdot
\end{equation}

Although the HJB PDE~\eqref{eqn:hjb_pde} is nonlinear, the numerical scheme given by Eqs,~\eqref{eqn:linear_update} to \eqref{eqn:numerical_controls} is explicit because the cost-to-go function at time $t_{n-1}$ is completely determined by its value and derivatives at time $t_n$. However, this relationship is nonlinear and cannot be vectorized. Moreover, the four-dimensional (4D)  constrained optimization problem in~\eqref{eqn:numerical_controls} must be solved at every grid point in $\tau$. Depending on the structure of $\phi_\mathrm{out}(t,\chi,P_\mathrm{tx})$, \eqref{eqn:numerical_controls} can be solved analytically or numerically. Linear interpolation is applied to extend the numerical solution $\bar{u}$ from the grid $\tau$ to the entire domain $[0,T] \times \Gamma$. Next, the dual function $\Theta(\lambda)$ is approximated by $\bar{\Theta}(\lambda) = \bar{u}(0,\mathbf{X}(0))$ for a given $\lambda$.  Algorithm~\ref{alg:hjb_solver} in Appendix~\ref{appendix:alg} presents a pseudo-algorithm of the numerical solver.

\begin{rem}{\textit{ (Discretization error of the PDE solver). }}
	\label{rem:hjb_error}
	The upwind finite-difference scheme in Section~\ref{sec:hjb_numerics} is a first-order explicit scheme with a discretization error of $\order{\Delta t + \Delta a + \Delta r + \Delta \chi}$. When solving the overall optimization problem up to a relative tolerance of $\mathrm{TOL}$ with respect to the dual function value $\Theta(\lambda)$, the PDE discretization error must be $\order{\abs{\Theta(\lambda)} \mathrm{TOL}}$. Thus, $\Delta t, \Delta a, \Delta r, \Delta \chi$ must be $\order{\abs{\Theta(\lambda)} \mathrm{TOL}}$, implying the computational work of \\
$\order{\abs{\Theta(\lambda)}^{-4} \mathrm{TOL}^{-4}}$ to solve the HJB PDE once.
\end{rem}

\subsection{Estimating Subgradient $\mathfrak{D}\Theta$}
\label{sec:subgradient_estimation}

Upon solving the HJB PDE~\eqref{eqn:hjb_pde} numerically and obtaining an approximate cost-to-go function in the domain $[0,T] \times \Gamma$, we can numerically simulate optimally controlled paths of the state variable $\mathbf{X}(t)$ forward in time, from the initial condition $\mathbf{X}(0)$ using dynamics~\eqref{eqn:gamma_fading}, \eqref{eqn:re_dynamics}, and~\eqref{eqn:battery_dynamics}. We use the  Euler--Maruyama time discretization of the SDEs in~\eqref{eqn:gamma_fading}, \eqref{eqn:re_dynamics}, and~\eqref{eqn:battery_dynamics}. We consider the discretization $0=\bar{t}_0<\bar{t}_1<\ldots<\bar{t}_{\bar{N}_t}=T$ of the time domain $[0,T]$ with $\bar{N}_t$ uniform time steps. It follows that $\bar{t}_n = n \times \Delta \bar{t}, \quad n=0,1,\ldots,\bar{N}_t$, and $\Delta \bar{t} = \frac{T}{\bar{N}_t}$. If $\bar{\mathbf{X}}^{\bar{N}_t}$ is the time-discretized version of the stochastic process $\mathbf{X}$, then the Euler--Maruyama time discretization of the SDEs is expressed as follows:

\begin{empheq}[left=\empheqlbrace, right = \cdot]{equation}
	\label{eqn:sde_euler}
	\begin{alignedat}{2}
		\bar{\mathbf{X}}^{\bar{N}_t}(\bar{t}_{n+1}) &= \bar{\mathbf{X}}^{\bar{N}_t}(\bar{t}_{n}) + F(\bar{t}_{n},\bar{\mathbf{X}}^{\bar{N}_t}(\bar{t}_{n}),\boldsymbol\phi^*(\bar{t}_{n},\bar{\mathbf{X}}^{\bar{N}_t}(\bar{t}_{n}))) \Delta \bar{t} \\
		&\quad + G(\bar{t}_{n},\bar{\mathbf{X}}^{\bar{N}_t}(\bar{t}_{n}),\boldsymbol\phi^*(\bar{t}_{n},\bar{\mathbf{X}}^{\bar{N}_t}(\bar{t}_{n}))) : \sqrt{\Delta \bar{t}} \boldsymbol\varepsilon^n, \quad n=0,\ldots,\bar{N}_t-1 \\
		\bar{\mathbf{X}}^{\bar{N}_t}(\bar{t}_{0}) &= \mathbf{X}(0),
	\end{alignedat}
\end{empheq}

where $\boldsymbol\varepsilon^n \sim \mathcal{N}(\mathbf{0},\mathbb{I}_3)$ for $n=0,\ldots,\bar{N}_t-1$, denotes independent and identically distributed (i.i.d.) $3$-dimensional standard normal random variables, and $\cdot : \cdot$ denotes elementwise multiplication of two vectors. The optimal controls $\boldsymbol\phi^*(\bar{t}_{n},\bar{\mathbf{X}}^{\bar{N}_t}(\bar{t}_{n}))$ in~\eqref{eqn:sde_euler} are obtained by minimizing~\eqref{eqn:numerical_controls} at the current point of the state $\bar{\mathbf{X}}^{\bar{N}_t}(\bar{t}_{n})$, which may be a point outside the grid $\tau$. The values of the numerical derivatives $\Delta_a^{\pm} \bar{u}$ of the cost-to-go function at those points are linearly interpolated from the computed derivative values at the nearest points in grid $\tau$. The subgradient of $\Theta$ at $\lambda^\ell \in \bar{\Lambda}^\ell$ is computed according to~\eqref{eqn:subgradient}, approximating the expected value using Monte Carlo sampling with $M_\text{SG}$ i.i.d. sample paths of the optimally controlled process $\bar{\mathbf{X}}^{\bar{N}_t}$. Moreover, the integral in~\eqref{eqn:subgradient} is approximated using a forward Euler approximation. For each $i \in [1,2,\ldots,\ell]$, we consider the discretization $t_{i-1}=\tilde{t}_0<\tilde{t}_1<\ldots<\tilde{t}_{\tilde{N}_t}=t_{i}$ of the time domain $[t_{i-1},t_{i}]$, corresponding to the $i^\mathrm{th}$ component of the subgradient, with $\tilde{N}_t$ uniform time steps. It follows that $\tilde{t}_n = n \times \Delta \tilde{t}, \quad n=0,1,\ldots,\tilde{N}_t$, and $\Delta \tilde{t} = \frac{t_{i}-t_{i-1}}{\tilde{N}_t}$. The numerical approximation of the $i^\mathrm{th}$ component of the $\ell$-dimensional subgradient is denoted by $(\bar{\mathfrak{D}} \Theta(\lambda^\ell))_i$ and given by

\begin{equation}
	\label{eqn:numerical_subgradient}
	\left( \bar{\mathfrak{D}} \Theta (\lambda^\ell) \right)_i = \frac{1}{M_\text{SG}} \sum_{m=1}^{M_\text{SG}} \sum_{n=0}^{\tilde{N}_t}  \left( \mathbbm{1}_{\{\phi_\text{out}(\tilde{t}_n,\bar{\mathbf{X}}^{\bar{N}_t}(\tilde{t}_n,\omega^{(m)}),\boldsymbol\phi(\tilde{t}_n,\bar{\mathbf{X}}^{\bar{N}_t}(\tilde{t}_n,\omega^{(m)}))) \geq \phi_\mathrm{th}\}} - \epsilon \right) \Delta \tilde{t}, \forall i \in [1,\ldots,\ell], \cdot
\end{equation}

where $\omega^{(m)}$ denotes the $m^\mathrm{th}$ i.i.d. realization of the random variables required to generate sample paths of $\bar{\mathbf{X}}^{\bar{N}_t}$. The time discretization for subgradient computation and the Lagrange multiplier function $\lambda^\ell$ generally need not coincide. In this case, Brownian bridge interpolation is necessary to evaluate the optimally controlled paths $\bar{\mathbf{X}}^{\bar{N}_t}$ at the time discretization points $\{\tilde{t}_n\}_{n=0}^{\tilde{N}_t}$ for each $i \in [1,\ldots,\ell]$. The same $M_\text{SG}$ realisations of the optimally controlled path $\bar{\mathbf{X}}^{\bar{N}_t}$ are used to estimate all $\ell$ components of $\bar{\mathfrak{D}} \Theta (\lambda^\ell)$. Algorithm~\ref{alg:estimate_subgradient} in Appendix~\ref{appendix:alg} presents the corresponding pseudo-algorithm.

Moreover~\eqref{eqn:numerical_subgradient} is a Monte Carlo approximation of~\eqref{eqn:subgradient}. It is subject to noisy  statistical error controlled in a probabilistic sense by the choice of $M_\text{SG}$. Thus, we only have access to noisy subgradients of the dual function. Hence, we apply the stochastic subgradient method (SSM) to solve the convex, nonsmooth optimization problem~\ref{prob:finite_dim_dual}.

\begin{rem}{ \textit{ (Subgradient estimation error). }}
	\label{rem:mc_error}
	The error in the subgradient estimation consists of two parts: a statistical Monte Carlo error from approximating the expected value in~\eqref{eqn:subgradient} using a sample average with $M_\mathrm{SG}$ samples and a discretization error from approximating the integral in the expectation in~\eqref{eqn:subgradient} using a forward Euler summation with $\tilde{N}_t$ steps. From the central limit theorem~\citep{Lefebvre:2007aa}, the Monte Carlo approximation error is $\order{M_\text{SG}^{-\frac{1}{2}}}$. The discretization error for the first-order forward Euler summation is $\order{\tilde{N}_t^{-1}}$, yielding a total error in the subgradient estimation of $\order{M_\text{SG}^{-\frac{1}{2}} + \tilde{N}_t^{-1}}$. When solving the overall optimization problem up to a relative tolerance of $\mathrm{TOL}$ with respect to the subgradient, then we must estimate the subgradients up to an absolute tolerance of $\epsilon \mathrm{TOL}$. This requirement implies that $M_\text{SG} = \order{\epsilon^{-2} \mathrm{TOL}^{-2}}$ and $\tilde{N}_t = \order{\epsilon^{-1} \mathrm{TOL}^{-1}}$, yielding a total computational work of $\order{\epsilon^{-3} \mathrm{TOL}^{-3}}$ for each subgradient estimation.
\end{rem}

\subsection{Dual Problem Optimization}
\label{sec:solve_dual}

After approximating the dual function $\Theta(\lambda^\ell)$ and its subgradient $\mathfrak{D}\Theta(\lambda^\ell)$ for a given Lagrange multiplier function $\lambda^\ell \in \bar{\Lambda}^\ell$, we must solve the nonsmooth convex dual optimization problem in  Problem~\ref{prob:finite_dim_dual}. We solve this in two stages: (i) by constructing an approximation $\lambda^\ell$, of $\lambda$ (as in~\eqref{eqn:finite_dim_multiplier}), and (ii) by numerically optimizing the amplitudes $\Upsilon^\ell$ of the approximated function $\lambda^\ell$. Algorithm~\ref{alg:dual_optimization} illustrates the overall procedure.

\begin{algorithm}
	\caption{Numerical dual optimization procedure}
	\label{alg:dual_optimization}
	\SetAlgoLined
	\KwIn{\( \mathrm{TOL} \), $\mathrm{TOL}_\text{init}$, \texttt{max-iter} , $\bar{N}_\text{iter}$, $N_\text{iter}$, $\Upsilon^1_1=1$, $\beta_F,M_\mathrm{SG},N_t,N_1,N_2,N_3,\bar{N}_t,\tilde{N}_t,\mathbf{X}(0)$}
	\KwOut{Optimal controls $\boldsymbol\phi^*$, optimal Lagrange multiplier function \( \lambda^\ell(t) \)}
	Construct $\lambda^1(t)$ with $\Upsilon_1^1$ using~\eqref{eqn:finite_dim_multiplier}; \\
	Obtain $\tilde{\Upsilon}_1^1$ using initialization Algorithm~\ref{alg:initialization} with inputs $\mathrm{TOL}_\text{init}$,$\Upsilon^1_1=1$, $\beta_F,M_\mathrm{SG},N_t,N_1,N_2,N_3,\bar{N}_t,\mathbf{X}(0)$; \\
	Construct $\lambda^1(t)$ with $\tilde{\Upsilon}_1^1$ using~\eqref{eqn:finite_dim_multiplier}; \\
	Obtain $\hat{\Upsilon}_1^1$ with the LMBM routine~\citep{Karmitsa:2007aa} with starting point $\tilde{\Upsilon}_1^1$, number of iterations $N_\text{iter}$ and parameters specified in Appendix~\ref{tab:lmbm_parameters}; \\
	Construct $\lambda^1(t)$ with $\hat{\Upsilon}_1^1$ using~\eqref{eqn:finite_dim_multiplier}; \\
	\While{$\ell < \bar{N}_t$}{
		$\ell \leftarrow 2 \ell$; \\
		Compute $\bar{\Theta}(\lambda^\ell) = \bar{u}(0,\mathbf{X}(0))$ by solving~\eqref{eqn:hjb_pde} using Algorithm~\ref{alg:hjb_solver} with $\lambda^\ell(t)$ and parameters $N_t,N_1,N_2,N_3$; \\
		Estimate $\bar{\mathfrak{D}}\Theta(\lambda^\ell)$  using Algorithm~\ref{alg:estimate_subgradient} with parameters $\bar{N}_t,M_\mathrm{SG},\tilde{N}_t$; \\
		$k = 0$; \\
		\While{$k <$ \texttt{max-iter} and $\norm{\bar{\mathfrak{D}}\Theta(\lambda^\ell)} > \mathrm{TOL} \epsilon$}{
			$\begin{cases}
			\Upsilon^\ell \leftarrow \Upsilon^\ell + C_\mathrm{SSM} \frac{\bar{\mathfrak{D}}\Theta(\lambda^\ell)}{\norm{\bar{\mathfrak{D}}\Theta(\lambda^\ell)}}, \quad k \leq \bar{N}_\text{iter} \\
			\Upsilon^\ell \leftarrow \Upsilon^\ell + \frac{C_\mathrm{SSM}}{k+1} \frac{\bar{\mathfrak{D}}\Theta(\lambda^\ell)}{\norm{\bar{\mathfrak{D}}\Theta(\lambda^\ell)}}, \quad k > \bar{N}_\text{iter}
			\end{cases}$; \\
			Construct $\lambda^\ell(t)$ with $\Upsilon^\ell$ using~\eqref{eqn:finite_dim_multiplier}; \\
			Compute and store $\varpi^{(k)} = \bar{\Theta}(\lambda^\ell) = \bar{u}(0,\mathbf{X}(0))$ by solving~\eqref{eqn:hjb_pde} using Algorithm~\ref{alg:hjb_solver} with $\lambda^\ell(t)$ and parameters $N_t,N_1,N_2,N_3$; \\
			Estimate $\bar{\mathfrak{D}}\Theta(\lambda^\ell)$  using Algorithm~\ref{alg:estimate_subgradient} with parameters $\bar{N}_t,M_\mathrm{SG},\tilde{N}_t$; \\
			$k \leftarrow k + 1$; \\
		}
		Save $\Upsilon^\ell$ corresponding to $\max \{ \varpi(1), \ldots, \varpi(k)\}$; \\
		Construct $\lambda^{2\ell}(t)$ with $\Upsilon^\ell$ using~\eqref{eqn:finite_dim_multiplier}; \\
	} 
\end{algorithm}

\subsubsection{Lagrange Multiplier Refinement}
\label{sec:ell_refinement}

This work constructs a finite-dimensional piecewise constant approximation, $\lambda^\ell$, of the Lagrange multiplier function, $\lambda$, as detailed in Section~\ref{sec:finite_dim_dual}. Increasing the refinement level $\ell$ enhances the approximation quality of $\lambda$ by $\lambda^\ell$, but also increases the optimization dimension (also $\ell$). This work determines a sufficient $\ell$ that controls the violation of the relaxed constraint. We start with $\ell=1$ and keep uniformly doubling it, until the violation of the  constraint in~\ref{constraint:prob} is sufficiently controlled pointwise in time. 

\subsubsection{Numerical Optimization}
\label{sec:ssm_algorithm}

This work solves the optimization of $\Theta(\lambda^\ell)$ for $\lambda^\ell \in \bar{\Lambda}^\ell$ with respect to the amplitudes $\Upsilon^\ell \in (\mathbb{R}^+)^\ell$ as discussed in Section~\ref{sec:finite_dim_dual}. For a given evaluation point $\lambda^\ell$, we have access to the dual function value (approximated in Section~\ref{sec:hjb_numerics}) and a noisy subgradient (approximated in~\eqref{eqn:numerical_subgradient}). Hence, the SSM~\citep{Boyd:2008aa} is employed to solve the nonsmooth optimization problem. The SSM is a subgradient-based search algorithm using a dual function evaluation and noisy evaluation of an arbitrary subgradient at each evaluation point. At each refinement level $\ell$, the SSM runs until the subgradient norm is below a prescribed relative tolerance $\mathrm{TOL}$ or the number of iterations exceeds a threshold \texttt{max-iter}. The algorithm output at each $\ell$ is the optimal amplitude vector $\Upsilon^\ell$. However, this optimal value is not necessarily reached at the final iteration of the SSM because the objective function to be maximized does not increase at every step of the subgradient method~\citep{Boyd:2008aa}. Hence, we track the highest value of the dual function and store its corresponding $\Upsilon^\ell$ as the optimal value.

The choices for the starting point and step-size are crucial for quick convergence of the SSM, especially in higher dimensions. This work uses a nonsummable diminishing step-size of $\order{\frac{1}{k+1}}$, where $k$ denotes the iteration number. The convergence of the SSM with step-size is proven in~\citep{Boyd:2008aa}. The associated constant, denoted by $C_\mathrm{SSM}$ in Algorithm~\ref{alg:dual_optimization}, is tuned once such that the step-size taken is $\order{\Upsilon^\ell}$. $C_\mathrm{SSM}$ need not be tuned again for a different set of parameters with the current problem structure. For enhanced performance, the step-size at each iteration is adjusted by dividing the current subgradient by its norm. Furthermore, the SSM runs for the first $\bar{N}_\text{iter}$ iterations using a constant step-size before reverting to a diminishing step size. The Lagrange multiplier function $\lambda^{\ell-1}(t)$ constructed using the computed optimal amplitudes $\Upsilon^{\ell-1}$ is used as the starting point of the SSM at level $\ell$. However, we must still choose a good starting point at level $\ell=1$. This work devises an initialization algorithm and complements it with a nonsmooth deterministic optimizer.

\begin{rem}{\textit{ (SSM convergence rate). }}
	\label{rem:ssm_convergence}
	A convergence rate of $\order{k^{-\frac{1}{2}}}$ for the SSM has been proven for a class of convex functions with suitable step-size choice~\citep{Grimmer:2019aa}. When solving the optimization problem up to a relative tolerance of $\mathrm{TOL}$ with respect to the dual function value $\Theta(\lambda)$, the SSM requires an expected number of $\order{\abs{\Theta(\lambda)}^{-2} \mathrm{TOL}^{-2}}$ iterations to converge under the assumptions in~\citep{Grimmer:2019aa}.
\end{rem} 

\begin{rem}{\textit{ (Primal feasible solution). }}
	\label{rem:primal_feasible}
	Note that Algorithm~\ref{alg:dual_optimization} uses the subgradient as a stopping criterion, instead of a duality gap. This is because the output of Algorithm~\ref{alg:dual_optimization} is already a primal feasible solution up to given relative tolerance $\tol$. That is, the optimal controls produce a solution that minimizes the primal cost while only violating Constraint~\ref{constraint:prob} lesser than a small value of $\epsilon \tol$.
\end{rem}

\subsubsection{LMBM-boosted Initialization}
\label{sec:initialization}

Level $\ell=1$ implies that the Lagrange multiplier function $\lambda^1(t)$ is constant in time, with the constant denoted by $\Upsilon_1^1$ (see \eqref{eqn:finite_dim_multiplier}). First, this work devises an initialization algorithm starting with the arbitrary point $\Upsilon_1^1 = 1$ and estimating the dual function value and a subgradient at that point. Then, if the subgradient at $\Upsilon_1^1=1$ is positive, $\Upsilon_1^1$ is continually increased by a factor $\beta_F$ until the subgradient becomes negative or its norm reaches a prescribed relative tolerance $\mathrm{TOL}_\text{init}$. Conversely, if the subgradient at $\Upsilon_1^1=1$ is negative, $\Upsilon_1^1$ is continually decreased by a factor $\beta_F$, until the subgradient becomes positive or its norm reaches relative tolerance $\mathrm{TOL}_\text{init}$. The final obtained point is stored as $\tilde{\Upsilon}_1^1$. Algorithm~\ref{alg:initialization} in Appendix~\ref{appendix:alg} presents the pseudo-algorithm of the devised initialization procedure.

To obtain a better starting point, $\tilde{\Upsilon}_1^1$ is used as a starting point for a deterministic nonsmooth optimization routine called the LMBM~\citep{Haarala:2004aa,Karmitsa:2012aa}. Bundle methods are more robust than a simple subgradient-based method because they approximate the entire subdifferential of the objective function, enhancing its convergence speed~\citep{Makela:2002aa}. However, the LMBM requires access to deterministic estimates of the objective function value and its subgradient at each evaluation point. Although this work only has access to noisy Monte Carlo estimates of subgradients~\eqref{eqn:numerical_subgradient}, we can still employ the LMBM routine to reach a good starting point for the SSM, especially in the regime where changes in the objective function value considerably outweigh the noise in subgradient estimates. This work runs the LMBM for a fixed number of iterations $N_\text{iter}$. Table~\ref{tab:lmbm_parameters} in Appendix~\ref{appendix:sim_parameters} specifies the parameters required to run the LMBM routine~\citep{Karmitsa:2007aa}. The output of this routine $\hat{\Upsilon}_1^1$ is applied as a starting point for the SSM algorithm described in Section~\ref{sec:ssm_algorithm}.

\begin{rem}{\textit{ (Scale of $\lambda$). }}
	\label{rem:lambda_scale}
	With respect to the parameter $\epsilon$ that determines the confidence level of satisfying Constraint~\ref{constraint:prob}, the optimal primal cost $\mathcal{U}(\boldsymbol\phi)$ associated with Problem~\ref{prob:finite_dim_dual} is an $\order{1}$ term, while the optimal subgradient would be $\order{\epsilon}$ (ideally 0). Equation~\eqref{eqn:finite_dim_lagrangian} then implies that the optimal Lagrange multiplier $\lambda$ is $\order{\frac{1}{\epsilon}}$. This yields a good guess for an initial point $\Upsilon_1^1$ for the LMBM-boosted initialization procedure.
\end{rem}

\section{Numerical Experiments and Results}
\label{sec:results}

This section presents a model example of a cellular base station powered by the German power grid. This section describes the system and all dynamics driving the operation of the base station, and provides the results of applying the proposed numerical approach to solve the optimal power procurement problem for the system.

\subsection{Description of Model Cellular Base Station System}
\label{sec:model_description}

Figure~\ref{fig:base_station} schematically illustrates the considered base station model. Table~\ref{tab:model_parameters} provides the  descriptions and numerical values of all coefficients used to describe the model. 

The daily mobile user traffic profile primarily drives the power demand of a cellular base station. The description of the daily traffic profile consists of two facets: (i) the number of people connected to the network $N_u(t)$, and (ii) the physical distribution of users around the base station $\rho_z(t)$. The numerical experiment applies the following sinusoidal profile to model $N_u(t)$:

\begin{equation}
	\label{eqn:sinusoidal_traffic_profile}
	N_u(t) = \max \left[ \underline{N}_u, \bar{N}_u \frac{1}{2^\varrho} \left( 1 + \sin \left(\frac{\pi t}{6} + \pi \right) \right)^\varrho \right] \cdot
\end{equation}

\begin{figure}[h!]
	\centering
	\includegraphics[width=0.6\textwidth]{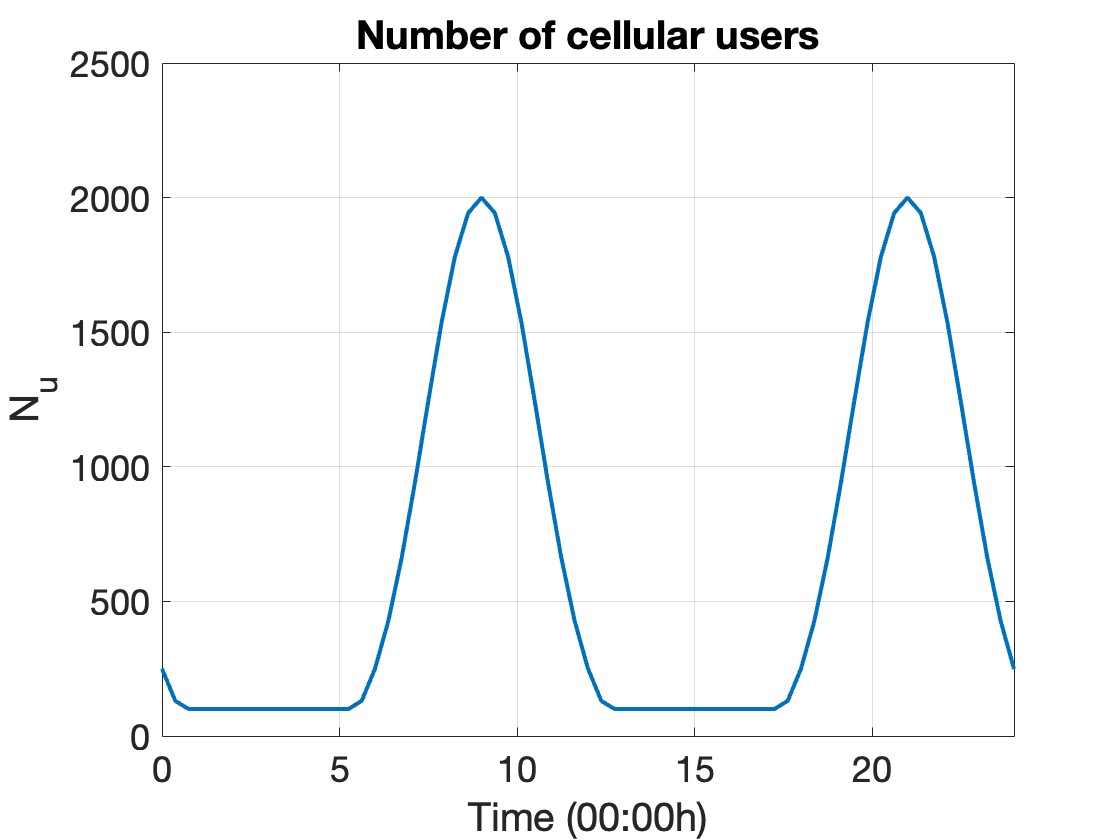}
	\caption{Typical daily cellular user traffic profile, described by \eqref{eqn:sinusoidal_traffic_profile}.}
	\label{fig:user_profile}
\end{figure}

This model and the stochastic versions of it have been empirically demonstrated to approximate practical user patterns closely for calibrated values of the smoothness parameter $\varrho$~\citep{Marsan:2010aa,Rached:2018aa}. The parameter $\varrho$ determines the rate of increase or decrease in $N_u(t)$ during the day. A higher value of $\varrho$ indicates a steeper increase or decrease rate, and $\varrho=3$ in the numerical experiment. This model also considers two peak hours (about 09.00 to 10.00 and 20.00 to 21.00) and few  connected users at night (00.00 to 05.00). Moreover, $\underline{N}_u$ and $\bar{N}_u$ denote the minimum and maximum number of users, respectively, that are connected to the base station at any time. This work sets $\underline{N}_u = 100$ and $\bar{N}_u = 2000$. These values were selected using the data on the average population densities of German cities~\cite{destatis2025} and the geographical area served by an average cellular base station in Germany~\cite{statista5G2023}. Figure~\ref{fig:user_profile} depicts the model for $N_u(t)$. 

\begin{figure}[h!]
    \centering
    % First subfigure
    \begin{subfigure}[b]{0.45\textwidth}
        \centering
        \includegraphics[width=\textwidth]{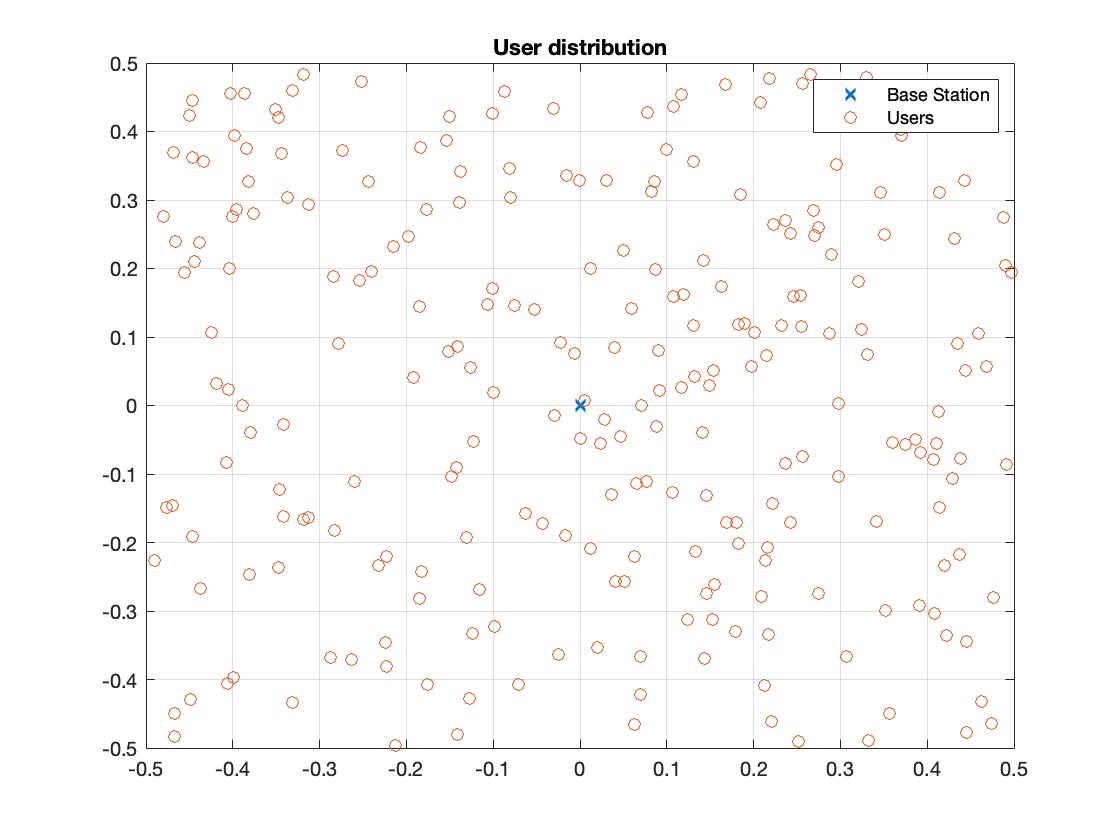}  % Replace with your image
        \caption{Uniform distribution in suburban areas.}
        \label{fig:uniform_user_dist}
    \end{subfigure}
    \hfill
    % Second subfigure
    \begin{subfigure}[b]{0.45\textwidth}
        \centering
        \includegraphics[width=\textwidth]{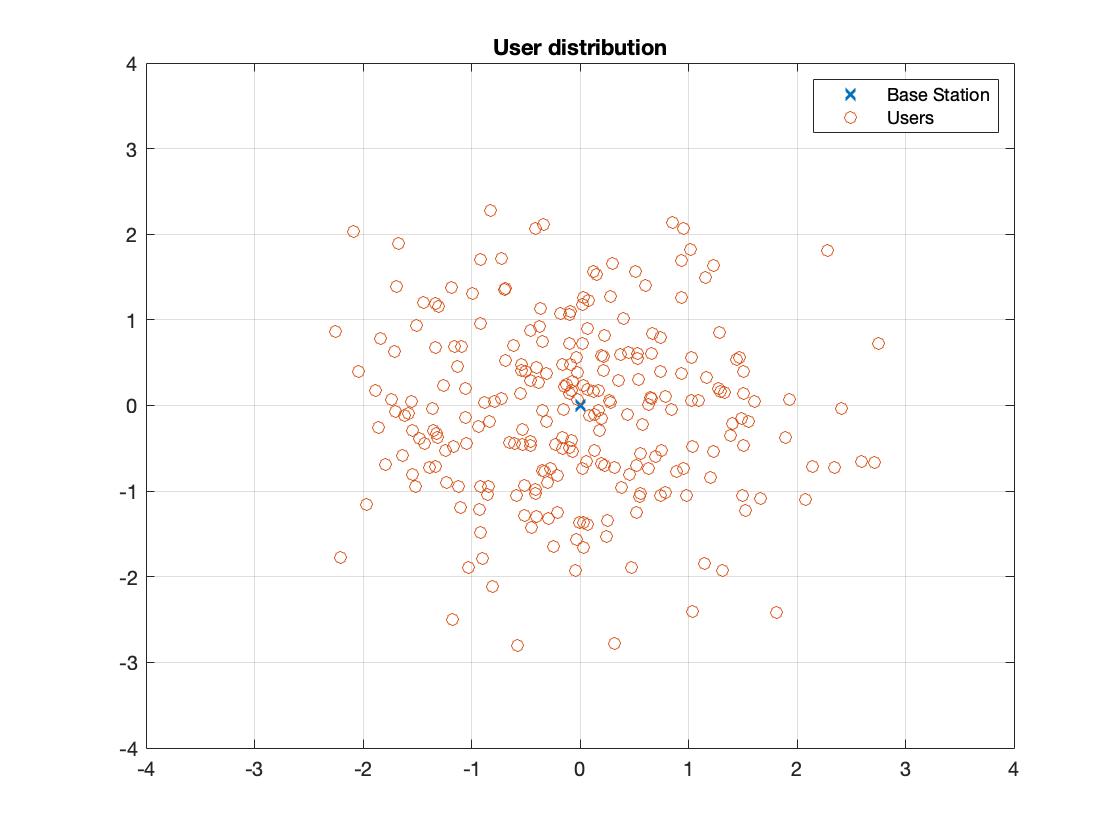}  % Replace with your image
        \caption{Gaussian distribution in industrial zones.}
        \label{fig:gaussian_user_dist}
    \end{subfigure}
    \caption{Visualizing the spatial distribution of mobile users in urban zones using simple analytical 2D distributions.}
    \label{fig:user_dist}
\end{figure}

The user distribution $\rho_z(t)$ crucially determines the estimation of the outage proportion $\phi_\text{out}(t)$ at given time $t$ (Section~\ref{sec:base_station_model}). Some simple 2D distributions have been employed to model the mobile user distribution in urban zones efficiently. For example, the uniform distribution can be employed to model the mobile user distribution in suburban areas very well, and the Gaussian distribution is a good model for mobile user distributions in industrial zones~\citep{Nguyen:2005aa}. Figure~\ref{fig:user_dist} visualizes this model. Analytical expressions can be derived for the corresponding $\phi_\text{out}(t)$ for these distributions. Appendix~\ref{appendix:outage_prop} provides the closed-form expressions. For the numerical experiments, we apply  the symmetric 2D Gaussian distribution centered on the site of the base station $\mathbf{x}_\mathrm{BS}$ with a diagonal covariance matrix with variance $\sigma_u^2 = 300^2 \text{m}^2$ along both directions. 

\begin{figure}[h!]
    \centering
    % First subfigure
    \begin{subfigure}[b]{0.45\textwidth}
        \centering
        \includegraphics[width=\textwidth]{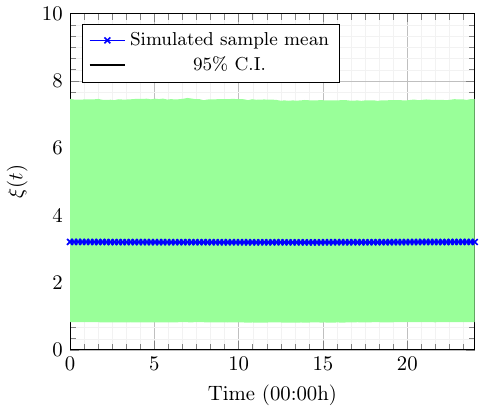}  % Replace with your image
        \caption{Plot of the sample mean and 95\% confidence intervals of sample paths of the SDE in~\eqref{eqn:gamma_fading}.}
        \label{fig:nakagami_paths}
    \end{subfigure}
    \hfill
    % Second subfigure
    \begin{subfigure}[b]{0.45\textwidth}
        \centering
        \includegraphics[width=\textwidth]{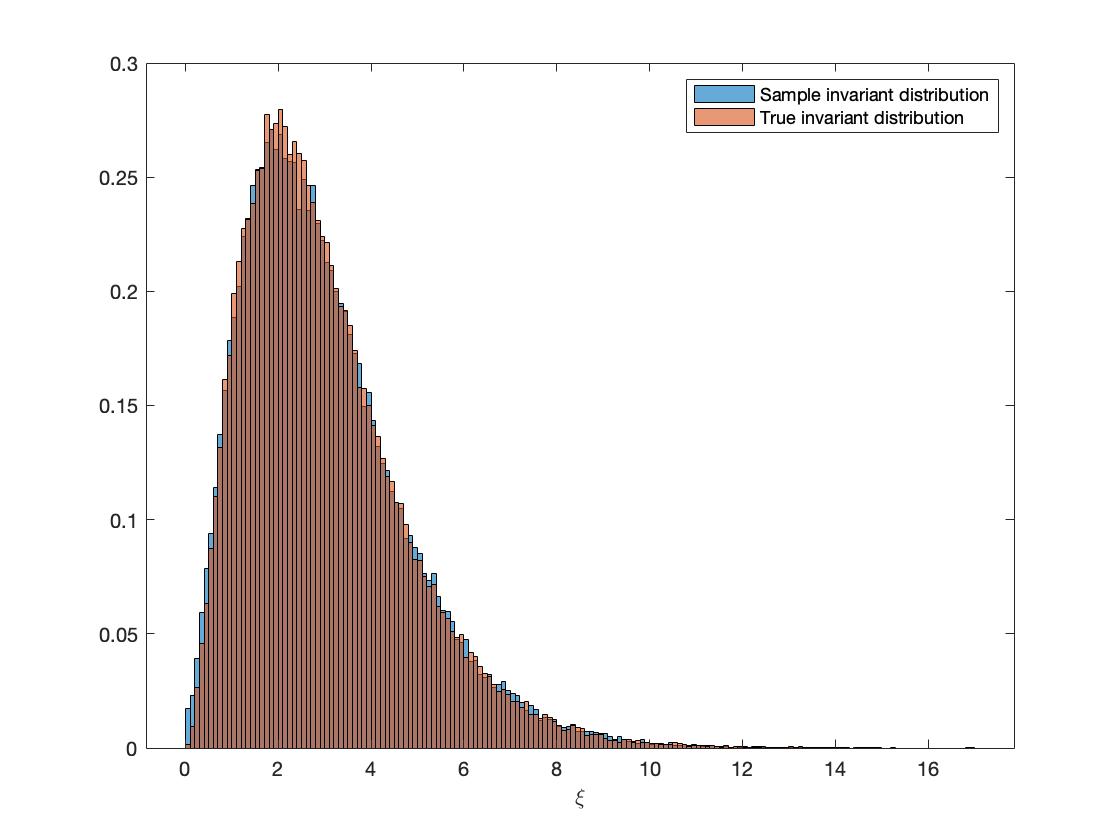}  % Replace with your image
        \caption{Visualization of the invariant distribution of sample paths of the SDE in~\eqref{eqn:gamma_fading}, verifying that the invariant distribution of the SDE in~\eqref{eqn:gamma_fading} is the shifted-gamma distribution with shape parameter $\mu=3$ and scale parameter $\theta=1$.}
        \label{fig:gamma_ergodic}
    \end{subfigure}
    \caption{Results of numerical simulation of the SDE in~\eqref{eqn:gamma_fading} governing wireless fading channel dynamics.}
    \label{fig:gamma_fading}
\end{figure}

We apply the SDE in~\eqref{eqn:gamma_fading} to model the instantaneous Nakagami fading channel, whose invariant distribution is a shifted-gamma distribution with the shape parameter $\mu$ and scale parameter $\theta$. We set $\mu=3$ and $\theta=1$ so that a quick, smooth convergence occurs to the invariant distribution and $\E{\xi(t)^{-2}}$ is bounded for all $t \in [0,T]$. We set the shift $\underline{\xi}$, corresponding to the value of $\xi$ that can still achieve an outage proportion of $\frac{\phi_\mathrm{th}}{2}$ in the worst-case scenario when $N_u(t) = \bar{N}_u$. Moreover, we set the initial distribution $\mu_0^\xi$ the same as the invariant distribution of process $\xi$. Figure~\ref{fig:gamma_fading} visualizes the samples paths and the invariant distribution of the SDE in~\eqref{eqn:gamma_fading}.

We set $K_\text{net}(t) = 0.01$ \texteuro/h per person for $t \in [0,T]$. The price was taken from the Vodafone 5G prices in Germany~\cite{vodafone5Gpricing}. The constant price assumption is generally valid because regulators set these prices and do not move within a day. 

We apply the SDE in~\eqref{eqn:re_dynamics} to model the instantaneous normalized generated renewable power. This work uses the 2023 German wind power forecast and production data with a 15-minute frequency from the operator \textit{50Hertz}~\citep{50hertzWindpower} to calibrate the SDE in~\eqref{eqn:re_dynamics}. The forecast data are employed to construct $p(t)$ and $\dot{p}(t)$, and the discrepancy between the forecast and true production data is applied to calibrate the parameters $\alpha$ and $\theta_0$ in~\eqref{eqn:re_dynamics}. Apart from these parameters, we also calibrate an additional parameter $\varsigma$. The discrepancy between the forecast and production at time $t=0$ is not typically zero in real-world applications. $\varsigma$ is defined as the time before $t=0$, for which the forecast error can be set to zero. This approach ensures that uncertainty also exists in renewable power production at $t=0$. We follow the calibration procedure detailed in~\citep{caballero21_derivative} (see~\citep{caballero21_derivative} for details). The mean-reversion parameter $\theta(t)$ in~\eqref{eqn:re_dynamics} is defined as follows:

\begin{equation*}
	\theta(t) = \max \left( \theta_0 , \frac{\alpha \theta_0 + \abs{\dot{p}(t)}}{\min (p(t),1-p(t))} \right) \cdot
\end{equation*}

\begin{figure}[h!]
    \centering
    % First subfigure
    \begin{subfigure}[b]{0.45\textwidth}
        \centering
        \includegraphics[width=\textwidth]{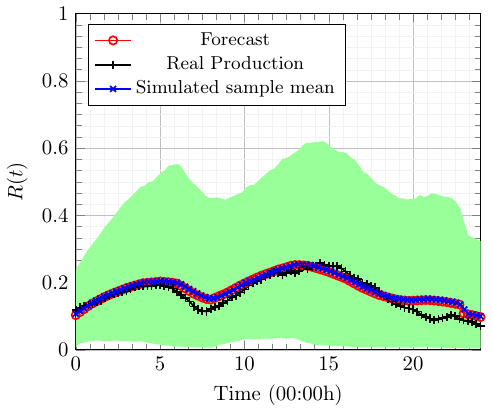}  % Replace with your image
        \caption{12.03.2024}
        \label{fig:wind_120324}
    \end{subfigure}
    \hfill
    % Second subfigure
    \begin{subfigure}[b]{0.45\textwidth}
        \centering
        \includegraphics[width=\textwidth]{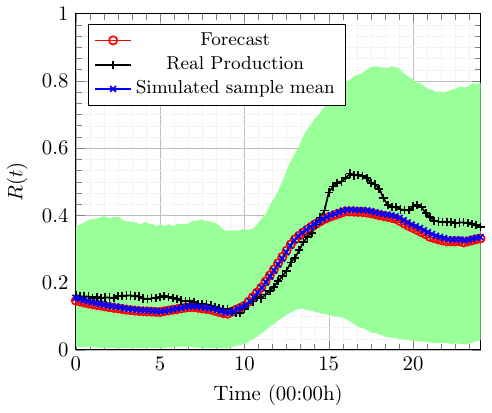}  % Replace with your image
        \caption{09.04.2024}
        \label{fig:wind_090424}
    \end{subfigure}
    \caption{Normalized wind power forecast and real production data in Germany in the region operated by \textit{50Hertz} in 2024. Mean path and 95\% confidence intervals of the SDE in~\eqref{eqn:re_dynamics} calibrated from 2023 data.}
    \label{fig:calibrated_re_dynamics}
\end{figure}

The calibrated values for the above data are $\alpha = 0.34$, $\theta_0 = 2.3948$, and $\varsigma=0.054$. Given the deterministic forecasts of wind power from the same operator in 2024, we can generate sample paths of the SDE in~\eqref{eqn:re_dynamics} with parameters $\alpha,\theta_0,\varsigma$ calibrated from 2023 data. Figure~\ref{fig:calibrated_re_dynamics} illustrates the results for two sample days in 2024.

\begin{figure}[h!]
	\centering
	\includegraphics[width=\textwidth]{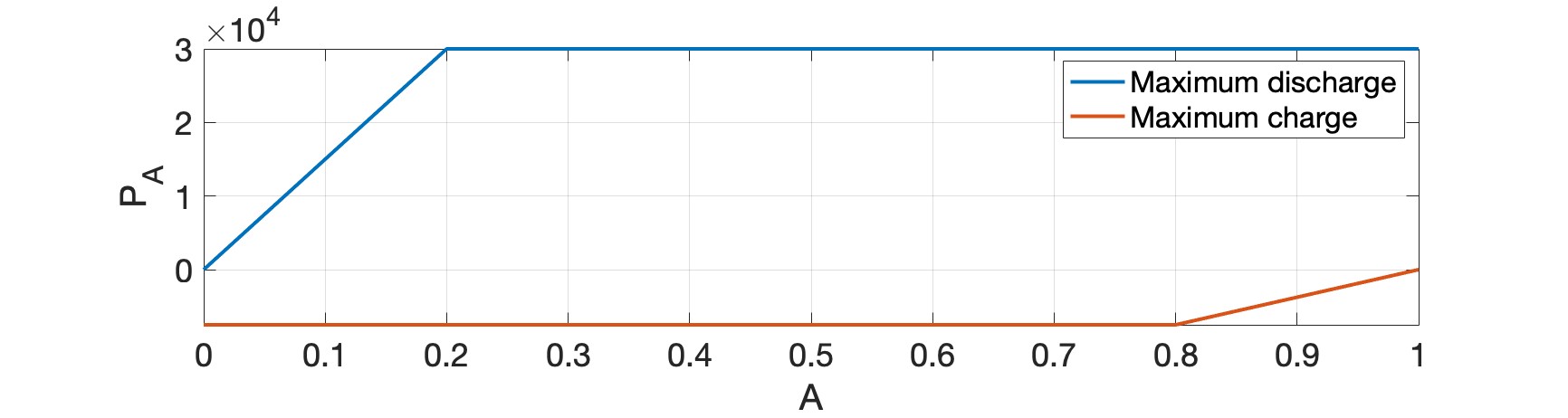}
	\caption{Simple characteristic curve of a battery, described by Eqs.~\eqref{eqn:battery_dynamics} and~\eqref{eqn:characteristic}.}
	\label{fig:battery_profile}
\end{figure}

This work applies~\eqref{eqn:battery_dynamics} to model the instantaneous normalized battery charge. As given in~\eqref{eqn:battery_dynamics}, the maximum charging capacity $\underline{P}_A$ and the maximum discharge capacity $\bar{P}_A$ are defined as functions of the stored battery charge. This function is typically called the characteristic curve of the battery. The numerical experiment employs a simple characteristic curve illustrated in Figure~\ref{fig:battery_profile}. The characteristic curve imposes~\eqref{eqn:characteristic}, implying that the battery cannot discharge power when it has zero charge and can no longer charge itself when fully charged. 
 
\begin{figure}[h!]
	\centering
	\includegraphics[width=0.6\textwidth]{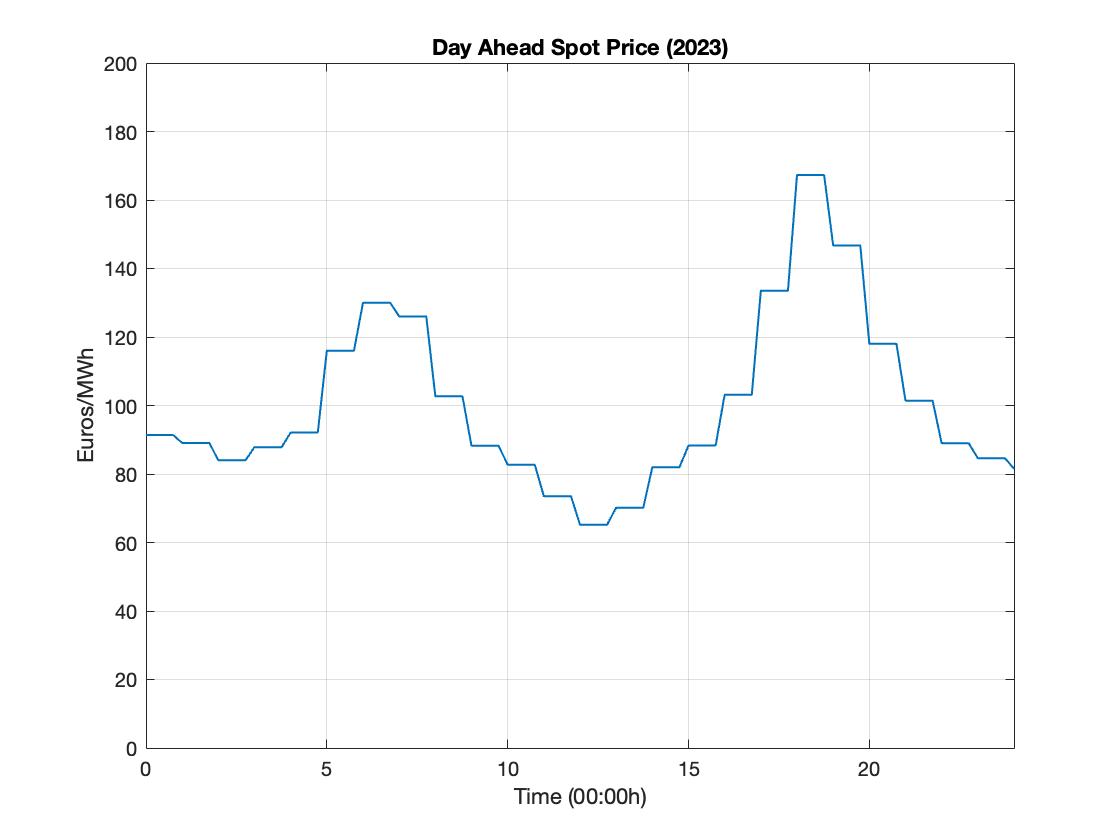}
	\caption{Day-ahead spot-price for grid power in Germany on 07.09.2023. Data are from~\citep{smardMarketData}.}
	\label{fig:spotprice_profile}
\end{figure}

Forecasts of buying and selling prices, $K_b(t)$ and $K_s(t)$, are constructed from publicly available day-ahead spot prices of grid power. We set $K_b(t) = K_s(t)$, implying that the buying and selling prices of power are the same. This assumption is valid because the market commonly ensures that this holds in real life. The numerical experiment employs the German day-ahead spot-price data from 2024 available in~\citep{smardMarketData}. However, these data are generally available only at a 1-hour frequency. Hence, we smoothed and extended the data to domain $[0,T]$ by fitting a polynomial curve to them and using the polynomial function as the input $K_b(t)$. Figure~\ref{fig:spotprice_profile} illustrates an example of the spot-price data.

\subsection{Numerical Results}
\label{sec:numerics}

We present the results of solving the stochastic optimal control problem for the model cellular base station (Section~\ref{sec:model_description}) using the numerical approach in Section~\ref{sec:numerical_methods}.

\subsubsection{Cost-to-go Function of the Dual Problem}
\label{sec:value_fxn}

The numerical solution of the HJB PDE in~\eqref{eqn:hjb_pde} for a given Lagrange multiplier function $\lambda^\ell(t)$~\eqref{eqn:finite_dim_multiplier} is obtained using the scheme in Section~\ref{sec:hjb_numerics}. We apply a uniform rectangular spatial grid with $N_1 = N_2 = N_3 = N_x = 10$, where $N_t$ is set as the lowest value that satisfies the CFL condition~\eqref{eqn:cfl}. In this case, we obtain $N_t = 800$. The presence of second-order derivatives in the HJB PDE in~\eqref{eqn:hjb_pde} yields $N_t = \order{N_x^2}$ for a numerically stable solution, explaining the high value obtained for $N_t$. These discretization parameter choices  ensure that the relative numerical PDE discretization error is below $1\%$. We set $\bar{\chi}$ as the upper $95\%$ quantile of the invariant distribution of process $\xi(t)$. 

\begin{figure}[h!]
    \centering
    % First subfigure
    \begin{subfigure}[b]{0.45\textwidth}
        \centering
        \includegraphics[width=\textwidth]{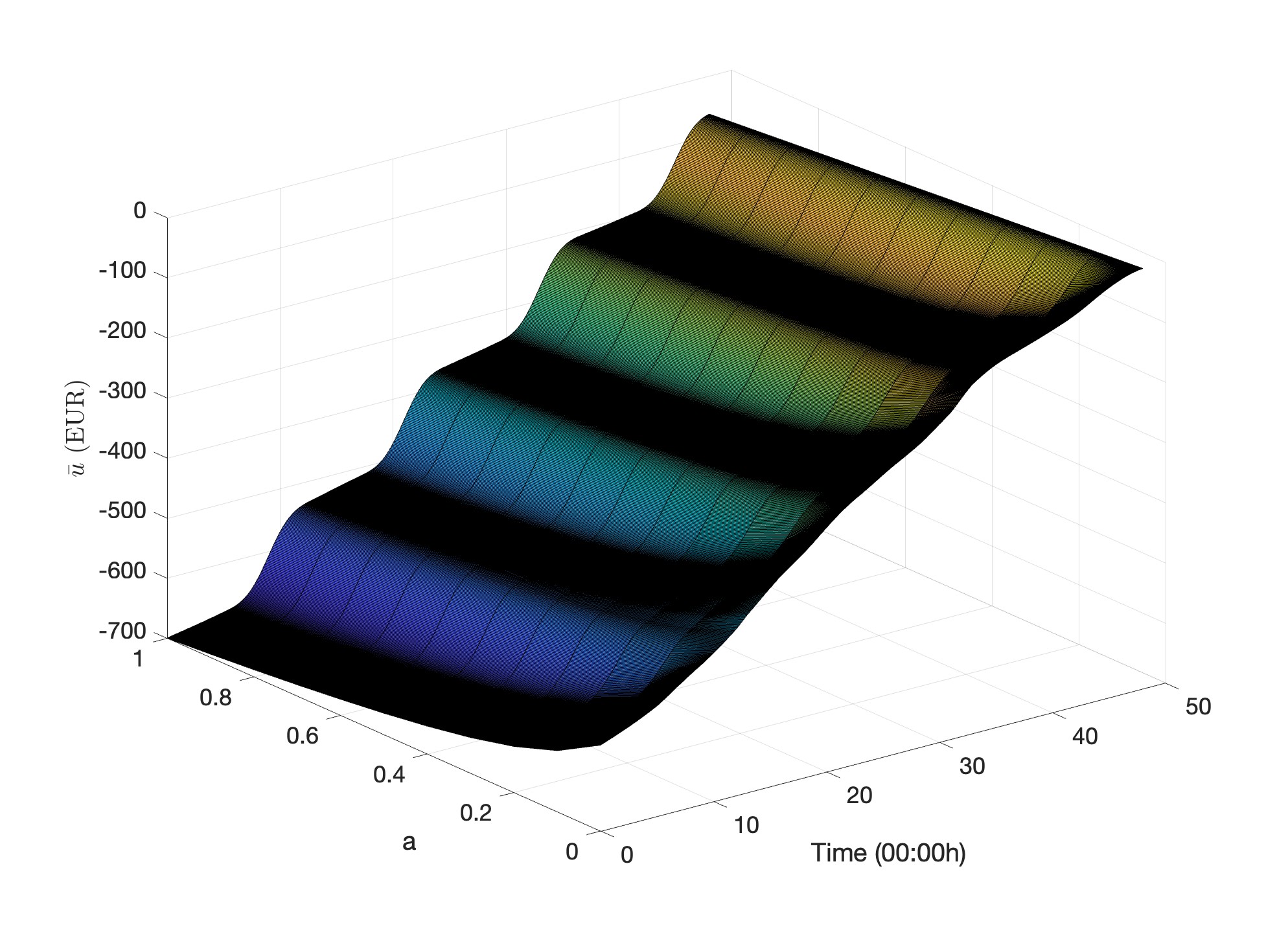}  % Replace with your image
        \caption{$\bar{u}(t,[a,0,0])$}
        \label{fig:u_a_slice}
    \end{subfigure}
    \hfill
    % Second subfigure
    \begin{subfigure}[b]{0.45\textwidth}
        \centering
        \includegraphics[width=\textwidth]{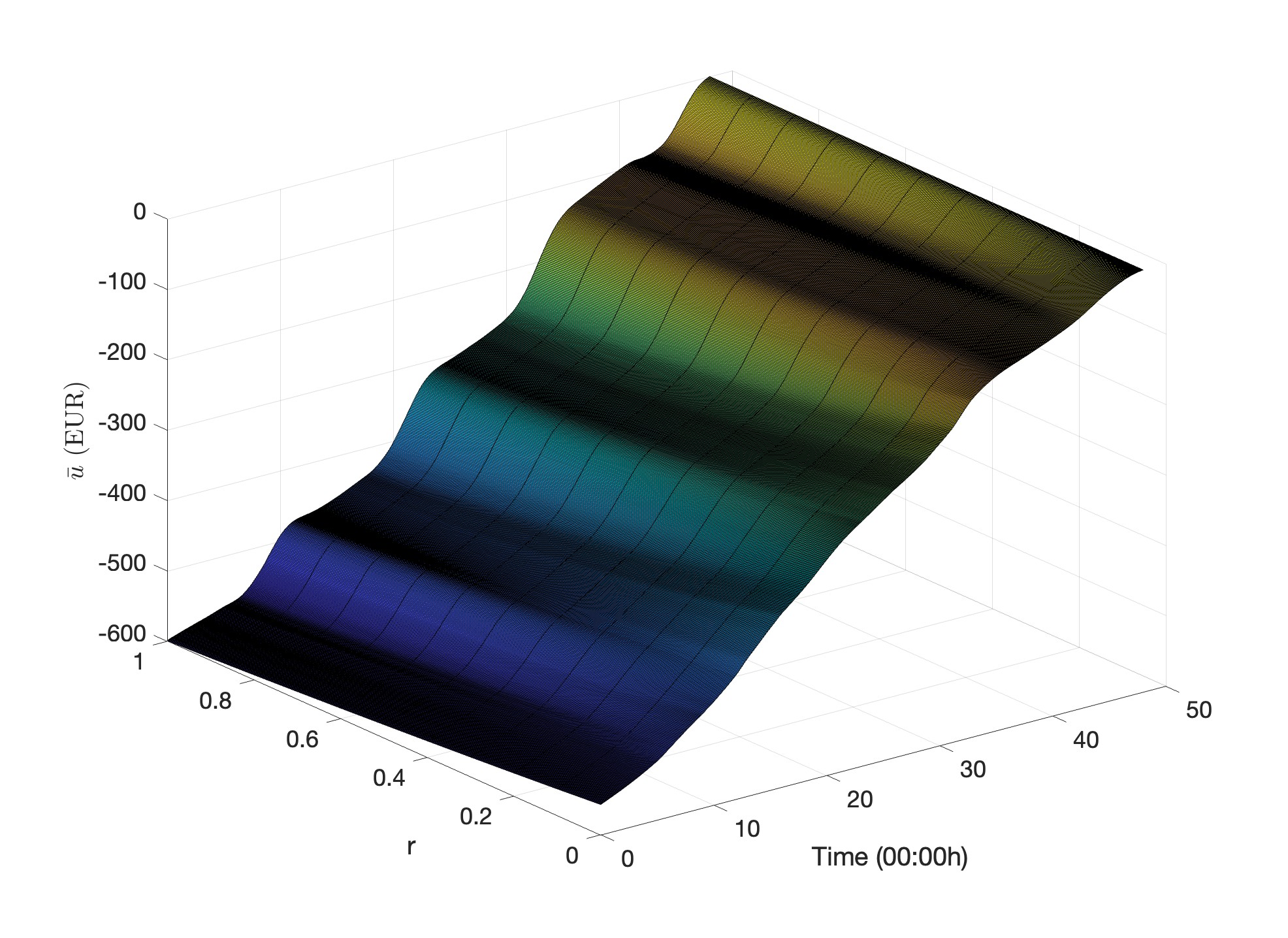}  % Replace with your image
        \caption{$\bar{u}(t,[0,r,0])$}
        \label{fig:u_r_slice}
    \end{subfigure}
    \hfill
    \vspace{1cm}
    % Third subfigure
    \begin{subfigure}[b]{0.45\textwidth}
        \centering
        \includegraphics[width=\textwidth]{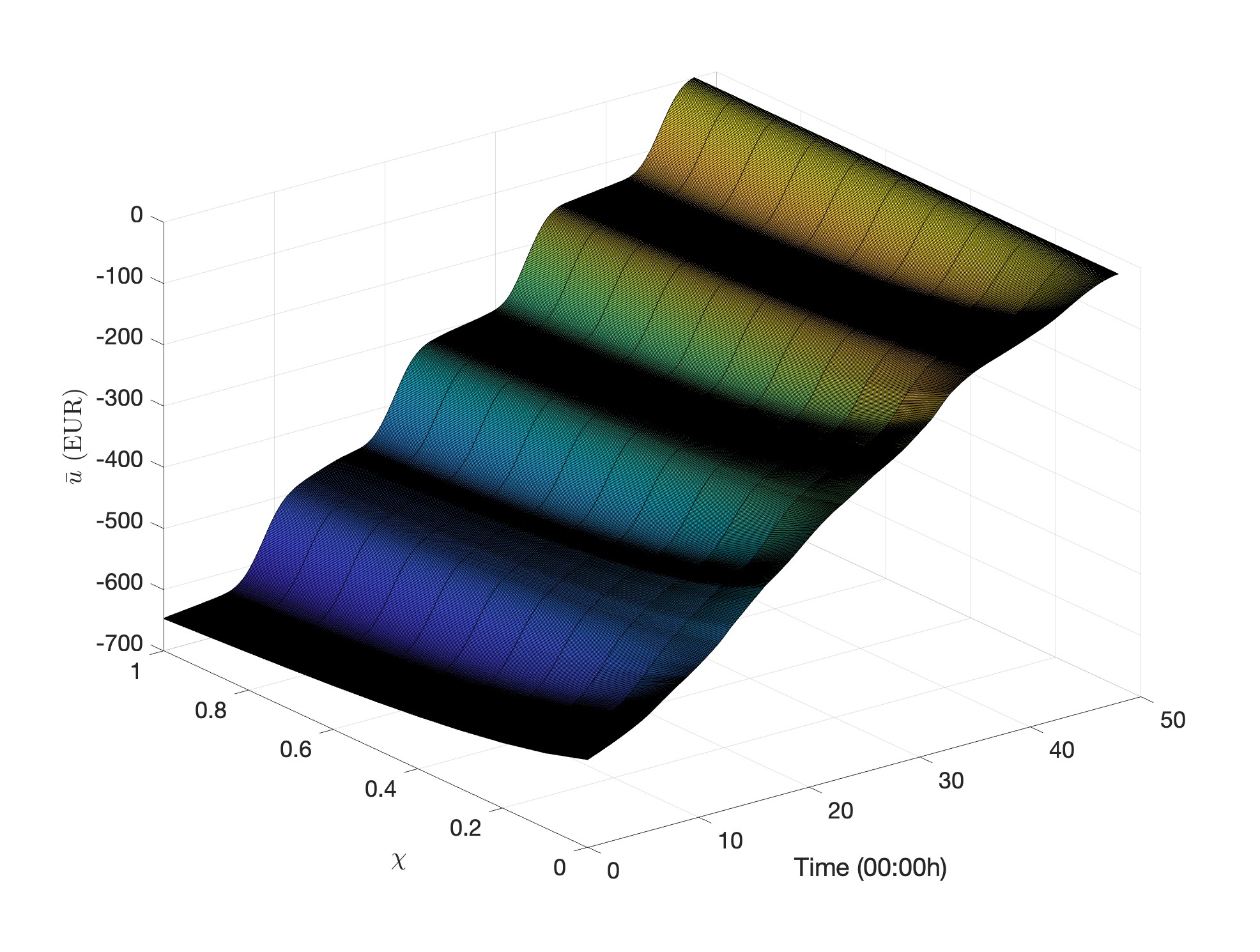}  % Replace with your image
        \caption{$\bar{u}(t,[0,0,\chi])$}
        \label{fig:u_xi_slice}
    \end{subfigure}
    \caption{Slices of the numerical cost-to-go function $\bar{u}$ on the axes of time and state variables. For all slices on each component of $\mathbf{x}$, the rest of the state variables are fixed at $0$.}
    \label{fig:cost_to_go_visual}
\end{figure}

The optimization problem in~\eqref{eqn:numerical_controls} must be solved at each point in the grid $\tau$. For the simple cases of uniform or Gaussian cellular user distribution, we obtain closed-form expressions for $\phi_\text{out}$ (Appendix~\ref{appendix:outage_prop}). The resulting $\phi_\text{out}$ ensures that~\eqref{eqn:numerical_controls} is a convex optimization problem and can be solved analytically. Appendix~\ref{appendix:hamiltonian_soln} presents additional details. Figure~\ref{fig:cost_to_go_visual} presents three slices of the numerical cost-to-go function $\bar{u}$ with the optimal Lagrange multiplier function $\lambda^\ell(t)$, each with two fixed state variables. The numerical dual function approximation $\bar{\Theta}(\lambda^\ell)$ is obtained by evaluating $\bar{u}(0,\mathbf{X}(0))$. 

\subsubsection{Obtaining Reference Costs for Comparison}
\label{sec:reference}

This work compares the values for the dual cost above with reference values obtained from two related, easy-to-solve stochastic optimal control problems. This comparison was conducted to ensure that the output of the numerical algorithm makes sense with respect to these reference values. First, we solve the same optimization problem as in Problem~\ref{prob:primal}, but with the probabilistic constraint replaced with a deterministic a.s. constraint, solving the following problem.

\prblm{ \textit{ (Primal problem with the a.s. QoS constraint)}.}
\label{prob:primal_almost_sure}
Given the initial data $\mathbf{X}(0) = [A_0, R_0 \sim \mu_0^R, \xi_0 \sim \mu_0^\xi]$, and deterministic forecasts for the daily cellular user traffic profile $N_u(t)$, energy spot prices $K_b(t), K_s(t)$, renewable energy forecast $p(t)$, and mobile network usage price $K_\text{net}(t)$ for $0 \leq t \leq T$, we solve the following:

\begin{equation}
\label{eqn:primal_almost_sure}
	\boldsymbol\phi^* = \argmin_{\{\boldsymbol\phi \in \bar{\mathcal{A}}_1\}} \mathcal{U}(\boldsymbol\phi) = w \left( \mathcal{C}_1 - \mathcal{R}_1 - \mathcal{R}_2 - \mathcal{R}_3 \right) + (1-w) \mathcal{C}_2,
\end{equation}

where the minimization is subject to dynamics \eqref{eqn:gamma_fading}, \eqref{eqn:re_dynamics}, and \eqref{eqn:battery_dynamics}, and is done over all controls that satisfy the constraint in~\ref{constraint:as} and the following a.s. constraint (instead of the  constraint in~\ref{constraint:prob}).

\begin{equation}
\label{eqn:qos_constraint_almost_sure}
	\phi_\mathrm{out}(t,\xi(t),P_\mathrm{tx}(t)) = \phi_\mathrm{th}, \quad \forall t \in [0,T],
\end{equation}
where \eqref{eqn:qos_constraint_almost_sure} provides a solution for the control $P_\mathrm{tx}(t)$ for all $t \in [0,T]$. 

Hence, the control vector for this problem has only three components, 
$\boldsymbol\phi(t) = [P_A(t), \allowbreak P_F(t),P_S(t)]$. Problem~\ref{prob:primal_almost_sure} is a standard continuous-time stochastic optimal control problem; hence the dynamic programming principle can be used to derive and numerically solve the associated HJB PDE, analogous to Section~\ref{sec:hjb}. The cost-to-go function for this problem is expected to be higher than the dual cost of Problem~\ref{prob:finite_dim_dual} because the network operator would need to buy more power to keep $\phi_\mathrm{out}$ constant at all times, whereas the operator is allowed to relax the constraint probabilistically in Problem~\ref{prob:finite_dim_dual} as described in the constraint in~\ref{constraint:prob}. Next, we solve the same optimization problem as Problem~\ref{prob:primal}, but without the constraint in~\ref{constraint:prob}.

\prblm{ \textit{ (Primal problem with no QoS constraint)}.}
\label{prob:primal_no_qos}
Given the initial data $\mathbf{X}(0) = [A_0, R_0 \sim \mu_0^R, \xi_0 \sim \mu_0^\xi]$, and deterministic forecasts for the daily cellular user traffic profile $N_u(t)$, energy spot prices $K_b(t), K_s(t)$, renewable energy forecast $p(t)$, and mobile network usage price $K_\text{net}(t)$ for $0 \leq t \leq T$, we solve the following:

\begin{equation}
\label{eqn:primal_no_qos}
	\boldsymbol\phi^* = \argmin_{\{\boldsymbol\phi \in \bar{\mathcal{A}}^\text{rel}\}} \mathcal{U}(\boldsymbol\phi) = w \left( \mathcal{C}_1 - \mathcal{R}_1 - \mathcal{R}_2 - \mathcal{R}_3 \right) + (1-w) \mathcal{C}_2,
\end{equation}

where the minimization is taken over all controls satisfying the constraint in~\ref{constraint:as} with dynamics \eqref{eqn:gamma_fading}, \eqref{eqn:re_dynamics}, and \eqref{eqn:battery_dynamics}. 

Problem~\ref{prob:primal_no_qos} is also a standard continuous-time stochastic optimal control problem; hence the dynamic programming principle can be applied to derive and numerically solve the associated HJB PDE. Without the Lagrangian relaxation, the cost-to-go function associated with this problem is the same as in~\eqref{eqn:value_function}. Thus, this problem can be solved using the numerical scheme in Section~\ref{sec:hjb_numerics}, setting $\lambda(t)=0$. The cost-to-go function for this problem is expected to be less than the dual cost of Problem~\ref{prob:finite_dim_dual} because the network operator only maximizes revenue without ensuring that QoS is achieved. The HJB PDEs corresponding to  Problems~\ref{prob:primal_almost_sure} and~\ref{prob:primal_no_qos} must be numerically solved only once because no associated dual problem needs to be solved.

\subsubsection{Dual Optimization Results}
\label{sec:dual_opt_results}

\begin{figure}[h!]
    \centering
    % First subfigure
    \begin{subfigure}[b]{0.45\textwidth}
        \centering
        \includegraphics[width=\textwidth]{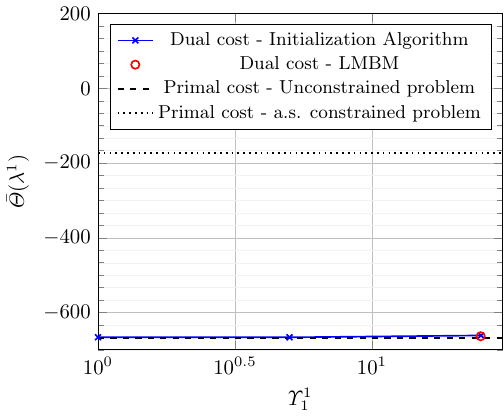}  % Replace with your image
        \caption{Dual cost $\bar{\Theta}(\lambda^1)$}
        \label{fig:init_lmbm+ssm_dual}
    \end{subfigure}
    \hfill
    % Second subfigure
    \begin{subfigure}[b]{0.45\textwidth}
        \centering
        \includegraphics[width=\textwidth]{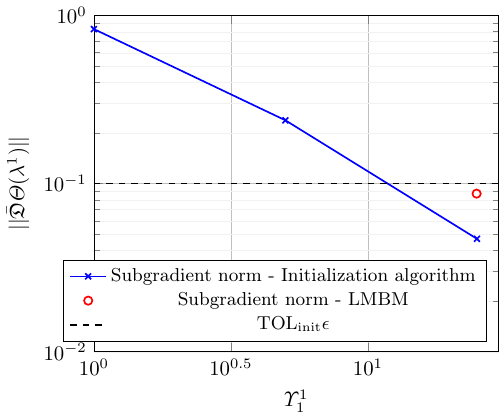}  % Replace with your image
        \caption{Subgradient norm $\norm{\bar{\mathfrak{D}}\Theta(\lambda^1)}$.}
        \label{fig:init_lmbm+ssm_subg}
    \end{subfigure}
    \caption{Dual cost and subgradient norm computed with $\lambda^1$ constructed from $\Upsilon_1^1$ attained at each iteration of Algorithm~\ref{alg:initialization}, along with the corresponding quantities using the final output of LMBM routine~\citep{Karmitsa:2007aa}.}
    \label{fig:init_lmbm+ssm}
\end{figure}

This work runs Algorithm~\ref{alg:dual_optimization} with parameters set to values in Appendix~\ref{appendix:sim_parameters}. Notably,  we set $\phi_\mathrm{th} = 10^{-3}$ and $\epsilon = 0.1$, implying that the network operator must ensure that the proportion of users in outage $\phi_\text{out}$ should be less than $10^{-3}$ with a probability of $90\%$ for all $t \in [0,T]$. Figure~\ref{fig:init_lmbm+ssm} presents the results of the initialization Algorithm~\ref{alg:initialization} and the LMBM routine. Figure~\ref{fig:init_lmbm+ssm_dual} reveals the dual cost versus $\Upsilon_1^1$ at each iteration of Algorithm~\ref{alg:initialization}. The dual cost increases with the number of iterations, and it is within the two reference values attained by solving Problems~\ref{prob:primal_almost_sure} and~\ref{prob:primal_no_qos}. Figure~\ref{fig:init_lmbm+ssm_subg} plots the corresponding estimated subgradient values. The initialization algorithm reaches the prescribed tolerance in three iterations. Figure~\ref{fig:init_lmbm+ssm} reveals that the LMBM offers a more optimal point than the output of Algorithm~\ref{alg:initialization},  providing an excellent starting point for the dual optimization in Algorithm~\ref{alg:dual_optimization}.

\begin{figure}[h!]
    \centering
    % First subfigure
    \begin{subfigure}[b]{0.45\textwidth}
        \centering
        \includegraphics[width=\textwidth]{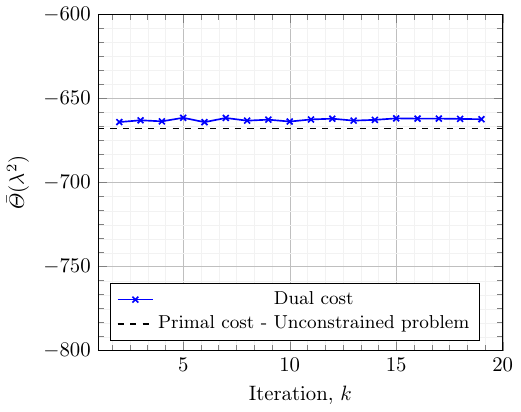}  % Replace with your image
        \caption{Dual cost $\bar{\Theta}(\lambda^2)$}
        \label{fig:main_ell8_dual}
    \end{subfigure}
    \hfill
    % Second subfigure
    \begin{subfigure}[b]{0.45\textwidth}
        \centering
        \includegraphics[width=\textwidth]{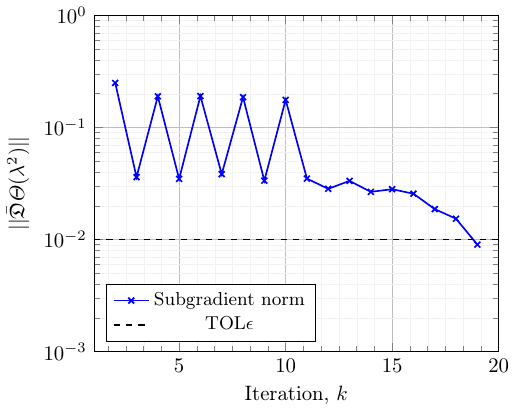}  % Replace with your image
        \caption{Subgradient norm $\norm{\bar{\mathfrak{D}}\Theta(\lambda^2)}$.}
        \label{fig:main_ell8_subg}
    \end{subfigure}
    \caption{Dual cost and subgradient norm computed with $\lambda^2$ constructed from $\Upsilon^2$ at each SSM iteration of Algorithm~\ref{alg:dual_optimization} for $\ell=2$.}
    \label{fig:main_ell8}
\end{figure}

Figure~\ref{fig:main_ell8} illustrates the evolution of the dual function values and its corresponding subgradient at each iteration of the SSM for a fixed $\ell$. Figure~\ref{fig:main_ell8_dual} verifies that not every step of the SSM is taken towards maximization. This behavior illustrates the difficulty of solving the nonsmooth optimization in Problem~\ref{prob:finite_dim_dual}. Due to the satisfactory  choice of starting values and step sizes (Section~\ref{sec:solve_dual}), Algorithm~\ref{alg:dual_optimization} converges for each $\ell$ within \texttt{max-iter} iterations with respect to the subgradient.

\begin{figure}[h!]
    \centering
    % First subfigure
    \begin{subfigure}[b]{0.45\textwidth}
        \centering
        \includegraphics[width=\textwidth]{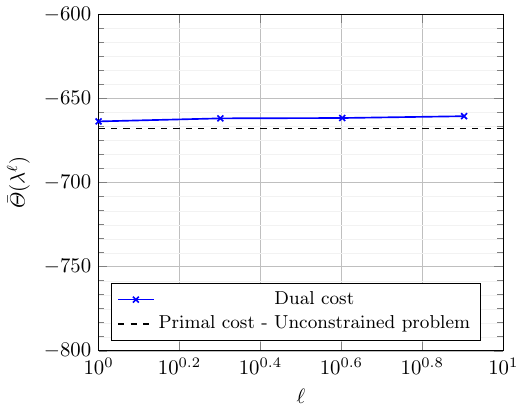}  % Replace with your image
        \caption{Dual cost $\bar{\Theta}(\lambda^\ell)$}
        \label{fig:adaptive_dual}
    \end{subfigure}
    \hfill
    % Second subfigure
    \begin{subfigure}[b]{0.45\textwidth}
        \centering
        \includegraphics[width=\textwidth]{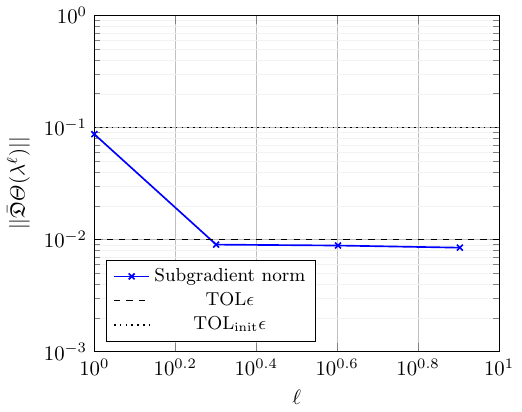}  % Replace with your image
        \caption{Subgradient norm $\norm{\bar{\mathfrak{D}}\Theta(\lambda^\ell)}$.}
        \label{fig:adaptive_subg}
    \end{subfigure}
    \caption{Dual cost and subgradient norm computed with $\lambda^\ell$ constructed from the optimal $\Upsilon^\ell$ obtained from  the SSM at each refinement level $\ell$ of Algorithm~\ref{alg:dual_optimization}.}
    \label{fig:adaptive}
\end{figure}

Figure~\ref{fig:adaptive} plots the evolving optimal dual function $\bar{\Theta}(\lambda^\ell)$ and its corresponding subgradient $\bar{\mathfrak{D}}\Theta(\lambda^\ell)$ constructed using the optimal amplitudes $\Upsilon^\ell$ at each level $\ell$ of Algorithm~\ref{alg:dual_optimization}. Table~\ref{tab:dual_cost} lists the optimal dual cost at each stage of Algorithm~\ref{alg:dual_optimization}. Table~\ref{tab:dual_cost} presents minor increases in the dual function evaluations as the refinement level $\ell$ increases, implying that we successively obtain more optimal solutions as the Lagrange multiplier function $\lambda^\ell(t)$ becomes finer. We normalize the subgradient vector norm by the number of dimensions for a fair comparison between the  norm of vectors of varying dimensions. 

Table~\ref{tab:dual_cost} also highlights the advantage of applying a probabilistic constraint to achieve the QoS in cellular networks, instead of doing it a.s. for all $t \in [0,T]$ (as in Problem~\ref{prob:primal_almost_sure}). The optimal dual cost achieved by Algorithm~\ref{alg:dual_optimization} is much closer to the expected optimal cost of the unconstrained Problem~\ref{prob:primal_no_qos} and around four times less expensive than satisfying~\eqref{eqn:qos_constraint_almost_sure}, a.s. (Problem~\ref{prob:primal_almost_sure}). In rare occasions when an unreasonable surge in demand ($N_u(t)$) occurs, or a very low wireless fading channel ($\xi(t)$) exists due to environmental conditions or low availability of renewable power ($P_R(t)$), the numerical approach allows the network operator to relax the stringent QoS constraint to save substantial expenditure while still serving as many customers as possible. This approach contrasts the solution to Problem~\ref{prob:primal_almost_sure}, which enforces~\eqref{eqn:qos_constraint_almost_sure} even in such periods, significantly increasing operating expenditure.

\begin{table}[h!]
	\centering
	 \begin{tabular}{||c|c|c||} 
	 \hline
	 Description & $\ell$ & Optimal cost (\texteuro) \\ [0.5ex] 
	 \hline\hline
	 Problem~\ref{prob:primal_almost_sure} & & $-173.57$ \\
	 Problem~\ref{prob:primal_no_qos} & & $-667.83$ \\
	 Initialization Algorithm~\ref{alg:initialization} & $1$ & $-661.52$ \\
	 LMBM & $1$ & $-663.62$ \\
	 Dual optimization Algorithm~\ref{alg:dual_optimization} & $2$ & $-661.8$ \\
	 Dual optimization Algorithm~\ref{alg:dual_optimization} & $4$ & $-661.59$ \\
	 Dual optimization Algorithm~\ref{alg:dual_optimization} & $8$ & $-660.51$ \\
	 \hline
	 \end{tabular}
	 \caption{Comparing the optimal dual costs achieved in various stages of Algorithm~\ref{alg:dual_optimization} with those of  Problems~\ref{prob:primal_almost_sure} and~\ref{prob:primal_no_qos}.}
	 \label{tab:dual_cost}
\end{table}

\begin{figure}[h!]
	\centering
	\includegraphics[width=0.6\textwidth]{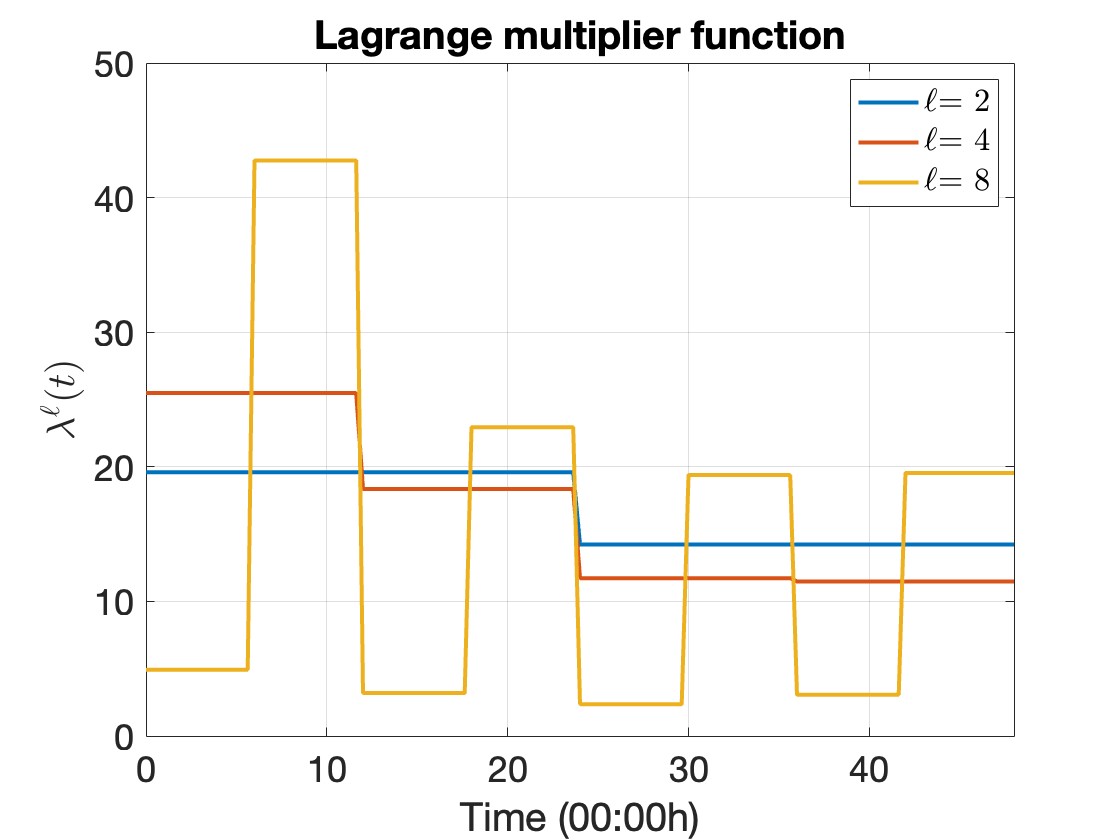}
	\caption{Evolution of the optimal Lagrange multiplier function $\lambda^\ell(t)$ at each refinement level $\ell$ in Algorithm~\ref{alg:dual_optimization}.}
	\label{fig:adaptive_lambda}
\end{figure}

Figure~\ref{fig:adaptive_lambda} presents the evolution of the optimal Lagrange multiplier function as the refinement level $\ell$ increases in Algorithm~\ref{alg:dual_optimization}. For the given example, Algorithm~\ref{alg:dual_optimization} reaches prescribed $\mathrm{TOL}$ at $\ell=8$, and does not refine the Lagrange multiplier function further. High-demand periods (Figure~\ref{fig:user_profile}) correspond to high Lagrange multiplier function values. This behaviour is expected because a high $\lambda^\ell$ forces the mobile operator to satisfy the constraint in~\ref{constraint:prob} in times of high demand, although it significantly increases costs. The oscillation frequency of the base station power demand dictates the oscillations in the Lagrange multiplier function.

\subsubsection{Optimal Power Procurement Policy}
\label{sec:optimal_policy}

\begin{figure}[H]
    \centering
    % First subfigure
    \begin{subfigure}[b]{0.45\textwidth}
        \centering
        \includegraphics[width=\textwidth]{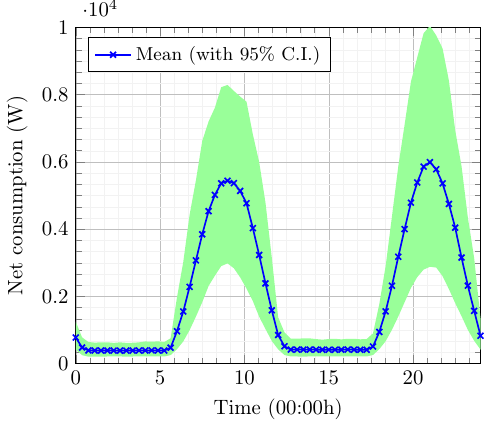}  % Replace with your image
        \caption{$C_\text{scal} N_u(t) P_\mathrm{tx}(t) + C_\text{offset}$}
        \label{fig:opt_consumption}
    \end{subfigure}
    \hfill
    % First subfigure
    \begin{subfigure}[b]{0.45\textwidth}
        \centering
        \includegraphics[width=\textwidth]{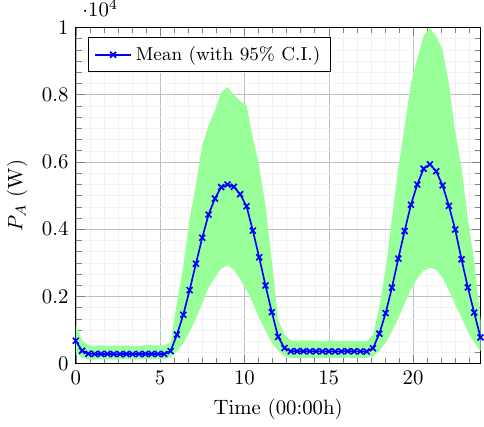}  % Replace with your image
        \caption{$P_A(t)$}
        \label{fig:opt_PA}
    \end{subfigure}
    \hfill
    \vspace{1cm}
    % First subfigure
    \begin{subfigure}[b]{0.45\textwidth}
        \centering
        \includegraphics[width=\textwidth]{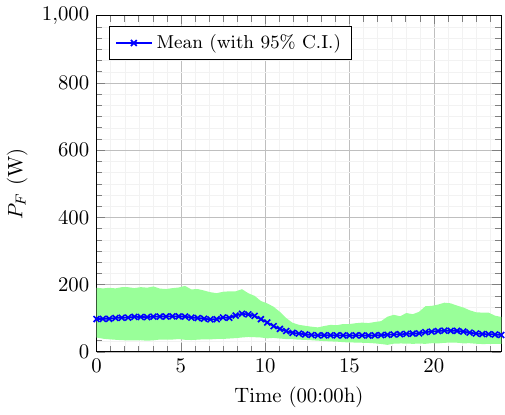}  % Replace with your image
        \caption{$P_F(t)$}
        \label{fig:opt_PF}
    \end{subfigure}
    \hfill
    % First subfigure
    \begin{subfigure}[b]{0.45\textwidth}
        \centering
        \includegraphics[width=\textwidth]{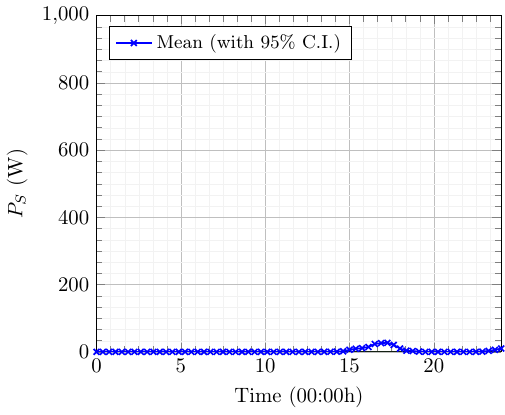}  % Replace with your image
        \caption{$P_S(t)$}
        \label{fig:opt_PS}
    \end{subfigure}
    \hfill
    \vspace{1cm}
    % Second subfigure
    \begin{subfigure}[b]{0.45\textwidth}
        \centering
        \includegraphics[width=\textwidth]{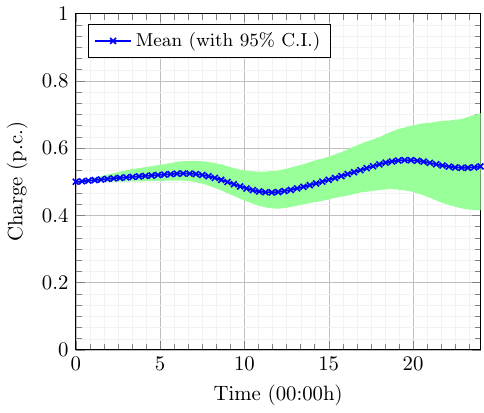}  % Replace with your image
        \caption{$A(t)$}
        \label{fig:opt_charge}
    \end{subfigure}
    \caption{Optimal power procurement policy, with 95\% confidence intervals, using the optimal cost-to-go function $\bar{u}$ from Algorithm~\ref{alg:dual_optimization}.}
    \label{fig:opt_policy}
\end{figure}

With the optimal Lagrange multiplier function $\lambda^\ell$ from Algorithm~\ref{alg:dual_optimization}, the dual cost is $-660.51$\texteuro      ~(Table~\ref{tab:dual_cost}). Figure~\ref{fig:opt_policy} illustrates the optimal power procurement policy for the model cellular base station for the next 24 hours using the optimal solution, including the (a) net power consumption of the base station $C_\text{scal} N_u(t) P_\mathrm{tx}(t) + C_\text{offset}$ (refer \eqref{eqn:power_balance_constraint}), (b) power drawn from the battery to run the base station $P_A(t)$, (c) power bought from the grid $P_F(t)$, (d) power sold back to the grid $P_S(t)$, and (e) optimal charge stored in battery $A(t)$. Figure~\ref{fig:opt_consumption} indicates that the base station consumes high power during periods of high demand (refer Figure~\ref{fig:user_profile}). Figures~\ref{fig:opt_PA}, \ref{fig:opt_PF}, and~\ref{fig:opt_PS} reveal that the battery provides most of the power to run the base station throughout the day. The clean, renewable energy stored in the battery powers a large proportion of base station operations so that the mobile operator does not need to buy expensive grid power. This approach reduces costs and the carbon footprint of the base station. The operator must only buy some power from the grid to satisfy high demand during peak hours (Figure~\ref{fig:opt_PF}). The base station sells a negligible amount of power back to the grid (Figure~\ref{fig:opt_PS}). The optimal policy also advises the operator to store incoming renewable power and charge the battery during low-demand hours before using that power to serve customers during high-demand hours (Figure~\ref{fig:opt_charge}). Figure~\ref{fig:opt_charge} also indicates that the battery has some stored energy at the end of the day that could be used the next day. This observation establishes the effect of implementing the running horizon framework (see Section~\ref{sec:running_horizon}). 

\begin{figure}[H]
    \centering
    % First subfigure
    \begin{subfigure}[b]{0.45\textwidth}
        \centering
        \includegraphics[width=\textwidth]{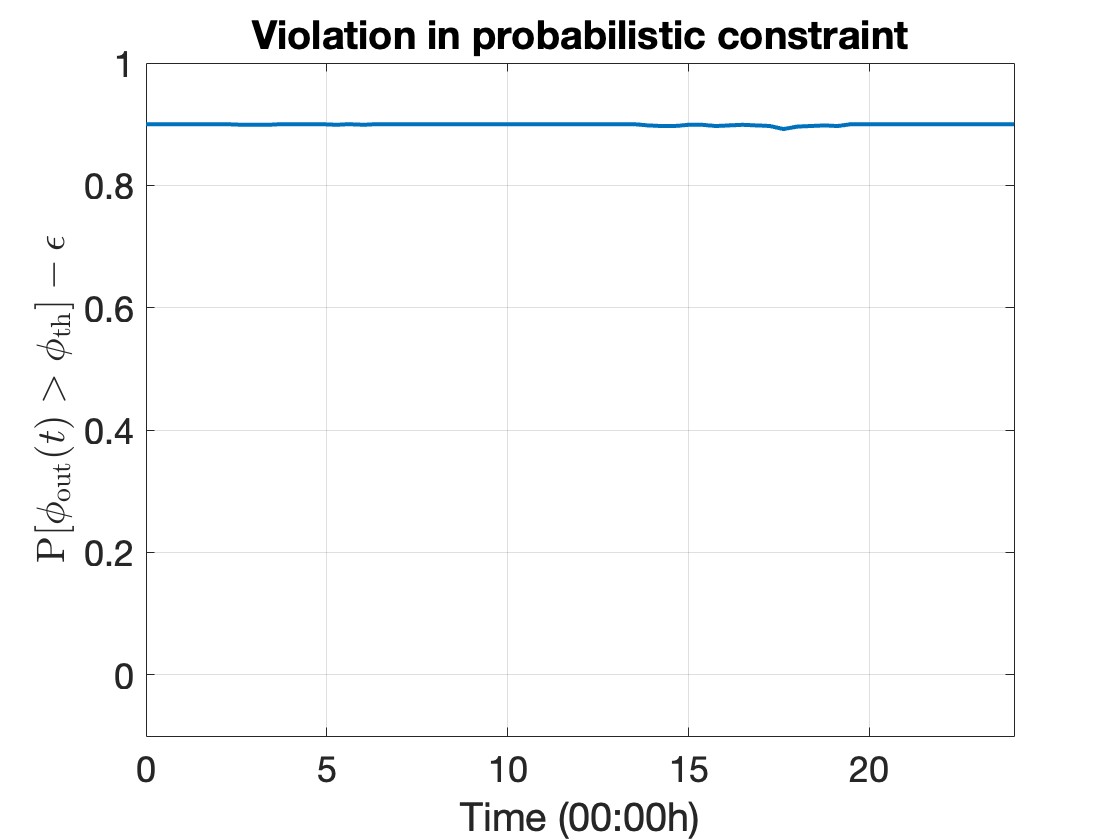}  % Replace with your image
        \caption{Solution to Problem~\ref{prob:primal_no_qos}}
        \label{fig:mcprob_noconstraint}
    \end{subfigure}
    \hfill
    % Second subfigure
    \begin{subfigure}[b]{0.45\textwidth}
        \centering
        \includegraphics[width=\textwidth]{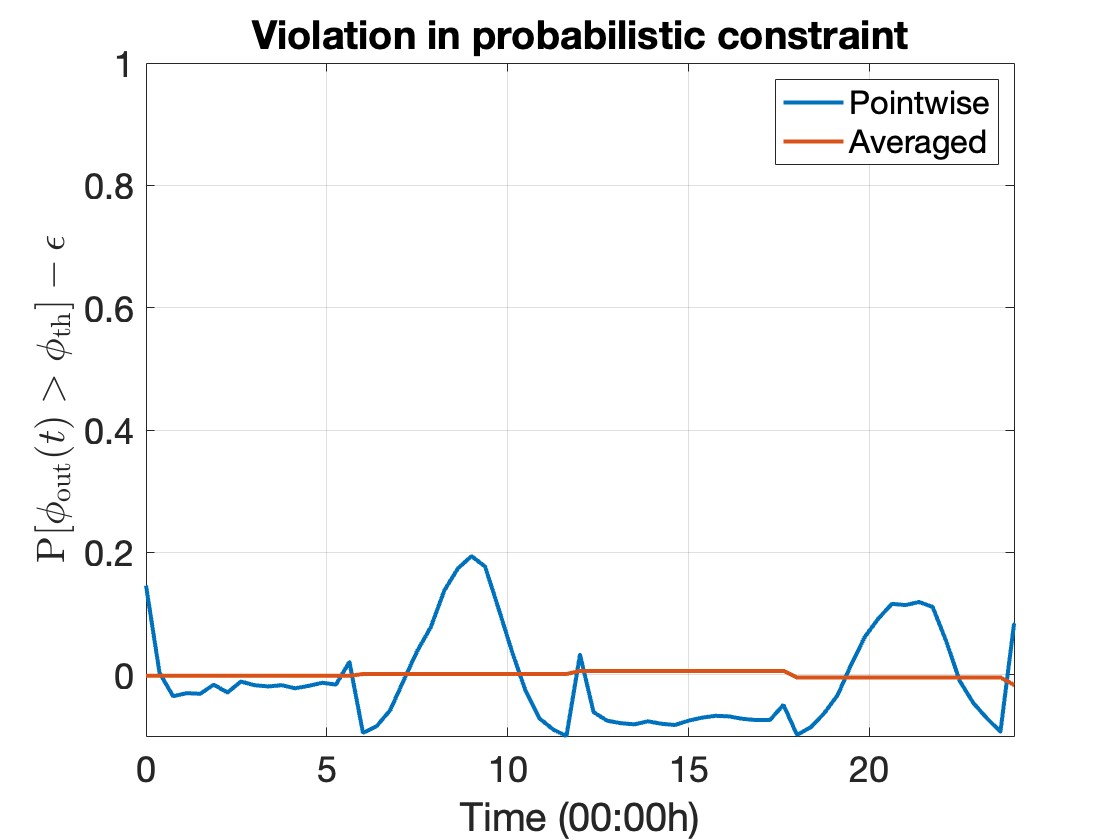}  % Replace with your image
        \caption{Solution to Problem~\ref{prob:finite_dim_dual}.}
        \label{fig:mcprob_adaptive}
    \end{subfigure}
    \caption{Monte Carlo estimates of the violation in constraint $\prob{\phi_\text{out}(t) > \phi_\mathrm{th}} - \epsilon$ for each $t \in [0,T]$ using optimal controls from solving (a) Problem~\ref{prob:primal_no_qos} and (b) Problem~\ref{prob:finite_dim_dual} using Algorithm~\ref{alg:dual_optimization}.}
    \label{fig:mcprob}
\end{figure}

We estimate the time-pointwise violation of the constraint in~\ref{constraint:prob} using the Monte Carlo method and compare it with the final subgradient vector $\bar{\mathfrak{D}} \Theta(\lambda^\ell)$ from Algorithm~\ref{alg:dual_optimization} to validate the optimality of the solution further. We run approximate sample numerical paths of the optimally controlled dynamics, according to~\eqref{eqn:sde_euler}, with the optimal policy above. Figure~\ref{fig:mcprob_adaptive} plots this information for a discrete set of points in $t \in [0,T]$. We also compare this with the optimal policy attained from solving the unconstrained problem in  Problem~\ref{prob:primal_no_qos} in Figure~\ref{fig:mcprob_noconstraint}. Figure~\ref{fig:mcprob_noconstraint} reveals that the user outage proportion $\phi_\text{out}$ for the optimal policy with no constraints on the QoS is almost consistently greater than the threshold $\phi_\mathrm{th}$. In contrast, the proposed optimal solution controls the probabilistic constraint in an integral (averaged) sense (refer~\eqref{eqn:subgradient}), even though there are some oscillations around zero in the time-pointwise violation of the constraint in~\ref{constraint:prob}. Figure~\ref{fig:mcprob_adaptive} illustrates this behavior, with the time-averaged violation in the probabilistic constraint in~\ref{constraint:prob} close to zero up to the prescribed relative tolerance $\tol$.

\subsubsection{Sensitivity Analysis}

This work investigates the performance of the proposed algorithm under certain scenarios that could be of interest to network operators,  which are listed below.

\begin{itemize}
	\item \textbf{Scenario A}: No incoming renewable power for a day ($P_R(t) = 0$ for all $t \in [0,T]$).
	\item \textbf{Scenario B}: The wireless fading channel is substantially low due to extreme weather ($\underline{\xi}$ is halved).
	\item \textbf{Scenario C}: High weight is assigned to minimizing operating expenditure ($w = 0.99$).
	\item \textbf{Scenario D}: The price users pay to connect to the network substantially reduces ($K_\text{net} = 0.001$ \texteuro /h per person).
	\item \textbf{Scenario E}: The probabilistic QoS constraint~\ref{constraint:prob} must be satisfied with low confidence ($\epsilon=0.2$).
\end{itemize}

All parameters and dynamics are set as in Table~\ref{tab:model_parameters} and Section~\ref{sec:model_description}, except for the prescribed parameter change to simulate each scenario. For each scenario, we solve three problems: (i) dual optimization (Problem~\ref{prob:finite_dim_dual}), (ii) Problem~\ref{prob:primal_almost_sure} with a.s. QoS constraint, and  (iii) Problem~\ref{prob:primal_no_qos} with no QoS constraint. Table~\ref{tab:sensitivity_comparison_tab2} reports the optimal costs in each case for all five scenarios.

\begin{table}[h!]
    \centering
    \begin{tabular}{||l||c|c|c||} 
    	\hline \hline
        \textbf{} & \textbf{Dual Algorithm~\ref{alg:dual_optimization}} & \textbf{Problem~\ref{prob:primal_no_qos}} & \textbf{Problem~\ref{prob:primal_almost_sure}} \\ 
        \hline 
        Scenario A  & $-441.8$ & $-564.9$ & $1941.9$ \\ 
        Scenario B  & $-649.28$ & $-660.27$ & $4160.4$ \\ 
        Scenario C  & $-1334.7$ & $-1338.8$ & $1980.5$ \\ 
        Scenario D & $-80.52$ & $-92.69$ & $1209.5$ \\ 
        Scenario E & $-664.03$ & $-667.83$ & $-173.57$ \\ 
        \hline
    \end{tabular}
    \caption{Comparing the optimal dual costs (in \texteuro) achieved with the dual optimization Algorithm~\ref{alg:dual_optimization} for each scenario with the optimal primal costs of the corresponding unconstrained and a.s. QoS constrained problems.}
    \label{tab:sensitivity_comparison_tab2}
\end{table}

In all scenarios, the optimal dual cost lies between the primal costs of the corresponding unconstrained and a.s. QoS constrained problems. Moreover, the optimal dual cost is closer to the primal cost of the unconstrained problem in all cases because it is costly to satisfy the QoS constraint a.s. at all times. The optimal dual cost in Scenario~A is 33\% higher than the optimal dual cost in Table~\ref{tab:dual_cost} because there is no "free" incoming renewable energy for the operator to run the base station, increasing the  expenditure. The optimal dual cost of Scenario~C is 102\% higher than the optimal dual cost in Table~\ref{tab:dual_cost}, directly implying the high weight assigned to financial profits, with no regard for its carbon footprint. The optimal cost of Scenario~E confirms the intuition that it is less expensive to satisfy the QoS with low confidence (80\%), compared to 90\% confidence. To quantify differences in the optimal procurement policy for each scenario, we produced plots analogous to Figure~\ref{fig:opt_policy}, and approximated the area below the mean curve for each plot to obtain the expected energy balance.  Table~\ref{tab:sensitivity_comparison_tab1} lists these values. The base scenario in Table~\ref{tab:sensitivity_comparison_tab1} refers to the original problem described in Section~\ref{sec:model_description}. 

\begin{table}[h!]
    \centering
    \begin{tabular}{||l||c|c|c|c||} 
        \hline \hline 
        \textbf{Expected energy} & \textbf{Consumed} & \textbf{Battery} & \textbf{Bought} &  \textbf{Sold} \\ 
        \hline
        Base scenario & $4199.6$ & $4079.9$ & $119.63$ & $49.64$ \\ 
        Scenario A  & $4040$ & $3295.7$ & $744.36$ & $29.94$ \\ 
        Scenario B  & $4584.3$ & $4440.6$ & $143.68$ & $43.5$ \\ 
        Scenario C  & $4402.8$ & $2347.5$ & $2055.2$ & $5363.6$ \\ 
        Scenario D & $3883.3$ & $3811.6$ & $71.64$ & $50.99$ \\ 
        Scenario E & $3870.8$ & $3763.3$ & $107.48$ & $52.35$ \\ 
        \hline
    \end{tabular}
    \caption{Comparison of the expected energy balance (in Wh) balance using the optimal power procurement policy from the optimal $\bar{u}$ obtained from the dual optimization Algorithm~\ref{alg:dual_optimization} for each scenario.}
    \label{tab:sensitivity_comparison_tab1}
\end{table}

In Scenario~A, the operator buys six times more energy than the base scenario. The lack of incoming renewable power forces the network operator to buy more power from the grid to satisfy the QoS constraint. In Scenario~C, the operator buys almost half the energy consumed by the base station from the grid and sells a considerable amount of energy stored in the battery. The optimal policy is to make as much financial profit as possible because a high weight is assigned to minimizing financial costs. The base station consumes 9\% more energy in Scenario~B than the base scenario because the energy demand of satisfying QoS during periods of low fading channel (extreme weather) is higher. An 8\% decrease occurs in the energy consumed by the base station in Scenario~E compared with the base scenario. Less energy is required to satisfy the QoS constraint with lower confidence. 

\begin{figure}[H]
    \centering
    % First subfigure
    \begin{subfigure}[b]{0.45\textwidth}
        \centering
        \includegraphics[width=\textwidth]{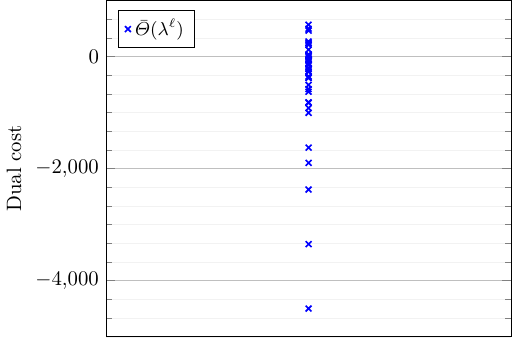}  % Replace with your image
        \caption{Dual cost $\bar{\Theta}(\lambda^\ell)$}
        \label{fig:scenarios_dualcost}
    \end{subfigure}
    \hfill
    % Second subfigure
    \begin{subfigure}[b]{0.45\textwidth}
        \centering
        \includegraphics[width=\textwidth]{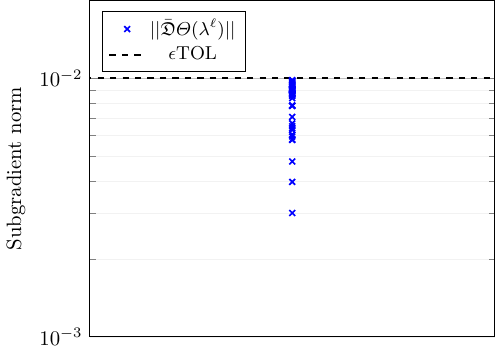}  % Replace with your image
        \caption{Subgradient norm $\norm{\bar{\mathfrak{D}} \Theta (\lambda^\ell)}$.}
        \label{fig:scenarios_subg}
    \end{subfigure}
    \caption{Optimal dual cost and subgradient norms obtained in $50$ independent runs of Algorithm~\ref{alg:dual_optimization} with randomized model and algorithmic parameters as specified in Table~\ref{tab:scenarios}.}
    \label{fig:scenarios}
\end{figure}

This work also performs scenario simulations over probability distributions of both model and algorithmic parameters in Table~\ref{tab:model_parameters} and Table~\ref{tab:sim_parameters}, respectively. Table~\ref{tab:scenarios} lists the various probability distributions used to randomize the simulations. Algorithm~\ref{alg:dual_optimization} is run independently for $50$ i.i.d. samples of the various parameters with no additional tuning between runs. Figure illustrates the optimal dual cost and subgradient achieved in each case. Figure~\ref{fig:scenarios_subg} reveals that Algorithm~\ref{alg:dual_optimization} achieves the prescribed subgradient tolerance in all $50$ cases, demonstrating the robustness of the proposed numerical procedure.  These results emphasize the practicality of applying the proposed numerical approach as a decision-making tool for stakeholders in the wireless communications industry.

\section{Conclusions}

This work formulates a short-term, continuous-time stochastic optimal control problem that minimizes the operating expenditure and carbon footprint of a cellular wireless network powered by a hybrid energy system and equipped with an energy storage capacity under uncertainty in renewable energy and wireless fading channels. We constructed an SDE with the shifted-gamma invariant distribution to model the stochastic Nakagami wireless fading channel. We employed an SDE model~\citep{caballero21_derivative} to model the uncertain renewable power, with parameters calibrated using publicly available German energy market data from 2023. The QoS of the wireless network was posed as a probabilistic constraint to address uncertainties efficiently in the system. We developed a novel Lagrangian relaxation procedure in continuous-time to relax the above probabilistic constraint. We developed a robust numerical procedure to solve the corresponding dual problem based on a finite-dimensional approximation of the Lagrange multiplier function, using an SSM algorithm driven by the numerical solution of the associated HJB PDE. An LMBM-boosted initialization procedure and an adaptive refinement of the Lagrange multiplier function enhanced algorithm performance. A running horizon framework resolved the spurious results at the end of the day due to the short-term horizon problem formulation. The practical application of the numerical procedure was demonstrated on a model cellular base station driven by the German power system and daily cellular traffic demand. The results illustrate the viability and efficiency of the proposed approach for real-world problems. 

A potential future research direction is incorporating stochasticity in other dynamics, such as energy prices $K_b(t), K_s(t)$, and the  cellular user traffic profile $N_u(t)$~\citep{Rached:2018aa}. This research would build a more realistic algorithm that is robust to random fluctuations in these dynamics but would increase the dimension of the associated HJB PDE. Another limitation of the current work is that small values of $\epsilon$ ($\ll 10^{-3}$) cannot be achieved, which modern-day cellular networks aim to achieve, because the naive Monte Carlo method is employed to estimate the subgradient. For a small $\epsilon$, the Fokker--Planck PDE for the joint distribution of the optimally controlled dynamics could be solved to get a deterministic estimate of the subgradient. Alternatively, importance sampling could be applied to reduce the relative statistical error of the Monte Carlo subgradient estimates. A further improvement is devising a better strategy for refining the Lagrange multiplier function $\lambda(t)$. The computed time-pointwise violation in constraint values could be employed to construct a clever strategy to refine $\lambda$ instead of equally spaced refinements.

\appendix
\section[\appendixname~\thesection]{Algorithms}
\label{appendix:alg}

\subsection[\appendixname~\thesubsection]{Upwind Finite-difference Numerical Solver for HJB PDE~\eqref{eqn:hjb_pde}}

\begin{algorithm}[H] 
    \caption{HJB numerical solver}
\label{alg:hjb_solver}
    \SetAlgoLined
    \textbf{Input: } $\lambda(t)$, $N_t$, $N_1$, $N_2$, $N_3$;  \\
    \textbf{Initialization: } $\bar{u}_{(i,j,k)}^{N_t} = -P_K \bar{A} a_i$, $\Delta^\pm_a \bar{u}_{(i,j,k)}^{N_t} = -P_K \bar{A}$ at grid $\tau_a \times \tau_r \times \tau_\chi$ ; \\
    \For{$n=N_t,\ldots,2$}{
    \For{$i=1,\ldots,N_1$}{
    \For{$j=1,\ldots,N_2$}{
    \For{$k=1,\ldots,N_3$}{
    	Using $\left( \Delta_a^\pm \bar{u} \right)_{(i,j,k)}^n$, compute optimal controls $(\boldsymbol\phi^*)_{(i,j,k)}^n$ by solving~\eqref{eqn:numerical_controls}; \\
    	Using $\bar{u}_{(i,j,k)}^{n}$ and $(\boldsymbol\phi^*)_{(i,j,k)}^n$, compute $\bar{u}_{(i,j,k)}^{n-1}$ with update rule \eqref{eqn:linear_update}; \\
    }
    }
    }
    Using $\bar{u}_{(i,j,k)}^{n-1}$, compute and save gradients $\left( \Delta^\pm_a \bar{u } \right)_{(i,j,k)}^{n-1}$ using~\eqref{eqn:numerical_derivative} at grid $\tau_a \times \tau_r \times \tau_\chi$; \\
    }
    \textbf{Output: } $\bar{u}$ at grid $\tau$.
\end{algorithm}

\subsection[\appendixname~\thesubsection]{Euler--Maruyama Monte Carlo for Subgradient~\eqref{eqn:numerical_subgradient} Estimation}

\begin{algorithm}[H] 
    \caption{Numerical subgradient estimation}
\label{alg:estimate_subgradient}
    \SetAlgoLined
    \textbf{Input: } $\ell$, $\bar{N}_t$, $\tilde{N}_t$, $M_\mathrm{SG}$;  \\
    \textbf{Initialization: } $\bar{\mathbf{X}}^{\bar{N}_t}(\bar{t}_0) = \mathbf{X}(0)$; \\
    \For{$m=1,\ldots,M_\mathrm{SG}$}{
    \For{$n=0,\ldots,\bar{N}_t-1$}{
    	Generate $m^\mathrm{th}$ i.i.d. realization of random variable $\boldsymbol\varepsilon^n$; \\
    	Estimate $\Delta^\pm_a \bar{u} (\bar{t}_n,\bar{\mathbf{X}}^{\bar{N}_t}(\bar{t}_n))$ by linearly interpolating from $\left( \Delta^\pm_a \bar{u} \right)$ values computed in Algorithm~\ref{alg:hjb_solver} at grid $\tau$; \\
    	Using $\bar{\mathbf{X}}^{\bar{N}_t}(\bar{t}_n)$ and $\Delta^\pm_a \bar{u} (\bar{t}_n,\bar{\mathbf{X}}^{\bar{N}_t}(\bar{t}_n))$, compute optimal controls $\boldsymbol\phi^*(\bar{t}_n,\bar{\mathbf{X}}^{\bar{N}_t}(\bar{t}_n))$ by solving~\eqref{eqn:numerical_controls}; \\
    	Using $\boldsymbol\varepsilon^n$ and $\boldsymbol\phi^*(\bar{t}_n,\bar{\mathbf{X}}^{\bar{N}_t}(\bar{t}_n))$, compute $\bar{\mathbf{X}}^{\bar{N}_t}(\bar{t}_{n+1})$ using~\eqref{eqn:sde_euler}; \\
    }
    \For{$i=1,\ldots,\ell$}{
    	Construct uniform time grid $\tilde{\tau}: t_{i-1} = \tilde{t}_0 < \tilde{t}_1 < \ldots < \tilde{t}_{\tilde{N}_t} = t_{i}$ with $\Delta \tilde{t} = \frac{t_i - t_{i-1}}{\tilde{N}_t}$; \\
    	\For{$n=0,\ldots,\tilde{N}_t$}{
    		Interpolate $\bar{\mathbf{X}}^{\bar{N}_t}(\tilde{t}_n)$ from $\bar{\mathbf{X}}^{\bar{N}_t}$ using Brownian bridge interpolation;
    		Estimate $\Delta^\pm_a \bar{u} (\tilde{t}_n,\bar{\mathbf{X}}^{\bar{N}_t}(\tilde{t}_n))$ by linearly interpolating from $\left( \Delta^\pm_a \bar{u} \right)$ values computed in Algorithm~\ref{alg:hjb_solver} at grid $\tau$; \\
    		Using $\bar{\mathbf{X}}^{\bar{N}_t}(\tilde{t}_n)$ and $\Delta^\pm_a \bar{u} (\tilde{t}_n,\bar{\mathbf{X}}^{\bar{N}_t}(\tilde{t}_n))$, compute optimal controls $\boldsymbol\phi^*(\bar{t}_n,\bar{\mathbf{X}}^{\bar{N}_t}(\tilde{t}_n))$ by solving~\eqref{eqn:numerical_controls}; \\
    		Compute $\varsigma^{(m)}(n) = \left( \mathbbm{1}_{\{\phi_\text{out}(\tilde{t}_n,\bar{\mathbf{X}}^{\bar{N}_t}(\tilde{t}_n),\boldsymbol\phi(\tilde{t}_n,\bar{\mathbf{X}}^{\bar{N}_t}(\tilde{t}_n))) \geq \phi_\mathrm{th}\}} - \epsilon \right) \Delta \tilde{t}$; \\
    	}
    	$\varphi^{(m)}(i) = \sum_{n=0}^{\tilde{N}_t} \varsigma^{(m)}(n)$; \\
    }
    }
    $\left( \bar{\mathfrak{D}} \Theta (\lambda^{\ell}) \right)_i = \sum_{m=1}^{M_\mathrm{SG}} \varphi^{(m)}(i)$ for each $i=1,\ldots,\ell$; \\
    \textbf{Output: } $\bar{\mathfrak{D}} \Theta (\lambda^{\ell})$.
\end{algorithm}

\subsection[\appendixname~\thesubsection]{Initialization Algorithm for $\ell=1$}

\begin{algorithm}[H]
	\SetAlgoLined
	\KwIn{\( \Upsilon_1^1 = 1 \), \( \mathrm{TOL}_\text{init} \), $\beta_F,M_\mathrm{SG}$,$N_t,N_1,N_2,N_3,\bar{N}_t,M_\mathrm{SG}$,$\mathbf{X}(0)$}
	\KwOut{$\tilde{\Upsilon}_1^1$}
	Construct $\lambda^1(t)$ with $\Upsilon^1_1$ using~\eqref{eqn:finite_dim_multiplier}; \\
	Compute $\bar{\Theta}(\lambda^1) = \bar{u}(0,\mathbf{X}(0))$ by solving~\eqref{eqn:hjb_pde} using Algorithm~\ref{alg:hjb_solver} with $\lambda^1(t)$ and parameters $N_t,N_1,N_2,N_3$; \\
	Estimate $\bar{\mathfrak{D}}\Theta(\lambda^1)$  using Algorithm~\ref{alg:estimate_subgradient} with $\ell=1$,$\tilde{N}_t=\bar{N}_t$ and parameters $\bar{N}_t,M_\mathrm{SG}$; \\
	\eIf{$\bar{\mathfrak{D}}\Theta(\lambda^1) > 0$}{
	\While{$\abs{\bar{\mathfrak{D}}\Theta(\lambda^1)}>\mathrm{TOL}_\text{init} \epsilon$ or $\bar{\mathfrak{D}}\Theta(\lambda^1) > 0$}{
		$\Upsilon_1^1 \leftarrow \beta_F \Upsilon_1^1$; \\
		Construct $\lambda^1(t)$ with $\Upsilon^1_1$ using~\eqref{eqn:finite_dim_multiplier}; \\
		Compute $\bar{\Theta}(\lambda^1) = \bar{u}(0,\mathbf{X}(0))$ by solving~\eqref{eqn:hjb_pde} using Algorithm~\ref{alg:hjb_solver} with $\lambda^1(t)$ and parameters $N_t,N_1,N_2,N_3$; \\
	Estimate $\bar{\mathfrak{D}}\Theta(\lambda^1)$  using Algorithm~\ref{alg:estimate_subgradient} with $\ell=1$,$\tilde{N}_t=\bar{N}_t$ and parameters $\bar{N}_t,M_\mathrm{SG}$; \\
	}		
	}{
	\While{$\abs{\bar{\mathfrak{D}}\Theta(\lambda^1)}>\mathrm{TOL}_\text{init} \epsilon$ or $\bar{\mathfrak{D}}\Theta(\lambda^1) < 0$}{
		$\Upsilon_1^1 \leftarrow \frac{\Upsilon_1^1}{\beta_F}$; \\
		Construct $\lambda^1(t)$ with $\Upsilon^1_1$ using~\eqref{eqn:finite_dim_multiplier}; \\
		Compute $\bar{\Theta}(\lambda^1) = \bar{u}(0,\mathbf{X}(0))$ by solving~\eqref{eqn:hjb_pde} using Algorithm~\ref{alg:hjb_solver} with $\lambda^1(t)$ and parameters $N_t,N_1,N_2,N_3$; \\
	Estimate $\bar{\mathfrak{D}}\Theta(\lambda^1)$  using Algorithm~\ref{alg:estimate_subgradient} with $\ell=1$,$\tilde{N}_t=\bar{N}_t$ and parameters $\bar{N}_t,M_\mathrm{SG}$; \\
	}
	}
	$\tilde{\Upsilon}_1^1 = \Upsilon_1^1$; \\
	Construct $\lambda^1(t)$ with $\tilde{\Upsilon}^1_1$ using~\eqref{eqn:finite_dim_multiplier}; \\
\label{alg:initialization}
\caption{Initialization}
\end{algorithm}

\section[\appendixname~\thesection]{Details of Numerical Problem Described in Section~\ref{sec:results}}

\subsection[\appendixname~\thesection]{Model Parameters}
\label{appendix:model_parameters}

\begin{table}[H]
	\centering
	 \begin{tabular}{||c|c|c|c||} 
	 \hline
	 Parameter & Unit & Description & Value \\ [0.5ex] 
	 \hline\hline
	 $C_\mathrm{scal}$ & & Base station power-loss scaling factor & $7.84$ \\
	 $C_\text{offset}$ & Watt (W) & Base station offset power & $71.5$ \\
	 $\bar{P}_\mathrm{tx}$ & Watt (W) & Maximum base station transmission limit & $5 \times 10^3$ \\
	 $\mathbf{x}_\mathrm{BS}$ & & Base station location & $[0,0]$ \\
	 $\kappa$ &  & Path loss constant & $1$ \\
	 $\eta$ & & Path loss exponent & $2$ \\
	 $\text{SNR}_\mathrm{th}$ & Decibel (dB) & Signal-to-noise ratio threshold & $15$ \\
	 $\sigma_0$ & Watt (W) & Ambient transmission noise & $3.1623 \times 10^{-8}$ \\
	 $\bar{P}_R$ & Watt (W) & Maximum renewable power-production capacity & $10^4$ \\
	 $\bar{A}$ & Watt-hour (Wh) & Maximum battery charge capacity & $10^4$ \\
	 $\underline{P}_A$ & Watt (W) & Maximum battery charge capacity & $7.5 \times 10^3$ \\
	 $\bar{P}_A$ & Watt (W) & Maximum battery discharge capacity & $3 \times 10^4$ \\
	 $C_1$ & \texteuro/Wh  & Pollutant emission Coefficient 1 & $4 \times 10^{-4}$ \\
	 $C_2$ & \texteuro/$\text{W}^2$h & Pollutant emission Coefficient 2 & $10^{-4}$ \\
	 $P_K$ & \texteuro/Wh & Fictitious cost per unit battery charge & $0.0064$ \\
	 $w$ & & Pareto parameter & $0.5$ \\
%	 $\rho_Z$ & $\mathcal{N}(z_\mathrm{BS},\sigma_u^2 \mathbb{I}_2)$ \\
	% $\phi_\text{th}$ & $0.1$ \\
	% $\epsilon$ & $0.05$ \\
%	 $\alpha$ & $0.05$ \\
%	 $\theta_0$ & $1.6$ \\
%	 $B$ & $1$ \\
%	 $\sigma$ & $3$ \\ 
	% $A_\Omega$ & $4.9 \times 10^5$ \\
%	 $N_u$ & $[100,2000]$ \\
%	 $\sigma_u$ & $300$ \\
%	 $\pi$ & $0.5$ \\
	 \hline
	 \end{tabular}
	 \caption{Parameters and coefficients for modeling the cellular base station.}
	 \label{tab:model_parameters}
\end{table}

\subsection[\appendixname~\thesection]{Outage Proportion for Simple Cellular User Distributions}
\label{appendix:outage_prop}

Section~\ref{sec:base_station_model} demonstrates that the outage proportion for a given distribution of cellular users $\rho_z(t)$ at time $t$ can be written as follows:

\begin{equation*}
	\phi_\mathrm{out}(t) = \int_{\mathbb{R}^2} \mathbbm{1}_{\left\{ \norm{z - \mathbf{x}_{\mathrm{BS}}} > \left( \frac{P_{\mathrm{tx}}(t) \xi(t) \kappa}{\sigma_0 10^{\frac{\mathrm{SNR}_{\mathrm{th}}}{10}}} \right)^{\frac{1}{\eta}} \right\}} \rho_z(t) \dd z \cdot 
\end{equation*}

First, we consider the mobile user distribution to be a 2D uniform distribution in the geographical domain $\Pi$ served by the base station. The domain area is denoted as $A_\Pi$. Then,

\begin{equation}
	\label{eqn:uniform_outage}
	\phi_\mathrm{out}(t) = 1 - \frac{\pi}{A_\Pi} \left( \frac{P_{\mathrm{tx}}(t) \xi(t) \kappa}{\sigma_0 10^{\frac{\mathrm{SNR}_{\mathrm{th}}}{10}}} \right)^{\frac{2}{\eta}}, \quad \text{for } \rho_z(t) = \mathfrak{U}(\Pi),\quad \forall t \in [0,T] \cdot
\end{equation}

Next, we consider the case in which the mobile user distribution is a two-dimensional Gaussian distribution centered at the site of the base station $\mathbf{x}_\mathrm{BS}$ with a diagonal covariance matrix with variance $\sigma_u^2$ along both directions:

\begin{equation}
	\label{eqn:gaussian_outage}
	\phi_\mathrm{out}(t) = \exp \left( - \frac{1}{2 \sigma_u^2} \left( \frac{P_{\mathrm{tx}}(t) \xi(t) \kappa}{\sigma_0 10^{\frac{\mathrm{SNR}_{\mathrm{th}}}{10}}} \right)^{\frac{2}{\eta}} \right), \quad \text{for } \rho_z(t) = \mathcal{N}(\mathbf{x}_\mathrm{BS},\sigma_u^2 \mathbb{I}_2),\quad \forall t \in [0,T] ,
\end{equation}

where $\mathbb{I}_2$ denotes the 2D identity matrix. Furthermore, semi-analytical expressions for $\phi_\mathrm{out}(t)$ can be derived for other Gaussian distributions that are not centered on $\mathbf{x}_\mathrm{BS}$ or with general  covariance matrices.

\subsection[\appendixname~\thesection]{Analytical Solution of the Hamiltonian~\eqref{eqn:hamiltonian} for Simple Cellular User Distributions}
\label{appendix:hamiltonian_soln}

For each point $(t,\mathbf{x}) \in [0,T] \times \Gamma$, we must solve the optimization problem in~\eqref{eqn:hamiltonian}. For uniform or Gaussian cellular user distributions with closed-form expressions for $\phi_\mathrm{out}$ (Appendix~\ref{appendix:outage_prop}), this simplifies to solving the following constrained convex optimization problem.

\begin{empheq}[left=\empheqlbrace, right = \cdot]{equation}
\label{eqn:hamiltonian_analytical}
	\begin{alignedat}{2}
		&\min_{P_F,P_S,P_\mathrm{tx}} a_1 P_F + a_2 P_F^2 + b_1 P_S + c_1 P_\mathrm{tx} + c2 \exp{-c_3 P_\mathrm{tx}}  \\
		&\qquad + \lambda(t) \mathbbm{1}_{\{\phi_\text{out}(t,P_\mathrm{tx},\xi) > \phi_{\mathrm{th}}\}} \\
		\text{s.t. }& P_\mathrm{tx} - d_3 \leq 0 \\
		& -P_F + P_S + d_1 P_\mathrm{tx} + d_2 \leq 0 \\
		& P_F - P_S - d_1 P_\mathrm{tx} + d_4 \leq 0 \\
		& -P_F \leq 0 \\
		& -P_\mathrm{tx} \leq d_5 \\
		& -P_S \leq 0 \\
		& P_F - d_1 P_\mathrm{tx} - C_\text{offset} \leq 0
	\end{alignedat}
\end{empheq}
where
\begin{align*}
	a_1 &= \frac{1}{\bar{A}} \frac{\partial u}{\partial a} (t,\mathbf{x}) + w K_B(t) + (1-w) C_1, \\
	a_2 &= (1-w) C_2, \\
	b_1 &= - \frac{1}{\bar{A}} \frac{\partial u}{\partial a} (t,\mathbf{x}) - w K_B(t), \\
	c_1 &= \begin{cases} 
		- \frac{C_\text{scal} N_u(t)}{\bar{A}} \frac{\partial u}{\partial a} (t,\mathbf{x}) - \frac{w K_\text{net}(t) \pi \chi \kappa}{A_\Pi \sigma_0 10^{\frac{\text{SNR}_\mathrm{th}}{10}}}  & \text{,if } \rho_z(t) = \mathfrak{U}(\Pi), \\
		- \frac{C_\text{scal} N_u(t)}{\bar{A}} \frac{\partial u}{\partial a} (t,\mathbf{x}) & \text{,if } \rho_z(t) = \mathcal{N}(\mathbf{x}_\mathrm{BS},\sigma_u^2 \mathbb{I}_2).
		\end{cases}, \\
	c_2 &= \begin{cases} 
		0 & \text{,if } \rho_z(t) = \mathfrak{U}(\Pi), \\
		w K_\text{net}(t) N_u(t) & \text{,if } \rho_z(t) = \mathcal{N}(\mathbf{x}_\mathrm{BS},\sigma_u^2 \mathbb{I}_2).
		\end{cases}, \\
	c_3 &= \begin{cases} 
		0 & \text{,if } \rho_z(t) = \mathfrak{U}(\Pi), \\
		\frac{\chi \kappa}{2 \sigma_u^2 \sigma_0 10^{\frac{\text{SNR}_\mathrm{th}}{10}}} & \text{,if } \rho_z(t) = \mathcal{N}(\mathbf{x}_\mathrm{BS},\sigma_u^2 \mathbb{I}_2).
		\end{cases}, \\
	d_1 &= C_\text{scal} N_u(t), \\
	d_2 &= -\bar{P}_A (a) - r \bar{P}_R + C_\text{offset}, \\
	d_3 &= \frac{\bar{P}_\mathrm{tx}}{N_u(t)}, \\
	d_4 &= \underline{P}_A (a) + r \bar{P}_R - C_\text{offset}, \\
	d_5 &= \begin{cases} 
		\frac{(1-\phi_\mathrm{th}) A_\Pi \sigma_0 10^{\frac{\text{SNR}_\mathrm{th}}{10}}}{\pi \chi \kappa} & \text{,if } \rho_z(t) = \mathfrak{U}(\Pi), \\
		-\frac{2 \sigma_u^2 \sigma_0 10^{\frac{\text{SNR}_\mathrm{th}}{10}} \log{\phi_\mathrm{th}}}{\xi \kappa} & \text{,if } \rho_z(t) = \mathcal{N}(\mathbf{x}_\mathrm{BS},\sigma_u^2 \mathbb{I}_2).
		\end{cases} \cdot \\
\end{align*}

The discontinuous indicator function in the objective function poses a difficulty while deriving the analytical solution to~\eqref{eqn:hamiltonian_analytical}. Decomposing~\eqref{eqn:hamiltonian_analytical} simplifies this into two constrained convex optimization subproblems, which are

\begin{empheq}[left=\empheqlbrace, right = \cdot]{equation}
	\begin{alignedat}{2}
		\mathcal{H}_1 = &\min_{P_F,P_S,P_\mathrm{tx}} a_1 P_F + a_2 P_F^2 + b_1 P_S + c_1 P_\mathrm{tx} + c2 \exp{-c_3 P_\mathrm{tx}} \\
		\text{s.t. }& -P_\mathrm{tx} \leq -d_5 \\
		& P_\mathrm{tx} - d_3 \leq 0 \\
		& -P_F + P_S + d_1 P_\mathrm{tx} + d_2 \leq 0 \\
		& P_F - P_S - d_1 P_\mathrm{tx} + d_4 \leq 0 \\
		& -P_F \leq 0 \\
		& -P_S \leq 0 \\
		& P_F - d_1 P_\mathrm{tx} - C_\text{offset} \leq 0 \cdot \nonumber 
	\end{alignedat}
\end{empheq}

and

\begin{empheq}[left=\empheqlbrace, right = \cdot]{equation}
	\begin{alignedat}{2}
		\mathcal{H}_2 = &\min_{P_F,P_S,P_\mathrm{tx}} a_1 P_F + a_2 P_F^2 + b_1 P_S + c_1 P_\mathrm{tx} + c2 \exp{-c_3 P_\mathrm{tx}} + \lambda(t) \\
		\text{s.t. }& P_\mathrm{tx} \leq d_5 \\
		& -P_F + P_S + d_1 P_\mathrm{tx} + d_2 \leq 0 \\
		& P_F - P_S - d_1 P_\mathrm{tx} + d_4 \leq 0 \\
		& -P_F \leq 0 \\
		& -P_\mathrm{tx} \leq 0 \\
		& -P_S \leq 0 \\
		& P_F - d_1 P_\mathrm{tx} - C_\text{offset} \leq 0 \cdot \nonumber 
	\end{alignedat}
\end{empheq}

The optimal controls for~\eqref{eqn:hamiltonian_analytical} are given as

\begin{equation*}
	\argmin_{P_F,P_S,P_\mathrm{tx}} \left( \mathcal{H}_1,\mathcal{H}_2 \right) \cdot
\end{equation*}

Analytical solutions to these subproblems can be derived by solving the corresponding Karush--Kuhn--Tucker conditions for constrained optimization problems~\citep{boyd2004convex}. 

\subsection[\appendixname~\thesection]{Simulation Parameters in Section~\ref{sec:results}}
\label{appendix:sim_parameters}

\begin{table}[H]
	\centering
	 \begin{tabular}{||c|c|c|c||} 
	 \hline
	 Parameter & Description & Value \\ [0.5ex] 
	 \hline\hline
	 $\phi_\text{th}$ & Mobile user outage proportion threshold & $10^{-3}$ \\
     $\epsilon$ & Confidence level of violating the constraint in~\ref{constraint:prob} & $0.1$ \\
	 $A_0$ & Normalized battery charge level at $t=0$ & $0.5$ \\
	 $\mu_0^\xi$ & Initial distribution of the wireless fading channel & $\gamma(3,1)$ \\
%	 $\mu_0^r$ & Initial distribution of the normalized wind power & $\mathcal{N} \left( p(-\varsigma) + \dot{p}(-\varsigma) \varsigma, 2 \alpha \theta_0 p(-\varsigma) (1-p(-\varsigma) ) \varsigma \right)$ \\
	 $N_1$ & Discretization of~\eqref{eqn:hjb_pde} in the $a$ domain & $10$ \\
	 $N_2$ & Discretization of~\eqref{eqn:hjb_pde} in the $r$ domain & $10$ \\
	 $N_3$ & Discretization of~\eqref{eqn:hjb_pde} in the $\chi$ domain & $10$ \\
	 $N_t$ & Discretization of~\eqref{eqn:hjb_pde} in the $t$ domain & $800$ \\
	 $\mathrm{TOL}$ & Prescribed relative tolerance for Algorithm~\ref{alg:dual_optimization} & $0.1$ \\
	 $\mathrm{TOL}_\text{init}$ & Prescribed relative tolerance for Algorithm~\ref{alg:initialization} & $1$ \\
	 \texttt{max-iter} & Prescribed maximum iterations in Algorithm~\ref{alg:dual_optimization} & $50$ \\
	 $\bar{N}_\text{iter}$ & Initial number of SSM iterations with a constant step-size & $10$ \\
	 $N_\text{iter}$ & Prescribed number of LMBM~\citep{Karmitsa:2007aa} iterations & $50$ \\
	 $\beta_F$ & Factor of increase/decrease in Algorithm~\ref{alg:initialization}  & $5$ \\
	 $M_\mathrm{SG}$ & Number of sample paths in Algorithm~\ref{alg:estimate_subgradient} & $10^3$ \\
	 $\bar{N}_t$ & Time discretization parameter in Algorithm~\ref{alg:estimate_subgradient}  & $64$ \\
	 $\tilde{N}_t$ & Time discretization parameter in Algorithm~\ref{alg:estimate_subgradient} & $\frac{\bar{N}_t}{\ell}$ \\
	 \hline
	 \end{tabular}
	 \caption{Simulation parameters to run Algorithm~\ref{alg:dual_optimization} and produce numerical results as described in Section~\ref{sec:results}.}
	 \label{tab:sim_parameters}
\end{table}

\begin{table}[H]
	\centering
	 \begin{tabular}{||c|c|c||} 
	 \hline
	 Parameter & Description & Value \\ [0.5ex] 
	 \hline\hline
	 \texttt{RPAR(1)} & Tolerance for changes in the function value & $0.1$ \\
     \texttt{RPAR(2)} & Second tolerance for changes in the function value & $-1$ (ignored) \\
	 \texttt{RPAR(3)} & Minimum acceptable function value & $0$ \\
	 \texttt{RPAR(4)} & Tolerance for the first termination parameter & $10^{-2}$ \\
	 \texttt{RPAR(5)} & Tolerance for the second termination parameter & $10^{-2}$ \\
	 \texttt{RPAR(6)} & Distance measure parameter & $0.5$ \\
	 \texttt{RPAR(7)} & Line search parameter & $0.2$ \\
	 \texttt{RPAR(8)} & Maximum step size & $10$ \\
	 \texttt{IPAR(1)} & Exponent for distance measure & $2$ \\
	 \texttt{IPAR(2)} & Maximum iterations & $50$ \\
	 \texttt{IPAR(3)} & Maximum function evaluations & $100$ \\
	 \texttt{IPAR(4)} & Maximum iterations with changes of function values smaller than \texttt{RPAR(1)} & $5$ \\
	 \texttt{IPAR(5)} & Printout specification & $-1$ \\
	 \texttt{IPAR(6)} & Selection of method & $0$ (LMBM) \\
	 \texttt{IPAR(7)} & Selection of scaling strategy & $0$ \\
	 \hline
	 \end{tabular}
	 \caption{Parameters to run the LMBM routine~\citep{Karmitsa:2007aa} in Algorithm~\ref{alg:dual_optimization} to produce numerical results described in Section~\ref{sec:results}. See~\citep{Karmitsa:2007aa} for more details.}
	 \label{tab:lmbm_parameters}
\end{table}

\subsection[\appendixname~\thesection]{Sensitivity analysis settings in Section~\ref{sec:results}}
\label{appendix:sensitivity_analysis}

\begin{table}[H]
	\centering
	 \begin{tabular}{||c|c||} 
	 \hline
	 Parameter & Distribution \\ [0.5ex] 
	 \hline\hline
	 $\bar{P}_\mathrm{tx}$ & $\log_{10}(\bar{P}_\mathrm{tx}) \sim \mathfrak{U}[3,4]$ \\
     $\sigma_0$ & $\sigma_0 \sim \mathcal{N} \left( 3.1623 \times 10^{-8}, 10^{-16} \right)$ \\
	 $\bar{P}_R$ & $\log_{10}(\bar{P}_R) \sim \mathfrak{U}[2,4]$ \\
	 $\bar{A}$ & $\log_{10}(\bar{A}) \sim \mathfrak{U}[3,4]$ \\
	 $w$ & $w \sim \mathfrak{U}[0,1]$ \\
	 $\phi_\mathrm{th}$ & $\log_{10}(\phi_\mathrm{th}) \sim \mathfrak{U}[-2,-4]$ \\
	 $A_0$ & $A_0 \sim \mathfrak{U}[0,1]$ \\
	 $K_\text{net}$ & $\log_{10}(K_\text{net}) \sim \mathfrak{U}[-1,1]$ \\
	 \hline
	 \end{tabular}
	 \caption{Distribution of model and algorithmic parameters used to run randomized multi-scenario simulations of Algorithm~\ref{alg:dual_optimization}.}
	 \label{tab:scenarios}
\end{table}

\bibliography{References}
\bibliographystyle{plainnat}
\setcitestyle{authoryear,open={[},close={]}}
\end{document}